\documentclass[reqno]{amsart}
\usepackage{amssymb}
\usepackage{amsmath}
\usepackage[mathscr]{euscript}
\usepackage{times}

\usepackage{color}

\makeatletter
\@addtoreset{equation}{section}
\makeatother

\newtheorem{theorem}{Theorem}[section]
\newtheorem{lemma}[theorem]{Lemma}
\newtheorem{proposition}[theorem]{Proposition}
\newtheorem{corollary}[theorem]{Corollary}

\newtheorem{definition}[theorem]{Definition}
\newtheorem{remark}[theorem]{Remark}

\newcommand{\mc}[1]{{\mathcal #1}}
\newcommand{\mf}[1]{{\mathfrak #1}}

\newcommand{\bb}[1]{{\mathbb #1}}
\newcommand{\bs}[1]{{\boldsymbol #1}}
\newcommand{\ms}[1]{{\mathscr #1}}

\newcounter{as}[section]

\renewcommand{\a}{\alpha}

\newcommand{\e}{\varepsilon}
\newcommand{\de}{\delta}

\newcommand{\ga}{\gamma}

\newcommand{\la}{\lambda}

\renewcommand{\P}{\mathbb{P}}

\newcommand{\De}{\Delta}
\newcommand{\Ga}{\Gamma}

\newcommand{\lan}{\langle}
\newcommand{\ran}{\rangle}

\newcommand{\R}{\mathbb{R}}

\newcommand{\N}{\mathbb{N}}
\newcommand{\Z}{\mathbb{Z}}

\newcommand{\T}{\mathbb{T}}
\newcommand{\M}{\mathcal{M}}

\newcommand{\cP}{\mathcal{P}}

\newcommand{\vte}{\varrho}

\newcommand{\<}{\langle}
\renewcommand{\>}{\rangle}

\begin{document}

\title{Static large deviations for a reaction-diffusion model}
\author{J. Farf\'an, C. Landim, K. Tsunoda}

\address{\noindent PUCP, Av. Universitaria cdra. 18,
  San Miguel, Ap. 1761, Lima 100, Per\'u. 
\newline e-mail: \rm
  \texttt{jfarfanv@pucp.edu.pe} }

\address{\noindent IMPA, Estrada Dona Castorina 110, CEP 22460 Rio de
  Janeiro, Brasil and CNRS UMR 6085, Universit\'e de Rouen, Avenue de
  l'Universit\'e, BP.12, Technop\^ole du Madril\-let, F76801
  Saint-\'Etienne-du-Rouvray, France.  \newline e-mail: \rm
  \texttt{landim@impa.br} }

\address{\noindent Department of Mathematics, Osaka university,
1-1, Machikaneyama-cho, Toyonaka, Osaka, 560-0043, Japan.
\newline e-mail:  \rm \texttt{k-tsunoda@math.sci.osaka-u.ac.jp}}



\begin{abstract}
  We consider the superposition of a symmetric simple exclusion
  dynamics, speeded-up in time, with a spin-flip dynamics in a
  one-dimensional interval with periodic boundary conditions.  We
  prove the large deviations principle for the empirical measure under
  the stationary state. We deduce from this result that the stationary
  state is concentrated on the stationary solutions of the
  hydrodynamic equation which are stable.
\end{abstract}

\maketitle

\section{Introduction}
\label{vsec00}

Nonequilibrium thermodynamics has aroused a lot of interest in the
last decades. Since the beginning of the 2000's, much attention has
been devoted to the investigation of nonequilibrium stationary states
which describe a steady flow through a system, \cite{d07,bdgjl15} and
references therein.

Over the last years, a general approach to examine nonequilibrium
stationary states, called the Macroscopic Fluctuation Theory, has been
developed based on a dynamical large deviations principle for the
empirical current \cite{bdgjl02,bd06,bdgjl06}. Among the major
achievements of the MFT was the deduction of a time-independent
variational formula for the quasi-potential, the functional obtained
by minimizing the dynamical large deviations rate functional over all
trajectories which start from the stationary density profile and
produces a fixed fluctuation \cite{dls,bdgjl03}, and the proof that
the quasi-potential is G\^ateaux differentiable at some density
profile if and only if the time-dependent variational formula which
defines the quasi-potential has a unique minimizer \cite{bdgjl11}.

At the same time, adapting to the infinite-dimensional setting the
strategy proposed by Freidlin and Wentzell \cite{fw1} for stochastic
perturbations of finite-dimensional dynamical systems, Bodineau and
Giacomin \cite{bg} and Farf\'an \cite{f} proved a large deviations
principle for the empirical measure under the nonequilibrium
stationary state for conservative dynamics in contact with reservoirs
in which the large deviations rate functional is given by the
quasi-potential.

We consider in this article the stochastic evolution obtained by
superposing a speeded-up symmetric simple exclusion process with a
spin-flip dynamics on a one-dimensional interval with periodic
boundary conditions.  The hydrodynamic equation induced by the
microscopic dynamics, the partial differential equation which
describes the macroscopic evolution of the density, is given by a
reaction-diffusion equation of type
\begin{equation}
\label{000}
\partial_t \rho \,=\,  (1/2)\, \De \rho  + B(\rho) - D(\rho) \;,
\end{equation}
where $\Delta$ represents the Laplacian and where $B$ and $D$ are
non-negative polynomials.

We investigate the static large deviations of the empirical measure
under the stationary state. In contrast with the previous dynamics
\cite{bg,f}, in which the hydrodynamic equation has a unique
stationary solution which is a global attractor of the dynamical
system generated by the PDE, in reaction-diffusion models, for
appropriate choices of the functions $B$, $D$, the hydrodynamic
equation possesses more than one stationary solution.

The existence of multiple stationary solutions to the PDE \eqref{000}
raises new problems and new questions. For instance, the conjecture
that the stationary state does not put mass on the unstable solutions
of the hydrodynamic equation \eqref{000}. For reaction-diffusion
models it has only been proved in \cite{lt} that the stationary state
is concentrated on the set of classical solutions to the semilinear
elliptic equation
\begin{equation}
\label{eqi1}
(1/2) \, \De \rho  + B(\rho) - D(\rho) \;=\; 0 \;.
\end{equation}

The main result of this article, Theorem \ref{sldp}, establishes a
large deviations principle for the empirical measure under the
stationary state. The quasi-potential, the rate functional of this
large deviations principle, is represented through a time-dependent
variational problem involving the dynamical large deviations rate
functional. The value of the quasi-potential at a measure $\vte$ is
given in terms of the infimum of the dynamical large deviations rate
functional over all trajectories which start from a stationary
solution of the hydrodynamic equation and end at $\vte$.

A consequence of this result is that the stationary measure is
concentrated at the stable, stationary solutions of the hydrodynamic
equation. This is the content of Theorem \ref{stablesol}, the second main
result of the article.

The proof of Theorem \ref{stablesol} is based on two properties of the
reaction-diffusion equation \eqref{000}. First, it is known from
\cite{cm89} that all solutions of \eqref{000} converge to solutions of
the semilinear elliptic equation \eqref{eqi1}. In particular, there
are no time-periodic solutions. Second, we assume that equation
\eqref{eqi1} has only a finite number of solutions, modulo
translations. This property holds when the polynomial $F(\rho) =
B(\rho) - D(\rho)$ satisfies the hypotheses of Lemma \ref{rd-l1}.
Theorem \ref{stablesol} further requires a characterization of the unstable
stationary solutions of \eqref{000}. Under the conditions of Lemma
\ref{rd-l2} on $F$, this set consists of all non-constant solutions and
all constant solutions associated to local maxima of the potential $V$,
where $V'(\rho) = - F(\rho)$.

\smallskip

We conclude this introduction with some comments.  This stochastic
dynamics has been introduced by De Masi, Ferrari and Lebowitz in
\cite{dfl}. The authors proved the hydrodynamic limit of the system by
duality arguments and the fluctuations of the density field.  The
dynamical large deviations principle for the empirical density
starting from a product measure appeared in \cite{jlv}, following the
ideas presented in \cite{kov}. Bodineau and Lagouge in \cite{bl2,bl}
proved the dynamical large deviations principle for the empirical
current, while two of the authors of this article extended in
\cite{lt} the dynamical large deviations principle to the case in
which the process starts from a deterministic configuration.

The stationary states of the symmetric simple exclusion process are
the Bernoulli product measures. The introduction of the spin-flip
dynamics creates long range correlations. Although the local
distribution of particles remains very close to a Bernoulli product
measure due to the speeding-up of the exclusion dynamics, the long
range correlations affect substantially the macroscopic behavior of the
system. The purpose of this article is to study this effect at the
level of the large deviations.

There is a huge literature on large deviations for
reaction-diffusion equations perturbed by Gaussian or L\'evy noise in
finite and infinite dimensions after the seminal paper by Faris and
Jona-Lasinio \cite{fj}. We refer to the recent books \cite{dz,dhi13}
for references on the subject. The noise created by the microscopic
spin-flip dynamics considered in this article is of a different
nature. This is reflected in the dynamical large deviations rate
functions in which singular exponential terms appear. This is one of
the sources of technical problems faced in order to prove the
regularity conditions of the dynamical rate functional needed to
derive the static large deviations principle.

We leave to the end of the next section technical comments and remarks
on the proofs and on the assumptions, and we mention here some open
problems for future research.  It would be interesting to extend to
this model the results described at the beginning of this introduction
which were obtained from the MFT for one-dimensional conservative
interacting particle systems in contact with reservoirs: an
alternative time-independent variational formula for the
quasi-potential, and a description of the optimal trajectory which
solves the time-dependent variational formula defining the
quasi-potential. This has been done in \cite{gjlv} in the case where
the reaction-diffusion model is reversible, but it remains an open
problem in the non-reversible setting. In this general situation the
only available information is an expansion of the quasi-potential
around a constant stable stationary point obtained by Basile and
Jona-Lasinio \cite{bj}. A description of the metastable behavior of
the reaction-diffusion model when the difference $B(\rho) - D(\rho)$
forms a double-well potential is also a challenging open problem.

\section{Notation and Results}
\label{sec2}

Throughout this article, we use the following notation.  $\N_{0}$
stands for the set $\{0,1,\cdots\}$.  For a function $f: X\to \bb R$,
defined on some set $X$, let $\|f\|_\infty = \sup_{x\in X}{|f(x)|}$.
We sometimes denote the interval $[0,\infty)$ by $\R_{+}$.

\subsection{Reaction-diffusion model}
Let $\T_{N}=\Z/N\Z$, $N\ge 1$, be the one-dimensional discrete torus
with $N$ points.  Denote by $X_{N}$ the set $\{0,1\}^{\T_{N}}$ and by
$\eta$ the elements of $X_{N}$, called configurations.  For each
$x\in\T_{N}$, $\eta(x)$ represents the occupation variable at site $x$
so that $\eta(x)=1$ if the site $x$ is occupied for the configuration
$\eta$, and $\eta(x)=0$ if the site is vacant.  For each $x \not
=y\in\T_{N}$, denote by $\eta^{x,y}$, resp. by $\eta^{x}$, the
configuration obtained from $\eta$ by exchanging the occupation
variables $\eta(x)$ and $\eta(y)$, resp. by flipping the occupation
variable $\eta(x)$:
\begin{align*}
\eta^{x,y}(z) \;=\; \begin{cases}
     \eta(y)  & \text{if $z=x$}\;, \\
     \eta(x)  & \text{if $z=y$}\;, \\
     \eta(z)    & \text{otherwise}\;,
\end{cases}
\quad
\eta^x(z) \;=\; \begin{cases}
     \eta(z)  & \text{if $z\neq x$}\;, \\
     1-\eta(z)  & \text{if $z=x$}\;. 
\end{cases}
\end{align*}

Consider the superposition of the speeded-up symmetric simple
exclusion process with a spin-flip dynamics. The generator of this
$X_{N}$-valued, continuous-time Markov chain acts on functions
$f:X_{N}\to\R$ as
\begin{align*}
\mathcal{L}_Nf \;=\; N^2 \mathcal{L}_Kf + \mathcal{L}_Gf\;,
\end{align*}
where $\mathcal{L}_K$ is the generator of a symmetric simple exclusion
process (Kawasaki dynamics),
\begin{align*}
(\mathcal{L}_Kf)(\eta) \;=\; (1/2) \sum_{x\in\T_N} [f(\eta^{x,x+1}) - f(\eta)]\;,
\end{align*}
and $\mathcal{L}_G$ is the generator of a spin-flip dynamics
(Glauber dynamics),
\begin{align*}
(\mathcal{L}_Gf)(\eta) \;=\; \sum_{x\in\T_N} c(\tau_x \eta)[f(\eta^x) - f(\eta)]\;.
\end{align*}
In the last formula, $c(\eta)$ represents a strictly positive,
cylinder function, that is, a function $c: \{0,1\}^\Z \to \bb R_+$ which
depends only on a finite number of coordinates $\eta(y)$.
For a sufficiently large $N$, $c$ can be regarded as a function on $X_N$.
$\{\tau_x : x\in \bb Z\}$ represents the group of translations defined
by $(\tau_x \eta)(y) = \eta(x+y)$, $y\in \bb T_N$, where the sum is
carried modulo $N$.

Note that the Kawasaki dynamics has been speeded-up by a factor
$N^{2}$, which corresponds to the diffusive scaling.  Setting the jump
rates of the Glauber part to be $0$, we retrieve the symmetric simple
exclusion dynamics speeded up by $N^{2}$, whose static large deviation
principle has been derived with several different boundary conditions
in \cite{bdgjl03,bg,f,bl2,bl}.

Fix a topological space $X$. Let $D(I,X)$, $I=[0,T]$, $T>0$, or
$I = \R_+$, be the space of right-continuous trajectories from $I$
to $X$ with left-limits, endowed with the Skorohod topology.  Let
$\{\eta_t^N:N\ge1\}$ be the continuous-time Markov process on $X_{N}$
whose generator is given by $\mathcal{L}_{N}$. For a probability
measure $\nu$ on $X_{N}$, denote by $\P_{\nu}$ the probability measure
on $D(\R_{+},X_{N})$ induced by the process $\eta^{N}_{t}$ starting
from $\nu$. The expectation with respect to $\P_{\nu}$ is represented
by $\bb E_{\nu}$. Denote by $\P_{\eta}$ the measure $\P_{\nu}$ when
the probability measure $\nu$ is the Dirac measure concentrated on the
configuration $\eta$. Analogously, $\bb E_{\eta}$ stands for the
expectation with respect to $\bb P_{\eta}$.

\subsection{Hydrodynamics}

Let $\T$ be the one-dimensional continuous torus $\T=\R/\Z=[0,1)$.
Denote by $L^p(\T)$, $p\ge1$, the space of all real $p$-th integrable
functions $G:\T\to\R$ with respect to the Lebesgue measure $d\theta$:
$\int_\T |G(\theta)|^p d\theta <\infty$. The corresponding norm is denoted by
$\|\cdot\|_{p}$:
\begin{equation*}
\|G\|_{p}^{p}\;:=\;\int_\T |G(\theta)|^p d\theta\;.
\end{equation*}
In particular, $L^2(\T)$ is a Hilbert space equipped with the inner
product
\begin{align*}
\lan G, H \ran \;=\; \int_\T G(\theta)H(\theta) d\theta\;.
\end{align*}
For a function $G$ in $L^{2}(\T)$, we also denote by $\lan G \ran$ the
integral of $G$ with respect to the Lebesgue measure: $\lan G \ran
:=\int_{\T} G(\theta) d\theta$.

Let $\M_{+}=\M_{+}(\T)$ be the space of all nonnegative measures on $\T$
with total mass bounded by $1$, endowed with the weak topology.  For a
measure $\vte$ in $\mathcal{M}_+$ and a continuous function
$G:\T\to\R$, denote by $\lan \vte, G \ran$ the integral of $G$ with
respect to $\vte$:
\begin{equation*}
\lan \vte, G \ran \;=\; \int_\T G(\theta) \vte(d\theta) \;.
\end{equation*}
The space $\M_+$ is metrizable. Indeed, if $e_{0}(\theta)=1$, $e_{k}(\theta)=
\sqrt 2\cos(2\pi k \theta)$ and $e_{-k}(\theta)= \sqrt 2\sin(2\pi k\theta)$,
$k\in\N$, one can define a distance $d$ on $\M_{+}$ by
\begin{equation}
\label{01}
d(\vte_{1},\vte_{2}) \;:=\; \sum_{k\in\Z} \dfrac{1}{2^{|k|}}
\, |\lan\vte_{1},e_{k}\ran - \lan\vte_{2}, e_{k}\ran| \; ,
\end{equation}
and one can check that the topology induced by this distance
corresponds to the weak topology.  

Note that $\M_{+}$ is compact under the weak topology, and that, by
Schwarz inequality, for all density profiles $\gamma$, $\gamma':\bb
T\to [0,1]$,
\begin{equation}
\label{05}
d(\gamma, \gamma') \;\le\; 3 \Vert \gamma - \gamma'\Vert_2\;.
\end{equation}
In the previous formula we abuse of notation by writing
$d(\gamma, \gamma')$ for $d(\gamma(\theta) d\theta, \gamma'(\theta) d\theta)$.

Denote by $C^m(\T)$, $m$ in $\N_{0}\cup\{\infty\}$, the set of all
real functions on $\T$ which are $m$ times differentiable and whose
$m$-th derivative is continuous.  Given a function $G$ in $C^2(\T)$,
we shall denote by $\nabla G$ and $\De G$ the first and second
derivatives of $G$, respectively.

Let $\nu_\rho=\nu^N_{\rho}$, $0\le\rho\le 1$, be the Bernoulli product
measure on $X_N$ with the density $\rho$.  Define the continuous
functions $B,D:[0,1]\to\R$ by
\begin{align*}
B(\rho)\;=\;\int [1-\eta(0)]\, c(\eta)\, d\nu_\rho\;, \quad
D(\rho)\;=\;\int \eta(0)\, c(\eta)\, d\nu_\rho \;.
\end{align*}

Let $\pi^N: X_N\to\M_+$ be the function which associates to a
configuration $\eta$ the positive measure obtained by assigning mass
$N^{-1}$ to each particle of $\eta$,
\begin{align*}
\pi^N(\eta)\;=\;\frac{1}{N}\sum_{x\in \T_N}\eta(x)\delta_{x/N}\;,
\end{align*}
where $\delta_\theta$ stands for the Dirac measure which has a point mass at
$\theta\in\T$. Let $\pi^{N}_{t} = \pi^N (\eta^N_t)$, $t\ge 0$.  The next
result was proved by De Masi, Ferrari and Lebowitz in \cite{dfl} for
the first time. We refer to \cite{dfl,jlv,kl} for its proof.

\begin{theorem}
\label{hydrodynamics}
Fix a measurable function $\ga:\T\to[0,1]$.  Let $\nu_N$ be a sequence
of probability measures on $X_N$ associated to $\ga$, in the sense
that
\begin{equation}
\label{24}
\lim_{N\to\infty}
\nu_N \Big(|\lan\pi^{N},G\ran -\int_{\T} G(\theta)\ga(\theta)d\theta|>\delta\Big)\;=\;0\;,
\end{equation}
for every $\delta>0$ and every continuous function $G:\T\to\R$.  Then,
for every $t\ge 0$, every $\delta>0$ and every continuous function
$G:\T\to\R$, 
\begin{equation*}
\lim_{N\to\infty}
\mathbb{P}_{\nu_N}
\Big(|\lan\pi_{t}^{N},G\ran -\int_{\T} G(\theta)\rho(t,\theta)d\theta|>\delta\Big)\;=\;0\;,
\end{equation*}
where $\rho:[0,\infty)\times\T\to[0,1]$ is the unique weak solution of
the Cauchy problem
\begin{equation}
\label{rdeq}
\begin{cases}
     \partial_t\rho \;=\; (1/2) \De \rho + F(\rho) \ \text{ on }\ \T\;,\\
     \rho(0,\cdot)\;=\;\ga(\cdot)\;,
\end{cases}
\end{equation}
where $F(\rho)=B(\rho)-D(\rho)$.
\end{theorem}

The definition, the existence and the uniqueness of weak solutions of
the Cauchy problem \eqref{rdeq} are discussed in Section \ref{sec3}.

\subsection{The reaction-diffusion equation}

We present in this subsection the results on the reaction-diffusion
equation \eqref{rdeq} needed in this section. Let $S$ be the set of
all classical solutions of the semilinear elliptic equation:
\begin{equation}
\label{seeq}
(1/2)\De\rho + F(\rho) \;=\; 0 \ \text{ on }\ \T\;.
\end{equation}
Classical solution means a $[0,1]$-valued function $\rho$ in $C^{2}(\T)$ which
satisfies the equation \eqref{seeq} for any $\theta\in\T$.  We also denote
by $\M_{\rm sol}$ the set of all absolutely continuous measures whose
density is a classical solution of \eqref{seeq}:
\begin{equation*}
\M_{\rm sol}\;:=\;\{\bar\vte\in\M_+:\bar\vte(d\theta)=\bar\rho(\theta)d\theta,\  
\bar\rho\in S \}\;.
\end{equation*}
Next lemma is Theorem D of \cite{cm89}.

\begin{lemma}
\label{rd-l0}
Let $\rho:[0,\infty)\times\T\to[0,1]$ be the unique weak solution of
the Cauchy problem \eqref{rdeq}.  Then, there exists a density profile
$\rho_{\infty}$ in $S$ such that $\rho_{t}$ converges to
$\rho_{\infty}$ as $t\to\infty$ in $C^{2}(\T)$.
\end{lemma}

This result excludes the existence of time-periodic solutions of equation
\eqref{rdeq}, a phenomenon which occurs if the function $F$ is allowed
to depend on $\nabla\rho$ as well (cf. \cite{frw} and references
therein).

We turn to the description of the set $S$.  Denote by $\mf R$ the set
of roots of $F$ in $[0,1]$.  It is clear that for all $r\in \mf R$,
the constant function $\rho: \bb T \to \bb R$ given by $\rho(\theta) = r$,
$\theta\in\bb T$, is an element of $S$. There might be also non-constant
periodic solutions.

Let $V:[0,1] \to \bb R$ be a potential such that $F(\rho) = -
V'(\rho)$. If the polynomial $F$ has degree $1$, as $V'(0)<0<V'(1)$,
equation \eqref{seeq} has a unique solution, which is a global
attractor for the dynamical system induced by the reaction-diffusion
equation \eqref{rdeq}, and given by $\rho(\theta) = r$, where $r$ is the
unique root of $F$.

Assume that the degree of $F$ is larger than or equal to $2$.  Denote
by $m_1, \dots, m_{\mf n}$ the local minima of $V$ in $[0,1]$ and by
$M_1, \dots, M_{\mf m}$ the local maxima in this interval. Since
$V'(0)<0<V'(1)$, $\mf n =\mf m+1 \ge 1$ and $m_1<M_1< \dots < M_{\mf
  n-1}<m_\mf n$.

Denote by $\sim$ the equivalence relation in $C^2(\bb T)$ defined by
$\rho \sim \rho'$ if there exists $\theta'\in \bb T$ such that $\rho'
(\theta) = \rho(\theta+\theta')$ for all $\theta\in\bb T$.  Of course,
if $\rho$ is a periodic solution and $\rho'\sim\rho$, then $\rho'$ is
also a solution.

\begin{lemma}
\label{rd-l1}
Suppose that all zeros of $F$ are real and that all critical points of
$V$ are local minima or local maxima.  Then, the elliptic equation
\eqref{seeq} with periodic boundary conditions has at most a finite
number of solutions, modulo the equivalence relation introduced above.
\end{lemma}

Since we could not find the previous result explicitly stated in the literature,
we sketch the proof of this result. The terminology employed can be
found in \cite{Sch}. By Proposition 1.5.2 in \cite{Sch}, $F$ is an
$A-B$ function on all intervals $(m_1,M_1)$, $\dots, (M_{\mf
  n-1},m_{\mf n})$. The diagram of the Hamiltonian system $\dot p =
q$, $\dot q = V'(p)$ shows that the periodic solutions of \eqref{seeq}
are bounded below and above by two consecutive minima of the potential
$V$.

Solutions of \eqref{seeq} with periodic boundary conditions can be
mapped to solutions of \eqref{seeq} with Dirichlet boundary
conditions. Indeed, fix two consecutive minima $m_j$, $m_{j+1}$ of $V$, and a
solution $\rho$ of \eqref{seeq} taking values in $[m_j, m_{j+1}]$. Let
$F_j(r) = F(M_j+r)$, so that $F_j(0)=0$ because $M_j$ is a local
maximum of $V$. Note that $F_j$ is an $A-B$ function on $(m_j -M_j, 0)
\cup (0, m_{j+1}-M_j)$. Let $\theta_0 = \min \{\theta\ge 0 :
\rho(\theta) = M_j\}$. Define $\phi:[0,1]\to [m_j-M_j, m_{j+1}-M_j]$
by $\phi(\theta) = \rho(\theta+\theta_0)-M_j$. It is clear that $\phi$
is a solution of $(1/2) \Delta v + F_j(v)=0$ with Dirichlet boundary
conditions.

Since $F_j$ is an $A-B$ function on the intervals $(m_j- M_{j},0)$ and
$(0, m_{j+1}-M_j)$, by Propositions 3.1.3, 3.1.4 and Theorem 3.1.9 in
\cite{Sch}, the time-map of the equation \eqref{seeq} with Dirichlet
boundary conditions is strictly convex and converges to $+\infty$ at
the boundary. In particular, for each branch there exist at most two
distinct solutions if $V''(M_j)=0$ and at most one solution if
$V''(M_j)<0$.  Since there is a finite number of branches whose
time-map takes value less than or equal to $1$, there is a finite
number of distinct solutions of \eqref{seeq} with Dirichlet boundary
conditions. As all solutions with periodic boundary conditions can be
mapped to solutions with Dirichlet boundary conditions, the lemma is
proved.

We turn to the heteroclinic orbits of \eqref{rdeq}. A complete
description has been obtained in \cite{frw}. We state here a partial
result which fulfills our needs. It asserts that all non-constant
stationary solutions are unstable, as well as all constant solutions
associated to local maxima of $V$. 

Fix two stationary solutions $\phi\neq\psi$ of \eqref{rdeq}. A
trajectory $\rho(t,\cdot)$, $t\in\bb R$, is called a heteroclinic
orbit from $\phi$ to $\psi$ if $\lim_{t\to-\infty} = \phi$,
$\lim_{t\to+\infty} = \psi$ and if $\rho$ solves \eqref{rdeq} for
every $t\in \bb R$. Convergences are meant in $C^1(\bb T)$.  A solution
$\phi$ of \eqref{seeq} is said to be unstable if there exist $\psi \not \sim \phi$
and a heteroclinic orbit from $\phi$ to $\psi$.

For a solution $\phi$ of \eqref{seeq}, denote by $\ms L_\phi$
the linear operator on $C^2(\bb T)$ given by
\begin{equation}
\label{rd-01}
\ms L_\phi h \; =\; (1/2) \Delta h \;-\;  V''(\phi) h \;.
\end{equation}
If $\phi$ is not constant, $\nabla\phi$ is an eigenfunction
associated to the eigenvalue $0$. A non-constant solution $\phi$ of
\eqref{seeq} is said to be hyperbolic if all eigenvalues of $\ms
L_\phi$ have non-zero real parts, except the eigenvalue $\lambda =0$,
whose associated eigenspace has dimension $1$.

The eigenvalue $0$ of the operator $\ms L_\phi$ is associated to the
orbit $\rho(t,\theta) = \phi(\theta + t)$. Actually, we prove in Lemma
\ref{lem5} that the cost of this orbit along a stationary set vanishes.
Moreover, the existence of a positive eigenvalue of
$\ms L_\phi$ is related to the existence of a heteroclinic orbit
starting from $\phi$ and, therefore, to the instability of $\phi$.

\begin{lemma}
\label{rd-l2}
Assume the conditions of Lemma \ref{rd-l1} and that all local maxima
of $V$ are non-degenerate: $V''(M_j)\not = 0$, $1\le j<\mf n$. Then,
for each non-constant solution $\phi$ of \eqref{seeq}, there exist
heteroclinic orbits from $\phi$ to $\phi_j$ and from $\phi$ to
$\phi_{j+1}$, where $\phi_k(\theta) = m_k$, $\theta\in\bb T$, and
$j=\max\{k < \mf n: m_k < \phi(\theta) \; \forall\, \theta\in \bb
T\}$. There exist also heteroclinic orbits from $\psi_j$ to $\phi_j$
and from $\psi_j$ to $\phi_{j+1}$, $1\le j < \mf n$, where
$\psi_j(\theta) = M_j$, $\theta\in \bb T$.
\end{lemma}

This result follows from Theorems 1.3 and 1.4 in \cite{frw}. We just
have to show that the hypotheses of these theorems are in force. As
$V'(0)<0<V'(1)$, the solutions are bounded below by $0$ and above by
$1$, so that $F$ is dissipative.

We claim that all non-constant solution $\phi$ of \eqref{seeq} are
hyperbolic. Indeed, fix such a function. As we have seen in the sketch
of the proof of Lemma \ref{rd-l1}, there exists $1\le j < \mf n$ such
that $m_j \le \phi(\theta) \le m_{j+1}$ for all $\theta\in \bb T$. Denote
by $\Pi$ the orbit map associated to the polynomial $F$
(cf. \cite[page 51]{Sch}). Since $V''(M_j)<0$, by \cite[Proposition
1.5.2]{Sch}, $F$ is an $A-B$ function in $(m_j,m_{j+1})$. Therefore, by
\cite[Theorem 2.1.3]{Sch}, $\Pi'(r)\not = 0$ for $r\not = M_j$. Hence,
by the proof of \cite[Lemma 4.4]{frw}, $\phi$ is hyperbolic.

We turn to the proof of the lemma. Fix $1\le j<\mf n$ and a
non-constant solution $\phi$ of \eqref{seeq} taking values in the
interval $[m_j, m_{j+1}]$.  We show that there exist heteroclinic
orbits from $\phi$ to $\phi_{j+1}$ and from $\psi_j$ to
$\phi_{j+1}$. Similar arguments permit to replace $\phi_{j+1}$ by
$\phi_{j}$.

The diagram of the Hamiltonian system $\dot p = q$, $\dot q = V' (p)$
shows that the periodic solutions of \eqref{seeq} which takes value in
the interval $[m_j, m_{j+1}]$ are either (i) $\phi_j$, $\phi_{j+1}$,
$\psi_j$ or (ii) a non-constant periodic solution whose maximal value
belongs to $(M_j, m_{j+1})$ and minimal value to $(m_{j}, M_j)$.
Moreover, if $\phi$, $\psi$ are such non-constant periodic solutions,
either
$\min_x \psi(x) < \min_x \phi(x)<M_j<\max_x \phi(x) < \max_x \psi(x)$
or the opposite.

We start with a heteroclinic orbit from $\psi_j$ to $\phi_{j+1}$.  We
may use the heteroclinic orbit from $M_j$ to $m_{j+1}$ for the ODE
$\dot x(t) = V'(x(t))$ to obtain a heteroclinic orbit from $\psi_j$ to
$\phi_{j+1}$ which remains constant in space.

Consider now a non-constant solution $\phi$ of \eqref{seeq} such that
$m_j \le \phi(\theta) \le m_{j+1}$. We repeat here the arguments of the
proof of Theorem 1.3 in \cite{frw} presented at the end of page
111. Let $z(h)$ be the number of strict sign changes of a function $h:
\bb T\to \bb R$. Since $z(\nabla\phi) \ge 2$, by
\cite[Proposition 3.1(b)]{frw}, the unstable dimension of $\phi$,
denoted by $i(\phi)$ in \cite{frw}, is larger than or equal to $1$. By
the positivity of the first eigenfunction of the operator $\mc
L_\phi$, one obtains a trajectory $\rho(t,\theta)$, $t\in \bb R$,
which solves \eqref{rdeq} and such that $\rho(t,\theta) >
\phi(\theta)$, $\lim_{t\to - \infty} \rho (t) = \phi$. Let $\psi=
\lim_{t\to +\infty} \rho (t)$, which exists in view of Lemma
\ref{rd-l0}. As $m_j\le\phi(\theta)\le m_{j+1}$, we have that $m_j \le
\psi(\theta) \le m_{j+1}$. By the Sturm property, $z(\rho(t) - \phi)$
decreases in time. Since it is equal to $0$ for $t$ close to
$-\infty$, $z(\psi - \phi)= 0$. Hence, $\psi$ can not be $\psi_j$ or
one of the non-constant solutions taking values in the interval $[m_j,
m_{j+1}]$. Thus, $\psi$ must be $\phi_j$ or $\phi_{j+1}$. Since $\psi
\ge \phi$, $\psi = \phi_{j+1}$, which proves the lemma.

We conclude this subsection with an example which fulfills the
assumptions of Lemma \ref{rd-l2}.  Fix $0<a<b$ and consider the
reaction-diffusion equation
\begin{equation}
\label{20}
\partial_t \rho \;=\; (1/2) \Delta \rho \;-\; V'(\rho)\;, \quad
\text{where}\quad V(\rho) \;=\; \frac{b}{4}\,(2\rho-1)^4 
\,-\, \frac{a}{2} \, (2\rho-1)^2 \;.
\end{equation}
This is the so-called Chafee-Infante equation \cite{ci}. It is clear
that the potential $V$ satisfies the assumptions of Lemma
\ref{rd-l2}. Actually, in this case all stationary solutions and all
heteroclinic orbits are known. We examine this example in Section
\ref{sec6}, where we present microscopic jump rates which fulfill the
hypotheses of Theorem \ref{sldp} below and whose hydrodynamic equation
is given by \eqref{20} with $0<a<b$.

\subsection{Hydrostatics}

Since the jump rate $c(\eta)$ is strictly positive, the Markov process
$\eta^{N}_{t}$ is irreducible in $X_N$. We denote by $\mu^{N}$ the
unique stationary probability measure under the dynamics.  We review
in this subsection the asymptotic behavior of the empirical measure
under the stationary state $\mu^N$.

Denote by $\mc P^N$ the probability measure on $\M_{+}$ defined by
$\cP^{N}:=\mu^{N}\circ(\pi^{N})^{-1}$.  The following theorem has been
established in \cite{lt}. It is a consequence of the law of large
numbers for the empirical measure, stated in Theorem
\ref{hydrodynamics}, and of the asymptotic behavior of the solutions
of the reaction-diffusion equation, stated in Lemma \ref{rd-l0}.

\begin{theorem}
\label{hsl}
The sequence of measures $\{\mathcal{P}^N:N\ge1\}$
is asymptotically concentrated on the set
$\M_{\rm sol}$. Namely, for any $\delta>0$, we have
\begin{align*}
\lim_{N\to\infty}\mathcal{P}^N\Big(\vte\in\M_+ : 
\inf_{\bar{\vte}\in\M_{\rm sol}} d(\vte , \bar{\vte})\ge\delta\Big) \;=\; 0\;.
\end{align*}
\end{theorem}

Note that this result does not exclude the possibility that the
stationary measure gives a positive weight to a neighborhood of an
unstable stationary solution of equation \eqref{rdeq}.

\subsection{Dynamical large deviations}

Let $\M_{+,1}$ be the closed subset of $\M_{+}$ consisting of all
absolutely continuous measures with density bounded by $1$:
\begin{equation*}
\mathcal{M}_{+,1}\;=\;\{\vte\in\mathcal{M}_+:
\vte(d\theta)=\rho(\theta)d\theta,\ 0\le\rho(\theta)\le1\ a.e.\ \theta\in\T\}\;.
\end{equation*}

Fix $T>0$, and denote by $C^{m,n}([0,T]\times\T)$, $m,n$ in
$\N_{0}\cup\{\infty\}$, the set of all real functions defined on
$[0,T]\times\T$ which are $m$ times differentiable in the first
variable and $n$ times in the second one, and whose derivatives are
continuous. Let $Q_{T,\eta}=Q^{N}_{T,\eta}$, $\eta\in X_{N}$, be the
probability measure on $D([0,T],\M_{+})$ induced by the measure-valued
process $\pi^{N}_{t}$ starting from $\pi^{N}(\eta)$.

For each path $\pi(t,d\theta)=\rho(t,\theta)d\theta$ in $D([0,T],\mathcal{M}_{+,1})$,
define the energy $\mathcal{Q}_{T}$ as
\begin{equation}
\label{i03}
\mathcal{Q}_{T}(\pi)\;=\;\sup_{G\in C^{0,1}([0,T]\times\T)}
\Big \{2\int_0^Tdt\ \lan\rho_t,\nabla G_t\ran
-\int_0^Tdt\int_{\T}d\theta\ G(t,\theta)^2 \Big\}\;.
\end{equation}
It is known (cf. \cite[Subsection 4.1]{blm2}) that the energy
$\mathcal{Q}_{T}(\pi)$ is finite if and only if $\rho$ has a
generalized derivative, denoted by $\nabla\rho$, and this generalized
derivative is square integrable on $[0,T]\times\T$:
\begin{equation*}
\int_{0}^{T} dt \ \int_{\T} d\theta \ |\nabla\rho(t,\theta)|^{2} <\infty\;.
\end{equation*}
Moreover, it is easy to see that the energy $\mathcal{Q}_{T}$
is convex and lower semicontinuous.

For each function $G$ in $C^{1,2}([0,T]\times\T)$, define the
functional $\bar{J}_{T,G}:D([0,T],\mathcal{M}_{+,1})\to\R$ by
\begin{align*}
\bar{J}_{T,G}(\pi) & \; =\;\lan\pi_T,G_T\ran -\lan\pi_{0},G_0\ran-
\int_0^Tdt\ \lan\pi_t, \partial_tG_t+\frac{1}{2}\De G_t\ran \\
&-\frac{1}{2}\int_0^Tdt\ \lan \chi(\rho_t), (\nabla G_t)^2\ran 
-\int_0^Tdt\ \big\{\lan B(\rho_t), e^{G_t}-1\ran + \lan D(\rho_t)\;, 
e^{-G_t}-1 \ran \big\},
\end{align*}
where $\chi(r)=r(1-r)$ is the mobility.  Let $J_{T,G}:
D([0,T],\mathcal{M}_+) \to[0,\infty]$ be the functional defined by
\begin{equation}
\label{07}
J_{T,G} (\pi) \;=\; 
\begin{cases}
\bar{J}_{T,G}(\pi)  & \text{if $\pi\in D([0,T],\mathcal{M}_{+,1})$}\;, \\
\infty  & \text{otherwise}\;,
\end{cases}
\end{equation}
and let $I_T:D([0,T],\mathcal{M}_+)\to[0,\infty]$ be the functional
given by
\begin{equation}
\label{i04}
I_T(\pi) \;=\; \begin{cases}
     \sup{J_{T,G}(\pi)} &
     \text{if $\mathcal{Q}_{T}(\pi)<\infty$}\;, \\
     \infty  & \text{otherwise}\;,
\end{cases}
\end{equation}
where the supremum is carried over all functions $G$ in
$C^{1,2}([0,T]\times\T)$. We sometimes abuse of notation by writing
$I_T(\rho)$ for $I_T(\pi)$ and we write $J_{G}$ for $J_{T,G}$ to
keep notation simple.

An explicit formula for the functional $I_T$ at smooth trajectories
was obtained in Lemma 2.1 of \cite{jlv}.  Let $\rho$ be a function in
$C^{2,3}([0,T]\times\T)$ with $c\le\rho \le 1-c$, for some $0< c
<1/2$.  Then, there exists a unique solution $H\in
C^{1,2}([0,T]\times\T)$ of the partial differential equation
\begin{align*}
\partial_t\rho \;=\; (1/2)\De\rho \,-\, \nabla(\chi(\rho)\nabla H)
\,+\, B(\rho)e^H \,-\, D(\rho)e^{-H} \;,
\end{align*}
and the rate functional $I_T(\rho)$ can be expressed as
\begin{align*}
I_T(\rho) \; & =\; \frac{1}{2}\int_0^Tdt\  \lan\chi(\rho_t), (\nabla H_t)^2\ran \\
& +\; \int_0^Tdt\  \lan B(\rho_t), 1-e^{H_t}+H_t e^{H_t} \ran \;+\; 
\int_0^Tdt\  \lan D(\rho_t), 1 - e^{-H_t} - H_t e^{-H_t} \ran \;.
\end{align*}

For a measurable function $\ga:\T\to[0,1]$, define the dynamical large
deviations rate function
$I_T(\cdot|\ga):D([0,T],\mathcal{M}_+)\to[0,\infty]$ as
\begin{equation*}
I_T(\pi|\ga) \;=\; \begin{cases}
     I_{T}(\pi) &
     \text{if $\pi(0,d\theta)=\ga(\theta)d\theta$}\;, \\
     \infty  & \text{otherwise}\;.
\end{cases}
\end{equation*}

The next result, which establishes a dynamical large deviations
principle for the measure-valued process $\pi_{\cdot}^{N}$ with rate
functional $I_{T}(\cdot|\ga)$ has been presented in \cite{lt} under
the assumption that the functions $B$ and $D$ are concave on
$[0,1]$. We refer to \cite{jlv,bl2,bl} for different versions.

\begin{theorem}
\label{dldp}
Assume that the functions $B$ and $D$ are concave in $[0,1]$.  Fix
$T>0$ and a measurable function $\ga:\T\to[0,1]$. Consider a sequence
$\eta^{N}$ of initial configurations in $X_{N}$ associated to $\ga$ in
the sense that $\lan\pi^{N}(\eta^N),G\ran$ converges to $\int_{\T}
G(\theta)\ga(\theta)d\theta$, as $N\uparrow\infty$, for all continuous function
$G:\bb T\to \bb R$. Then, the measure $Q_{T,\eta^{N}}$ on
$D([0,T],\mathcal{M}_{+})$ satisfies a large deviation principle with
the rate function $I_{T}(\cdot|\ga)$. That is, for each closed subset
$\mathcal{C}\subset D([0,T],\mathcal{M}_{+})$,
\begin{equation*}
\limsup_{N\to\infty}\frac{1}{N}\log{Q_{T, \eta^{N}}(\mathcal{C})} \;\le\; 
-\inf_{\pi\in\mathcal{C}}I_{T}(\pi|\ga)\;,
\end{equation*}
and for each open subset $\mathcal{O}\subset
D([0,T],\mathcal{M}_{+})$, 
\begin{equation*}
\liminf_{N\to\infty}\frac{1}{N}\log{Q_{T, \eta^{N}}(\mathcal{O})} \;\ge\; 
-\inf_{\pi\in\mathcal{O}}I_{T}(\pi|\ga)\;.
\end{equation*}
Moreover, the rate function $I_{T}(\cdot|\ga)$ is lower semicontinuous
and it has compact level sets.
\end{theorem}

\subsection{Static large deviations}\label{ssec:sldp}

We state in Theorem \ref{sldp} below the main result of this paper, a
large deviations principle for the empirical measure under the
stationary measure.  

Assume that the semilinear elliptic equation \eqref{seeq} admits at
most a finite number of solutions, modulo translations. More
precisely, assume that there exists $l\ge 1$ and density profiles
$\bar{\rho}_{1}, \cdots, \bar{\rho}_{l}$ in $C^2(\bb T)$, such that
\begin{equation*}
\M_{\rm sol}\;=\;\{ \bar\varrho_{i}(d\theta- \omega) 
= \bar\rho_{i}(\theta- \omega)d\theta : 1\le i\le l, \omega\in\T\}\;. 
\end{equation*}
Lemma \ref{rd-l1} provides conditions on the potential $V$ which
guarantee that this condition is in force.  Let $\M_{i}$, $1\le i\le
l$, be the subset of $\M_{\rm sol}$ given by
$\M_{i}=\{\bar\vte_i(d\theta - \omega) : \omega\in\T\}$.

Define the functionals $V_{i}:\M_{+}\to[0,\infty]$, $1\le i \le l$, by
\begin{equation}
\label{f01}
V_i(\vte) = \inf \Big\{I_T(\pi|\gamma):
T> 0 \,,\, \gamma (\theta)d\theta\in\M_{i} \,,\,
\pi\in D([0,T],\mc M_{+}) \,\text{ and }\, \pi_T = \vte \Big\}\;,
\end{equation}
which is the minimal cost to create the measure $\vte$ from the set
$\M_{i}$. We prove in Lemma \ref{lem5} that in the previous
variational formula we may replace the condition $\gamma
(\theta)d\theta\in\M_{i}$ by the more restrictive condition $\gamma(\theta) d\theta =
\bar\vte_{i}$ for some fixed $\bar\vte_{i} \in\M_{i}$: for all
$\bar\vte_i \in \mc M_i$,
\begin{equation}
\label{27}
V_i(\vte) = \inf \Big\{I_T(\pi|\gamma):
T> 0 \,,\, \gamma (\theta)d\theta = \bar\vte_i \,,\,
\pi\in D([0,T],\mc M_{+}) \,\text{ and }\, \pi_T = \vte \Big\}\;.
\end{equation}
By an abuse of notation, we sometimes write $V_i(\gamma)$ instead of
$V_i(\gamma (\theta) d\theta)$, where $\gamma:\bb T\to [0,1]$ is a density
profile.

By translation invariance, $V_{i}(\gamma) = V_{i}(\gamma')$ if
$\gamma' (\cdot) = \gamma(\cdot - \omega)$ for some $\omega\in \bb T$.
In particular, $V_i$ is constant on the set $\mc M_j$, $j\neq i$, and
$v_{ij}=V_{i}(\bar{\vte}_{j})$ is well defined, where $\bar\vte_j$ is
any element of $\mc M_j$.  Moreover, by choosing $T=1$ and
$\pi_t=(1-t)\bar\vte_i+t\bar\vte_j$, $t\in[0,1]$, in the infimum of
\eqref{27} yields that $v_{ij}$ is finite for any $i\neq j$.  Finally,
by Lemmata \ref{zero} and \ref{lem5}, $V_i(\vte)=0$ for any $\vte\in\M_i$.

Following \cite [Chapter 6]{fw1}, denote by $\ms T(i)$, $i\in \ms V
:=\{1, \cdots, l\}$, the set of all oriented, weighted, rooted trees
whose vertices are all the elements of $\ms V$ and whose root is
$i$. The edges are oriented from the child to the parent, and the
weight $v_{mn}$ is assigned to the oriented edge $(m,n)$.  Denote by
$\kappa (g)$ the sum of the weights of the tree $g\in \ms T(i)$ and by
$w_i$ the minimal weight of all trees in $\ms T(i)$:
\begin{equation*}
w_{i} \;=\; \min_{g\in \ms T (i)} \kappa (g) \;, \qquad
\kappa (g) \;=\; \sum_{(m, n)\in g} v_{mn}\;.
\end{equation*}
Since $v_{ij}$ is finite for $i\not = j$, so are $w_i$ and
$w=\min_{1\le i \le l} w_{i}$.  We will see below in \eqref{18} that
the non-negative parameter $w_i$ corresponds to the exponential weight
of a neighborhood of the set $\mc M_i$ under the stationary state.

Note that for all $i\not = j$,
\begin{equation}
\label{19}
w_i \;\le \; w_j \;+\; v_{ji} \;.
\end{equation}
Indeed, let $g$ be a graph in $\ms T(j)$ such that $w_j = \kappa(g)$.
Denote by $(a,b)$, $a\not = b\in \ms V$, the oriented edge where $a$
is the child and $b$ the parent.  Let $i'$ be the parent of $i$ in
$g$. Of course, $i'$ might be $j$. Denote by $g'$ the tree in $\ms
T(i)$ obtained from $g$ by adding the oriented edge $(j,i)$ and
removing the the edge $(i,i')$, and note that $\kappa(g) + v_{ji} =
\kappa(g') + v_{ii'}$. Since $w_i$ is the minimal value of $\kappa
(\tilde g)$, $\tilde g\in \ms T(i)$, $w_i \le \kappa (g') \le \kappa
(g') + v_{ii'} = \kappa (g) + v_{ji} = w_j + v_{ji}$.

For each $1\le i \le l$, define the functions
$W_{i}, W:\M_{+}\to[0,\infty]$ by
\begin{equation}
\label{06}
W_{i}(\vte) \;=\; w_{i} - w + V_{i}(\vte) \;, \quad
W(\vte) \;=\; \min_{1\le i \le l} W_{i}(\vte) \;.
\end{equation}
Note that for all $\vte\in \mc M_i$,
\begin{equation}
\label{18}
W(\vte) \;=\; \overline{w}_i \;:=\; w_i \,-\, w \;.
\end{equation}
Indeed, fix $\vte\in \mc M_i$. In view of the definition of $W$, we
have to show that $\min_{1\le j\le l} \{ \overline{w}_j + V_j(\vte)\} =
\overline{w}_i$. The minimum is less than or equal to $\overline{w}_i$
because $V_i(\vte)=0$. On the other hand, since $V_j(\vte)=v_{ji}$ and
since, by \eqref{19}, $w_i \le w_j + v_{ji}$, $\overline{w}_i \le \overline{w}_j +
V_j(\vte)$ for $j\not = i$.

The following theorem is the main result of this paper.

\begin{theorem}
\label{sldp}
Assume that the jump rates are strictly positive and that the
functions $B$ and $D$ are concave on $[0,1]$. Assume, furthermore,
that the semilinear elliptic equation \eqref{seeq} admits at most a
finite number of solutions, modulo translations. Then, the sequence of
probability measures $\{\cP^N; N\ge1\}$ satisfies a large deviation
principle on $\M_{+}$ with speed $N$ and rate function $W$.  Namely,
for each closed set $\mc C\subset\mc M_{+}$ and each open set $\mc
O\subset\mc M_{+}$,
\begin{align*}
& \limsup_{N\to\infty}\frac{1}{N}\log\mc P^N(\mc C) \;\le\; 
-\inf_{\vte\in \mc C} W(\vte) \; , \\
&\quad \liminf_{N\to\infty}\frac{1}{N}\log\mc P^N(\mc O) \;\ge\; 
-\inf_{\vte\in \mc O} W(\vte) \; .
\end{align*}
Moreover, the rate functional $W$ is bounded on $\M_{+,1}$, it is lower
semicontinuous, and it has compact level sets.
\end{theorem}

\subsection{The support of the stationary measure $\mu^N$} 

The next result improves on Theorem \ref{hsl} and asserts that the
stationary measure $\mu^N$ is concentrated on neighborhoods of stable
equilibria of the reaction-diffusion equation \eqref{rdeq}.  This
result has been conjectured in \cite{bl2,bl}.

As in the previous subsection, assume that the semilinear elliptic
equation \eqref{seeq} admits at most a finite number of solutions,
modulo translations.  Denote by $\mf M$ the set of local minima of
$V$, and by $\mc S \subset \mc M_{\rm sol}$ the set of associated
density profiles:
\begin{equation*}
\mc S \;=\; \{\varrho (d\theta) = m\, d\theta : m\in \mf M\}\;.
\end{equation*}
The elements of $\mc S$ are called stable solutions. Lemma \ref{l-f2}
justifies this terminology. It states that the quasi-potential
associated to each $\varrho \in \mc S$ is strictly positive outside
any neighborhood of $\varrho$. More precisely, for every
$\bar\varrho_i (d\theta) = \bar\rho_i(\theta) d\theta \in \mc S$ and $\varepsilon>0$,
there exists $c>0$ such that $\inf_{\gamma \not \in \mc
  B_\varepsilon(\bar\varrho_i)} V_i(\gamma) \ge c$, where $\mc
B_\varepsilon(\bar\varrho_i)$ represents a ball in $\mc M_+$ of radius
$\varepsilon$ centered at $\bar\varrho_i$.

Denote by $I_s$, $I_u\subset \{1, \dots, l\}$ the set of indices
associated to stable, unstable density profiles, respectively:
\begin{equation*}
I_s \;=\; \{ j : \bar\rho_j (\theta) \, d\theta \in \mc S\}\;, \quad
I_u \;=\; \{1, \dots, l\} \setminus I_s\;.
\end{equation*}

\begin{theorem}
\label{stablesol}
Assume that the hypotheses of Theorem \ref{sldp} are in force, and
that for all $i\in I_u$ there exists $j\in I_s$ such that
\begin{equation}
\label{ss-hyp}
v_{ij}\;=\; 0\;.
\end{equation}
Then, for all $\varepsilon>0$ there exist $c>0$ and $N_0\ge 1$ such that
for all $N\ge N_0$, 
\begin{equation*}
\mc P^N \Big( \M_+ \setminus\Big[ \bigcup_{j\in I_s} \mc
B_\varepsilon(\bar\varrho_j) \Big] \Big) \;\le\; e^{-c N} \;.
\end{equation*}
\end{theorem}

Of course, one expects the stationary measure $\mu^N$ to be
concentrated on neighborhoods of the density profiles associated to
the global minima of the potential $V$. This problem remains an open
question. A finer estimate than the one provided by the large
deviations might be needed to answer this open question. We refer to
\cite{fu07} for a similar problem in the context of the pinned Wiener
measure, and to \cite{bfo} and references therein for the study of the
concentration of measures in the situation where the rate functional
has more than one minimizer.

Lemma \ref{rd-l2} provides a set of sufficient conditions, expressed
in terms of the potential $V$, for assumption \eqref{ss-hyp} to
hold. Indeed, by Lemma \ref{l-f1}, if there exists a heteroclinic
orbit from $\mc M_i$ to $\mc M_j$, then $v_{ij}=0$. Hence, by Lemma
\ref{rd-l2}, if all zeros of $F$ are real, all critical points of $V$
are local minima or local maxima and all local maxima of $V$ are
non-degenerate, hypothesis \eqref{ss-hyp} is in force. 

Fix a stationary solution $\phi=\bar\rho_i$ of \eqref{seeq}.  There
are at least three different possible definitions of instability: (i)
the operator $\mc L_\phi$, defined in \eqref{rd-01}, has an eigenvalue
with positive real part.  (ii) there exists $\psi\not\sim \phi$ and a
heteroclinic orbit from $\phi$ to $\psi$.  (iii) $v_{ij}=0$ for some
$j\neq i$. We presented in the proof of Lemma \ref{rd-l2} a sketch of
the proof that (i) $\Rightarrow$ (ii) under some additional
hypotheses. Lemma \ref{l-f1} asserts that (ii) $\Rightarrow$ (iii). We
believe that the other implications hold, at least with some extra
assumptions, but we were not able to prove them.

\subsection{Comments and Remarks}

The characterization of the global attractor of the solutions of
reaction-diffusion equations \cite{cm89,frw} has only been achieved in
dimension $1$, not to mention the description of the heteroclinic
orbits. This is the main obstacle to extend the previous result to
higher dimensions. Although it is true that the dynamical large
deviations principle has been derived only in one dimension \cite{lt},
it should not be very difficult to extend it to higher dimensions.

The results presented in this article can be proved for
one-dimensional reaction-diffusion models with Dirichlet or Neumann
boundary conditions. The description of the heteroclinic orbits in
these contexts is simpler than the one with periodic boundary
conditions (cf. \cite{ci} for the case of Dirichlet boundary conditions).

As mentioned above, it is an open, and very appealing, problem to show
that the stationary measure $\mu^N$ is concentrated on neighborhoods
of the density profiles associated to global minima of the potential
$V$. To apply the method presented in the article to solve this
question would require a sharp estimate of the cost of the instanton,
the trajectory which drives the system from a stable equilibrium to
another. In view of Theorem \ref{tfrw} below, it is clear that the
instanton in the case of a double well potential with non-constant
stationary profiles is the trajectory which crosses the non-constant
stationary solution with one period. To estimate the cost of this
trajectory seems to be out of reach.

The previous questions lead us to the problem of the metastability of
the dynamics. It is challenging to describe the metastable behavior of
these reaction-diffusion models.

The hypothesis that the functions $B$ and $D$ are concave is only
needed in the proof of the dynamical large deviations principle
\cite{lt}, and we never use it in this paper. If one is able to prove this
dynamical result without the concavity assumption, the arguments
presented in this article provide a proof of the static large
deviations principle without the concavity assumption.

As mentioned in the introduction, the strategy of the proof consists
in adapting to our infinite-dimensional setting the Freidlin and
Wentzell approach \cite{fw1} to prove a large deviation principle for
the stationary state of a small perturbation of a dynamical system.
This has been done before in \cite{bg,f} for conservative evolutions
in contact with reservoirs. However, in the context of
reaction-diffusion models the existence of several stationary
solutions to the hydrodynamic equation introduces additional
difficulties.

The proof relies on a representation of the stationary state of the
reaction-diffusion model in terms of the invariant measure of a
discrete-time Markov chain induced by the successive visits to the
neighborhoods of the stationary solutions of the hydrodynamic
equation.

The proof of the static large deviations principle can be decomposed
in essentially three steps.  We first need to derive some regularity
properties of the dynamical large deviations rate functional. For
instance, that any trajectory which remains in a long time interval
far apart (in the $L^2$-topology) from the stationary solutions of the
hydrodynamic equation pays a strictly positive cost. Or that the
quasi-potential is lower semicontinuous in the weak topology.

The second step consists in obtaining sharp large deviations bounds
for the invariant measure of the discrete-time Markov chain. The final
step, whose proofs are similar to the ones presented in \cite{bg,f},
consists in estimating the minimal cost to create a measure starting
from a stationary solution of the hydrodynamic equation.

The topology is one of the main technical difficulties in the
argument.  The weak topology is imposed by the dynamical large
deviations principle which has been derived in this set-up, and one is
forced to prove all regularity properties of the rate functionals in
this topology. To overcome this obstacle, we systematically use the
smoothening properties of the hydrodynamic equation. Lemma \ref{conn}
is a good illustration of this strategy.

The article is organized as follows.  In Section \ref{sec3}, we
present the main properties of the weak solutions of the Cauchy
problem \eqref{rdeq}, and in Sections \ref{sec4a} and \ref{sec4b}, we
examine the dynamical and the static large deviations rate
functionals. These sections are purely analytical, and no
probabilistic argument is used. In Section \ref{sec5}, we prove the
static large deviations principle, and, in Section \ref{sec7}, the
concentration of the stationary measure $\mu^N$.  In Section
\ref{sec6} we present a reaction-diffusion model which fulfills the
hypotheses of Theorem \ref{stablesol}.

\section{The reaction-diffusion equation}
\label{sec3}

We present in this section several properties of the weak solutions of
the Cauchy problem \eqref{rdeq}. When we did not find a reference, we
present a proof of the result. Throughout this section and in the next
ones, $C_0$ represents a finite, positive constant which depends only
on $F$ and which may change from line to line.  As mentioned above,
this section and the following two ones are purely analytical, and no
probabilistic arguments appear.

We first define two concepts of solutions.

\begin{definition}
\label{wsol}
A measurable function $\rho:[0,T]\times\T\to[0,1]$ is said to be a
weak solution of the Cauchy problem \eqref{rdeq} in the layer
$[0,T]\times\T$ if for every function $G$ in $C^{1,2}([0,T]\times\T)$,
\begin{equation}
\label{weq}
\begin{aligned}
\lan \rho_{T}, G_{T}\ran - \lan \gamma, G_{0} \ran
- \int_{0}^{T} dt \, \lan & \rho_{t}, \partial_{t}G_{t} \ran  \\
&=\;\dfrac{1}{2} \int_{0}^{T} dt \, \lan \rho_{t}, \De G_{t} \ran 
+ \int_{0}^{T} dt \, \lan F(\rho_{t}), G_{t} \ran\;.  
\end{aligned}
\end{equation} 
\end{definition}

\begin{definition}
\label{msol}
A measurable function $\rho:[0,T]\times\T\to[0,1]$
is said to be a mild solution of the Cauchy problem \eqref{rdeq}
in the layer $[0,T]\times\T$ if for any t in $[0,T]$
\begin{align}
\label{meq}
\rho_{t} \;=\; P_{t}\gamma \;+\; \int_{0}^{t} P_{t-s}F(\rho_{s}) \, ds\;,
\end{align} 
where $\{P_{t} : t\ge 0\}$ stands for the semigroup on $L^{2}(\T)$ generated by
$(1/2)\De$.
\end{definition}

Next proposition asserts that the two notions of solutions are
equivalent. We refer to Proposition 6.3 of \cite{lt} for the proof.

\begin{proposition}
\label{exun}
Definitions \ref{wsol} and \ref{msol} are equivalent. Moreover, there
exists a unique weak solution of the Cauchy problem \eqref{rdeq}.
\end{proposition}

The next result is contained in Proposition 2.1 of \cite{ms95}.

\begin{proposition}
\label{reg}
Let $\rho$ be the unique weak solution of the Cauchy problem \eqref{rdeq}.
Then $\rho$ is infinitely differentiable over $(0,\infty)\times\T$.
\end{proposition}

Let $\bb Z_* = \bb Z\setminus \{0\}$, and let $c_0: \{0,1\}^{\bb Z_*}
\to \bb R_+$ be the cylinder function defined by $c_0(\xi) =
c(\xi^{(0)})$, where $\xi^{(0)}$ is the configuration of
$\{0,1\}^{\bb Z}$ defined by $\xi^{(0)}(x) = \xi(x)$, $x\not = 0$,
$\xi^{(0)}(0) = 0$. The cylinder function $c_1: \{0,1\}^{\bb Z_*} \to
\bb R_+$ is defined analogously with $\xi^{(0)}$ replaced by
$\xi^{(1)}$, where $\xi^{(1)}(0) = 1$. 

Note that $c_0$ and $c_1$ are strictly positive cylinder functions
because so is $c(\eta)$. Hence, if $\nu^*_\rho$ represents the
Bernoulli product measure on $\{0,1\}^{\bb Z_*}$ with density $\rho$,
the polynomial $\widehat B(\rho)$ defined by $\widehat B(\rho) =
E_{\nu^*_\rho}[c_0 (\eta)]$ is strictly positive. Similarly, the
polynomial $\widehat D(\rho)$ defined by $\widehat D(\rho) =
E_{\nu^*_\rho}[c_1 (\eta)]$ is strictly positive.

By definition, $B(\rho) \;=\; E_{\nu_\rho}[\{1-\eta(0)\} c(\eta)]
\;=\; (1-\rho) \, E_{\nu^*_\rho}[ c_0 (\eta)] = (1-\rho) \widehat
B(\rho)$, and $D(\rho) = \rho \widehat D(\rho)$. Hence,
\begin{equation}
\label{26}
B(\rho) \;=\; (1-\rho) \, \widehat B(\rho)\;, \quad 
D(\rho) \;=\; \rho \, \widehat D(\rho)\;,
\end{equation}
where $\widehat B(\rho)$ and $\widehat D(\rho)$ are strictly positive
polynomials. In particular, $F(0) = B(0) - D(0) = \widehat B(0) >0$
and $F(1) = B(1) - D(1) = - \widehat D(1) <0$.

Denote by $x_a(t)$, $0\le a\le 1$, the solution of the ODE 
\begin{equation}
\label{25}
\dot x(t) \;=\; F(x(t))  
\end{equation}
with initial condition $x(0) = a$.  Since $F(1)<0<F(0)$, $x_0(t)$
(resp. $x_1(t)$) is strictly increasing (resp. decreasing) and
$x_0(t)\to x_0$ (resp. $x_1(t)\to x_1$), where $x_0$ (resp. $x_1$) is
the smallest (resp. largest) solution of $F(x)=0$. The next result is a
simple application of the maximum principle.

\begin{lemma}
\label{lem1}
Let $\gamma:\bb T\to [0,1]$ be a density profile such that $a\le
\gamma(\theta) \le b$ for a.e. $\theta\in\bb T$. Denote by $\rho^\gamma(t,\theta)$ the
unique weak solution of \eqref{rdeq} with initial condition
$\gamma$. Then, $x_a(t) \le \rho(t,\theta) \le x_b(t)$ for all $t\ge 0$.
In particular, for any $t>0$, there exists $\varepsilon = \varepsilon
(t)>0$ such that $\varepsilon \le \rho^\gamma(t,\theta) \le 1- \varepsilon$
for all $\theta \in \bb T$ and all initial density profiles $\gamma : \bb T
\to [0,1]$.  Moreover, there exists $\delta>0$, depending only on $F$, such that
\begin{equation}
\label{11}
\delta \;\le\; \bar\rho_i \;\le\; 1-\delta
\end{equation}
for all $1\le i\le l$, $\bar\rho_i(\theta)\, d\theta \in \mc M_i$.
\end{lemma}

\begin{lemma}
\label{lem2}
There exists a finite constant $C_0$, depending only on $F$, such that
for any density profile $\gamma:\bb T \to [0,1]$ and any $t>0$,
\begin{equation*}
\Vert \rho_t \Vert^2_2 \;+\; \int_0^t 
\Vert \nabla\rho_s \Vert^2_2\, ds \;\le\; C_0 (1+t)\;,
\end{equation*}
where $\rho(t,\theta)$ stands for the unique weak solution of \eqref{rdeq} with
initial condition $\gamma$.
\end{lemma}

\begin{proof}
Fix a density profile $\gamma:\bb T \to [0,1]$. By Proposition
\ref{reg}, and since $\rho$ is the weak solution of \eqref{rdeq}, for
any $0<s<t$, by an integration by parts,
\begin{equation*}
\Vert \rho_t \Vert^2_2 \;=\; \Vert \rho_s \Vert^2_2 
\;-\; \int_s^t \Vert \nabla \rho_r \Vert^2_2 \, dr
\;+\; 2 \int_s^t dr \int_{\bb T}  \rho_r (\theta)  F(\rho_r(\theta))\, d\theta\;.
\end{equation*}
Since $\rho_r$ is absolutely bounded by $1$, we complete the proof of
the lemma by letting $s\downarrow 0$.
\end{proof}

A similar argument provides a bound on the distance between a
solution of the hydrodynamic equation and a constant stationary
solution. Recall that $\alpha \in (0,1)$ is an attractor of the ODE
\eqref{25} if there exists $\varepsilon>0$ such that the solution
$x(t)$ of the ODE with initial condition $x_0$ converges to $\alpha$
as $t\uparrow\infty$ if $|x_0-\alpha|<\varepsilon$.  Note in
particular that $F(\alpha)=0$ if $\alpha$ is an attractor.

\begin{lemma}
\label{lem10}
Let $\varepsilon>0$, let $\alpha$ be an attractor of the ODE
\eqref{25}, and let $\bar\rho_\alpha$ be the density profile given by
$\bar\rho_\alpha (\theta)=\alpha$, $\theta\in\bb T$.  There exists $\delta_{10}
= \delta_{10}(\varepsilon, \alpha)>0$ such that for any density
profile $\gamma: \bb T \to [0,1]$ such that $\Vert \gamma -
\bar\rho_\alpha \Vert_2 \le \delta_{10}$, $\rho_t$ converges in the
sup norm to $\bar\rho_\alpha$ as $t\uparrow\infty$, where $\rho_t(\theta) =
\rho(t,\theta)$ is the unique weak solution of \eqref{rdeq} with initial
condition $\gamma$. Moreover, for all $t\ge 1$, $\Vert \rho_t -
\bar\rho_\alpha\Vert_\infty \le \varepsilon$.
\end{lemma}

\begin{proof}
Fix $\varepsilon>0$, $\alpha\in (0,1)$ such that $F(\alpha)=0$, a
density profile $\gamma:\bb T \to [0,1]$ and recall the notation
introduced in the statement of the lemma.  Let $\rho$ be the weak
solution of \eqref{rdeq} with initial condition $\gamma$.  Repeating
the computation presented in the proof of Lemma \ref{lem2}, we obtain
that for every $0<s<t$,
\begin{equation*}
\Vert \rho_t - \bar\rho_\alpha \Vert^2_2 \;=\; 
\Vert \rho_s - \bar\rho_\alpha \Vert^2_2 
\;-\; \int_s^t \Vert \nabla \rho_r \Vert^2_2 \, dr
\;+\; 2 \int_s^t dr \int_{\bb T}  [\rho_r (\theta) - \bar\rho_\alpha(\theta)]\,   
F(\rho_r(\theta))\, d\theta\;.
\end{equation*}
Since $F(\alpha)=0$, we may subtract $F(\bar\rho_\alpha(\theta))$ from
$F(\rho_r(\theta))$ in the last integral and bound the product by $C_0
[\rho_r (\theta) - \bar\rho_\alpha(\theta)]^2$, where $C_0$ is the Lipschitz
constant of $F$. By letting $s\downarrow 0$ and then applying Gronwall
inequality, we obtain that
\begin{equation}
\label{22}
\Vert \rho_t - \bar\rho_\alpha \Vert^2_2 \;+\; 
\int_0^t \Vert \nabla \rho_r \Vert^2_2 \, dr
\;\le\; \Vert \rho_0 - \bar\rho_\alpha \Vert^2_2 \; e^{2C_0t}
\end{equation}
for all $t\ge 0$.

Choose $\varepsilon_0>0$ so that
$(\alpha-3\varepsilon_0, \alpha + 3\varepsilon_0)$ is contained in the
basin of attraction of $\alpha$ for the ODE \eqref{25} and set
$\varepsilon_1 = \min\{\varepsilon/2, \varepsilon_0\}$.  Choose
$\delta_{10} = \varepsilon_1 e^{-C_0}$ and let $\gamma$ be an initial
profile such that
$\Vert \gamma - \bar\rho_\alpha \Vert_2\le \delta_{10}$. By
\eqref{22}, for all $0\le t\le 1$,
\begin{equation}
\label{23}
\Vert \rho_t - \bar\rho_\alpha \Vert_2 \;\le\; \varepsilon_1
\quad\text{and}\quad
\int_0^1 \Vert \nabla \rho_r \Vert^2_2 \, dr \;\le\; \varepsilon^2_1\;.
\end{equation}
In particular, there exists $0\le s\le 1$ such that $\Vert \nabla
\rho_s \Vert_2 \le \varepsilon_1$, so that
\begin{equation*}
\sup_{\theta\not = \omega\in\bb T} \big| \rho(s,\theta) -
\rho(s,\omega) \big| \;\le\; 
\Vert \nabla \rho_s \Vert_1 \;\le\; 
\Vert \nabla \rho_s \Vert_2 \;\le\; \varepsilon_1\;.
\end{equation*}
Therefore, by \eqref{23}, for all $\theta \in\bb T$,
\begin{equation*}
\big| \rho(s,\theta) - \alpha \big| \;\le\;
\sup_{\theta\not = \omega\in\bb T} \big| \rho(s,\theta) - \rho(s,\omega) \big| \;+\;
\Vert \rho_s - \bar\rho_\alpha \Vert_1 \;\le\;
\varepsilon_1 \;+\; \Vert \rho_s - \bar\rho_\alpha \Vert_2 \;\le\;
\varepsilon
\end{equation*}
because $2\varepsilon_1 \le \varepsilon$.

Let $x_\pm(t)$, $t\ge 0$, be the solution of the ODE \eqref{25} with
initial condition $x_\pm (0) = \alpha \pm 2 \varepsilon_1$. By the
previous estimate, $x_-(0) \le \rho(s,\theta) \le x_+(0)$ for all $\theta\in\bb
T$. Hence, by Lemma \ref{lem1}, $x_-(t) \le \rho(s+t,\theta) \le x_+(t)$
for all $\theta\in\bb T$, $t\ge 0$. Since the basin of attraction of the
ODE is contained in $(\alpha-3\varepsilon_1, \alpha +
3\varepsilon_1)$, $x_\pm(t) \to \alpha$, as $t\to\infty$, and $x_-(0)
\le x_-(t) \le x_+(t) \le x_+(0)$ for all $t\ge 0$. In particular,
$\rho_t$ converges in the sup norm to $\bar\rho_\alpha$, as
$t\to\infty$, and $\Vert \rho_{s+t} - \bar\rho_\alpha \Vert_\infty \le
\Vert \rho_s - \bar\rho_\alpha \Vert_\infty \le \varepsilon$ for all
$t\ge 0$.  This completes the proof of the lemma because $s\le 1$.
\end{proof}

Similar arguments permit to estimate the distance between two solutions
of the reaction-diffusion equation \eqref{rdeq}.

\begin{lemma}
\label{est3.1}
There exists a constant $C_0>0$ such that for any weak solutions
$\rho^{j}$, $j=1,2$, of the Cauchy problem \eqref{rdeq} with initial
profile $\rho_{0}^{j}$ and for any $t>0$, 
\begin{equation*}
\|\rho_{t}^{1}-\rho_{t}^{2}\|_2 \;\le\;
e^{C_0 t}\|\rho_{0}^{1}-\rho_{0}^{2}\|_2 \;.
\end{equation*}
\end{lemma}

\begin{proof}
From \eqref{meq}, for any $t\ge0$ and $j=1,2$, 
\begin{equation*}
\rho_{t}^{j}\;=\;
P_{t}\rho_{0}^{j} + \int_{0}^{t}P_{t-s}F(\rho_{s}^{j}) ds \;.
\end{equation*}
Therefore
\begin{align*}
\|\rho_{t}^{1}-\rho_{t}^{2}\|_{2}
& \;\le\; \| P_{t}(\rho_{0}^{1}-\rho_{0}^{2})\|_{2}+
\int_{0}^{t}\| P_{t-s}(F(\rho_{s}^{1})-F(\rho_{s}^{2}))\|_{2} \, ds \\
& \;\le\; \|\rho_{0}^{1}-\rho_{0}^{2}\|_{2}+
\|F'\|_{\infty}\int_{0}^{t}\|\rho_{s}^{1}-\rho_{s}^{2}\|_{2} \,ds\;.
\end{align*}
In the last inequality, we used the fact that the operator norm of
$P_{t}$ is equal to $1$. To conclude the proof of the lemma, it
remains to apply Gronwall inequality.
\end{proof}

For each function $\rho$ in $L^{2}(\T)$, let $\mathbb
B_{\delta}(\rho)$, $\delta>0$, be the $\delta$-open neighborhood of
$\rho$ in $L^{2}(\T)$.  Recall also that we denote by $\mathcal
B_{\delta}(\vte)$ be the $\delta$-open neighborhood of $\vte$ in
$\M_{+}$. We sometimes represent the neighborhood $\mathcal
B_{\delta}(\vte)$ by $\mc B_\delta(\gamma)$ when $\vte (d\theta) =
\gamma(\theta) d\theta$.

\begin{lemma}
\label{rep}
Let $\bar\rho:\T\to[0,1]$ be a classical solution to the equation
\eqref{seeq}, and set $\bar\vte(d\theta)=\bar\rho(\theta)d\theta$. For
any $\varepsilon>0$ and $0<T<T'$, there exists
$\delta_{11} = \delta_{11}(\varepsilon, T,T') \in (0,\varepsilon)$
such that for any density profile $\rho_{0}:\T\to[0,1]$,
$\rho_{0}(\theta)d\theta$ in $\mathcal{B}_{\delta_{11}}(\bar\vte)$, it
holds that
$\rho(t,\theta) d\theta\in \mathcal{B}_{\varepsilon}(\bar{\vte})$ for
all $0\le t \le T'$ and that $\rho_t\in\bb B_\e(\bar\rho)$ for all
$T\le t \le T'$, where $\rho_{t}(\theta)=\rho(t,\theta)$ is the unique
weak solution of the Cauchy problem (\ref{rdeq}) with initial
condition $\rho_{0}$.
\end{lemma}

\begin{proof}
Fix $\varepsilon>0$ and $0<T<T'$.  Let $\zeta_1 = (1/3) \varepsilon
e^{-C_0T'}$, where $C_0$ is the constant appearing in Lemma
\ref{est3.1}. Fix a density profile $\rho_{0}:\bb T\to [0,1]$, and let
$\pi_t(d\theta) = \rho(t,\theta) d\theta$, where $\rho_{t}$ is the unique weak
solution of the Cauchy problem (\ref{rdeq}) with initial condition
$\rho_{0}$.  Recall the definition of the complete orthogonal normal
basis $\{ e_k;k\in\Z\}$ introduced just before \eqref{01}.  By
\eqref{01} and \eqref{weq}, for any $t\ge 0$,
\begin{equation*}
d(\pi_t, \bar\vte) \;\le\; d(\pi_0, \bar\vte) \;+\;
\sum_{k\in\bb Z} \frac 1{2^{|k|}}\, \Big| 
\dfrac{1}{2} \int_{0}^{t} ds \, \lan \rho_{s}, \De e_k \ran 
+ \int_{0}^{t} ds \, \lan F(\rho_{s}), e_k \ran\Big| \;.
\end{equation*}
The first term on the right hand side is bounded by
$\varepsilon/2$ if $\rho_0 \in \mathcal{B}_{\zeta_2}(\bar\vte)$,
where $\zeta_2=\varepsilon/2$, while the second one is less than
or equal to
\begin{equation*}
t \sum_{k\in\bb Z} \frac 1{2^{|k|}}\, \big\{ C_0 + (2\pi k)^2 \big\}
\;=\; C_0 \,t\;,
\end{equation*}
because $\rho_s$ is bounded by $1$, $F$ by a constant $C_0$ and
$\Vert e_k\Vert_2 =1$. Hence, if $T_1=\varepsilon/2C_0$ and $\pi_{0} \in
\mathcal{B}_{\zeta_2}(\bar\vte)$, 
\begin{equation}
\label{02}
\pi_t \;\in\; \mathcal{B}_{\varepsilon}(\bar\vte) \text{ for all $0\le
  t\le T_1$}\;.
\end{equation}

We turn to the $L^2$-estimate. From the equation \eqref{meq}, we have
\begin{equation}
\label{bound3.1}
\begin{aligned}
\|\rho_{t}-\bar\rho\|_{2}
&\;\le\; \| P_{t}(\rho_{0}-\bar\rho)\|_{2}+
\int_{0}^{t} \big \| P_{t-s} \big[F(\rho_{s})-F(\bar\rho)\big] \,
\big\|_{2} \, ds  \\
&\;\le\;  \| P_{t} ( \rho_{0}-\bar\rho) \|_{2}  +  t\|F'\|_{\infty} \;,
\end{aligned}
\end{equation}
since the operator norm of $P_{t}$ is equal to $1$.
Let $\tilde\rho=\rho_{0}-\bar\rho$ and, for each $t>0$,
$\tilde\rho_{t}=P_{t}\tilde\rho$.
It is easy to see that, for any $k\in\Z$,
\begin{equation*}
\lan \tilde\rho_{t}, e_{k}\ran \;=\;
\lan \tilde\rho, P_{t}e_{k} \ran \;=\;
e^{-2\pi^{2}k^{2}t}\lan \tilde\rho, e_{k} \ran \;.
\end{equation*}
Therefore, from Parseval's relation, 
\begin{equation}
\label{bound3.2} 
\|\tilde\rho_{t}\|^{2}_{2}\;=\;
\sum_{k\in\Z}e^{-4\pi^{2}k^{2}t}\lan\tilde\rho, e_{k}\ran^{2} \;.
\end{equation}

Set $T_2:= \min\{ (\zeta_1/2\|F'\|_{\infty}), T_1, T\}$ and choose a
large enough positive integer $k_{1}$ so that
\begin{equation*}
\sum_{|k|>k_1} e^{-4\pi^{2}k^{2}T_2} \;\le\;\zeta^{2}_1/8 \;.
\end{equation*}
To estimate the first terms of the series, observe that
\begin{equation*}
\sum_{|k|\le k_1} e^{-4\pi^{2}k^{2}t} \lan\tilde\rho,
e_{k}\ran^{2} \;\le\; 4^{k_1} \Big( \sum_{|k|\le k_1} 2^{-|k|} 
\, \big| \lan\tilde\rho, e_{k}\ran \big| \Big)^2 \;\le\;
4^{k_1} d(\rho_0,\bar\rho)^2\;.
\end{equation*}
Hence, if we set $\zeta_3 = \zeta_1/2^{k_1+2}$, this last expression
is bounded by $\zeta^2_1/16 \le \zeta^2_1/8$ provided $\rho_{0}(\theta)d\theta$
belongs to $\mathcal{B}_{\zeta_3}(\bar\vte)$. Therefore, by
\eqref{bound3.1}, \eqref{bound3.2} and the choice of $\zeta_3$,
\begin{equation}
\label{03}
\|\rho_{T_2}-\bar\rho\|_{2}  \;\le\; \zeta_1
\end{equation}
if $\pi_{0} \in \mathcal{B}_{\zeta_3}(\bar\vte)$.

Let $\delta_{11}=\min\{\zeta_2, \zeta_3\}$, and note that
$\delta_{11}$ depends only on $\varepsilon$, $T$, $T'$.  By
\eqref{02} and \eqref{03}, and since $T_2\le T_1$, for all $\pi_{0}
\in \mathcal{B}_{\delta_{11}}(\bar\vte)$,
\begin{equation*}
\pi_t \;\in\; \mathcal{B}_{\varepsilon}(\bar\vte) \text{ for all $0\le
  t\le T_2$} \quad\text{and}\quad 
 \|\rho_{T_2}-\bar\rho\|_{2}  \;\le\; \zeta_1\;.
\end{equation*}
By Lemma \ref{est3.1}, by the previous estimate and by definition of
$\zeta_1$, for all $T_2\le t\le T'$, $\pi_{0} \in
\mathcal{B}_{\delta_{11}}(\bar\vte)$,
\begin{equation*}
\|\rho_{t}-\bar\rho\|_{2}  \;\le\; e^{C_0 (t-T_2)}
\|\rho_{T_2}-\bar\rho\|_{2} \;\le\; e^{C_0 T'} \zeta_1 
\;\le\; \varepsilon/3\;.
\end{equation*}
Moreover, by this bound and by \eqref{05}, for all $T_2\le t\le T'$,
$\pi_{0} \in \mathcal{B}_{\delta_{11}}(\bar\vte)$,
\begin{equation*}
d(\pi_t, \bar\vte) \;\le\; 3\, \|\rho_{t}-\bar\rho\|_{2} 
\;\le\; \varepsilon\;.
\end{equation*}
This completes the proof of the lemma since $T_2\le T$.
\end{proof}

The previous results permit to strengthen Lemma \ref{lem10}.

\begin{lemma}
\label{lem11}
Let $\varepsilon>0$, let $\alpha$ be an attractor of the ODE
\eqref{25}, and let $\bar\vte_\alpha (d\theta) =
\bar\rho_\alpha(\theta) d\theta$, $\bar\rho_\alpha(\theta)=\alpha$,
$\theta\in\bb T$. There exists $\delta_{12}=\de_{12}(\e,\a)>0$ such
that for any density profile $\gamma: \bb T \to [0,1]$ such that
$\gamma(\theta) d\theta \in \mc B_{\delta_{12}}(\bar\vte_\alpha)$,
$\rho_t$ converges in the sup norm to $\bar\rho_\alpha$, as
$t\uparrow\infty$, where $\rho_t(\theta) = \rho(t,\theta)$ is the
unique weak solution of \eqref{rdeq} with initial condition
$\gamma$. Moreover, $\pi_t(d\theta) = \rho(t,\theta) d\theta$ belongs
to $\mc B_{\varepsilon}(\bar\vte_\alpha)$ for all $t\ge 0$.
\end{lemma}

\begin{proof}
Fix $\varepsilon>0$.  Denote by $\zeta_1$ the constant $\delta_{10}
= \delta_{10}(\varepsilon, \alpha)$ provided by Lemma
\ref{lem10}. Let $\zeta_2=\min \{\zeta_1, \varepsilon\}$, and let
$\delta_{12}$ be the constant $\delta_{11} = \delta_{11}(\zeta_2, 1/2,
2)$ provided by Lemma \ref{rep} with $\bar\rho = \bar\rho_\alpha$.

Fix $\gamma: \bb T \to [0,1]$ such that $\gamma(\theta) d\theta \in \mc
B_{\delta_{12}}(\bar\vte_\alpha)$. Denote by $\rho_t(\theta) = \rho(t,\theta)$ the
weak solution of the hydrodynamic equation with initial condition
$\gamma$. By Lemma \ref{rep}, $\Vert \rho_1 - \bar \rho_\alpha\Vert_2
\le \zeta_1$ and $\pi_t(d\theta) = \rho(t,\theta) d\theta$ belongs to $\mc
B_{\varepsilon}(\bar\vte_\alpha)$ for all $t\le 2$. 

Since $\Vert \rho_1 - \bar \rho_\alpha\Vert_2 \le \zeta_1$, by Lemma
\ref{lem10}, $\rho_t$ converges in the sup norm to $\bar\rho_\alpha$
as $t\uparrow\infty$, and $\Vert \rho_t - \alpha\Vert_\infty \le
\varepsilon$ for all $t\ge 2$. In particular, $\pi_t(d\theta) =
\rho(t,\theta) d\theta$ belongs to $\mc B_{\varepsilon}(\bar\vte_\alpha)$
for all $t\ge 2$.
\end{proof}

\section{The dynamical rate function}
\label{sec4a}

We present in this section some features of the dynamical rate
function needed to prove the properties of the static rate function
stated in the next section. The main result of the section asserts
that a trajectory can not remain too long far in the $L^2$-topology
from all stationary solutions of the hydrodynamic equation without
paying a fixed positive cost. 

The first four lemmata have been proved in \cite[Section 4]{lt} for
the rate functional $I_{T}(\cdot|\ga)$. The same arguments apply the
functional $I_{T}$. The first three extract information on the
trajectory $\pi(t,d\theta)$ from the finiteness of the large
deviations rate functional. Lemma \ref{regu} states that a trajectory
with finite rate function is a continuous path in
$D([0,T],\M_+)$. This lemma is repeatedly used in the rest of this
paper, and therefore, is used without any further mention.

\begin{lemma}
\label{regu}
Fix $T>0$. Let $\pi$ be a path in $D([0,T], \M_+)$ such that
$I_T(\pi)$ is finite.  Then $\pi$ belongs to $C([0,T], \M_{+,1})$.
\end{lemma}

Lemma \ref{energy} states that the density
$\rho(t,\theta)$ of a trajectory $\pi(t,d\theta)$ with finite rate function
belongs to $\mc H_1$ for almost all $t$, where $\mc H_1$ represents
the space of functions $f: \bb T \to \bb R$ which have a general
derivative in $L^2(\bb T)$.

\begin{lemma}
\label{energy}
There exists a finite constant $C_0>0$ such that for any $T>0$ and for
any path $\pi(t,d\theta)=\rho(t,\theta)d\theta$ in $D([0,T],\mathcal{M}_{+,1})$ with
finite energy, 
\begin{equation*}
\mc E_{T}(\rho):=\int_0^{T} dt\ \int_{\T} d\theta\ 
\frac{|\nabla\rho(t,\theta)|^{2}}{\chi(\rho(t,\theta))}
\;\le\; C_0\,\{ I_{T}(\pi)+ T + 1\}\;.
\end{equation*}
\end{lemma}

The next result characterizes the weak solutions of the hydrodynamic
equation as the trajectories at which the dynamical large deviations
rate functional vanishes.

\begin{lemma}
\label{zero}
Fix $T>0$.  The density $\rho$ of a path $\pi(t,d\theta)=\rho(t,\theta)d\theta$ in
$D([0,T], \M_{+,1})$ is the weak solution of the Cauchy problem
\eqref{rdeq} with initial profile $\ga$ if and only if $I_{T}(\pi|\ga)
= 0$.  Moreover, in that case
\begin{equation*}
\mc E_{T}(\rho) \;<\; \infty.
\end{equation*}
\end{lemma}

The next result is extremely useful. In the expression of the functionals
$J_{T,G}$ and $I_T$,  terms appearing such as $\int_0^T \< B(\rho_t), G_t\>\,
dt$, where $G$ is a smooth function, are not continuous for the
weak topology, but only for the $L^1$-topology. The next lemma
establishes that if the cost of a sequence $\pi^n$ of trajectories is
uniformly bounded and if this sequence converges weakly to some
trajectory $\pi$, then the sequence of density profiles converges in
$L^2$-topology. The proof of this result follows from the computations
presented in the proof of Theorem 4.7 in \cite{lt}.

\begin{lemma}
\label{conv}
Fix $T>0$.  Let $\{\pi^{n}(t,d\theta)=\rho^{n}(t,\theta)d\theta:n\ge1\}$ be a
sequence of trajectories in $D([0,T], \M_{+,1})$. Assume that there
exists a finite constant $C$ such that
\begin{equation*}
\sup_{n\ge1} \, I_{T}(\pi^{n}) \;\le\; C.
\end{equation*}
If $\rho^{n}$ converges to $\rho$ weakly in $L^{2}(\T\times[0,T])$,
then $\rho^{n}$ converges to $\rho$ strongly in $L^{2}(\T\times[0,T])$.
\end{lemma}

Recall the definition of the neighborhoods $\mathbb B_{\delta}(\rho)$ and
$\mathcal B_{\delta}(\vte)$ introduced just before the statement of Lemma
\ref{rep}.  For each $\delta>0$ and $T>0$, denote by $\mathbb{D}_{T,\delta}$
the set of trajectories $\pi(t,d\theta)=\rho(t,\theta)d\theta$ in $D([0,T],\M_{+,1})$
such that $\rho_{t}\notin\mathbb{B}_{\delta}(\bar{\rho})$ for all $0\le
t\le T$ and $\bar{\rho}\in S$.

Next lemma states that a trajectory can not stay a long time interval
far, in the $L^2$-topology, from all stationary solutions of the
hydrodynamic equation without paying an appreciable cost. This result
plays a fundamental role in the proof of the lower semicontinuity of
the functional $W$. To enhance its interest, note that $L^2$-neighborhoods
are much thinner than the neighborhoods of the weak
topology.

\begin{lemma}
\label{costs}
For every $\delta>0$ there exists $T=T(\delta)>0$ such that
\begin{equation*}
\inf_{\pi\in \mathbb{D}_{T,\delta}} I_{T}(\pi) \;>\; 0 \;.
\end{equation*}
\end{lemma}

\begin{proof}
Assume that the assertion of the lemma is false.
Then, there exists some $\delta>0$ such that,
for any $n\in\N$, 
\begin{equation*}
\inf_{\pi\in \mathbb{D}_{n,\delta}} I_{n}(\pi) \;=\; 0 \;.
\end{equation*}
In this case there exists a sequence of trajectories
$\{\pi^{n}(t,d\theta)=\rho^{n}(t,\theta)d\theta: n\ge1\}$, $\pi^n \in
\mathbb{D}_{n,\delta}$, such that $I_{n}(\pi^{n})\le 1/n$.  Since $I_{T}$
has compact level sets, by using a Cantor's diagonal argument and
passing to a subsequence if necessary, we obtain a path
$\pi(t,d\theta)=\rho(t,\theta)d\theta$ in $D(\R_{+}, \mathcal{M}_{+,1})$ such that
$\pi^{n}$ converges to $\pi$ in $D([0,T], \mathcal{M}_{+})$ for any
$T>0$. Moreover, by Lemma \ref{conv}, $\rho^{n}$ converges to $\rho$
strongly in $L^{2}([0,T]\times \T)$ for all $T>0$.

Since $I_T$ is lower semicontinuous, $I_{T}(\pi)=0$ for all $T>0$.
By Lemma \ref{zero}, the density of $\pi$, denoted by $\rho$ so that
$\pi(t,d\theta) = \rho(t,\theta) d\theta$, is the unique weak solution of the
equation \eqref{rdeq} with initial condition $\rho(0,\cdot)$. Hence,
by Lemma \ref{rd-l0}, $\rho_{t}$ converges in $C^2(\bb T)$ to
some density profile $\rho_{\infty} \in S$.  Therefore, there exists
some $T_0>0$ such that
\begin{equation*}
\|\rho_{t}-\rho_{\infty}\|_{2}\;\le\; \delta/2 \;,
\end{equation*}
for any $t\ge T_0$. Hence, since $\pi^n$ belongs to
$\mathbb{D}_{n,\delta}$, for $n\ge T_0+1$
\begin{align*}
\int_{0}^{T_0+1}\|\rho_{t}^{n}-\rho_{t}\|_{2} \, dt
& \;\ge\; \int_{T_0}^{T_0+1}\|\rho_{t}^{n}-\rho_{t}\|_{2} \, dt \\
& \;\ge\; \int_{T_0}^{T_0+1}\big(\|\rho_{t}^{n}-\rho_{\infty}\|_{2} -
\|\rho_{t}-\rho_{\infty}\|_{2}\big) \, dt \\
& \;\ge\; \delta - \delta/2 = \delta/2 \;,
\end{align*}
which contradicts the strong convergence of $\rho^{n}$ to $\rho$
in $L^{2}([0,T_0+1]\times \T)$ and we are done.
\end{proof}

Analogously, for each $\delta>0$ and $T>0$,
denote by $\mathcal{D}_{T,\delta}$ the set of trajectories
$\pi(t,d\theta)=\rho(t,\theta)d\theta$ in $D([0,T],\mathcal{M}_{+,1})$ such that
$\pi_{t}\notin\mathcal{B}_{\delta}(\bar\vte)$ for all $0\le t\le T$ and
$\bar\vte\in\M_{\rm sol}$.
A similar result also holds for the set $\mathcal{D}_{T,\delta}$.

\begin{corollary}
\label{costw}
For every $\delta>0$, there exists $T>0$ such that
\begin{equation*}
\inf_{\pi\in \mathcal{D}_{T,\delta}} I_{T}(\pi) \;>\; 0\;.
\end{equation*}
\end{corollary}

\begin{proof}
The assertion follows from Lemma \ref{costs} and the fact that
\begin{equation*}
\{\vte(d\theta)=\rho(\theta)d\theta : \rho\in\mathbb{B}_{\delta}(\bar{\rho})\}
\subset \mathcal{B}_{3 \delta}(\bar{\rho}),
\end{equation*}
for every $\bar{\rho}\in S$ and every $\delta>0$ in view of \eqref{05}.
\end{proof}

In the proof of the static large deviations principle, it will be
useful to estimate $V_i(\vte)$ for some measure $\vte(d\theta) = \gamma(\theta)
d\theta$. If $\gamma$ is a smooth density profile, this can be achieved by
joining $\bar\vte_i$ to $\vte$ through a linear interpolation $\pi_t =
(1-t/T) \bar\vte_i + (t/T) \vte$, $T>0$, and by estimating the cost of
the path $\pi$. This is the content of Lemma \ref{qpbound}. For a
general measure $\vte(d\theta) = \gamma(\theta) d\theta$, we need first to smooth the
density profile $\gamma$. We use the hydrodynamic equation to do
that. Fix $\varepsilon>0$ small, and denote by $\rho(t,\theta)$ the
solution of the hydrodynamic equation starting from $\gamma$, $0\le
t\le \varepsilon$. By Proposition \ref{reg}, $\rho(\varepsilon,
\cdot)$ is smooth. We may use the first part of this argument to joint
$\bar\vte_i$ to $\rho(\varepsilon, \theta) d\theta$. To connect
$\rho(\varepsilon, \theta) d\theta$ to $\gamma(\theta) d\theta$ we use the backward path
$\hat \pi_t(d\theta) = \rho(\varepsilon-t, \theta) d\theta$, $0\le t\le
\varepsilon$. The cost of this path is estimated in the next lemma.

\begin{lemma}
\label{bound4.1}
There exists a constant $C_0>0$ such that for any $T>0$, any weak
solution $\rho$ of \eqref{rdeq}, and any classical solution
$\bar{\rho}$ to the equation \eqref{seeq},
\begin{equation*}
I_T(\pi) \;\le\; C_0
\{T + \|\rho_T - \bar\rho\|_1 + \|\rho_0 - \bar\rho\|_{1}\}\;,
\end{equation*}
where $\pi$ is the trajectory defined by $\pi(t,d\theta) =\rho(T-t,\theta)d\theta$.
\end{lemma}

\begin{proof}
For any test function $G\in C^{1,2} ([0,T]\times\T)$,
$J_{T,G}(\pi)$ can be rewritten as
\begin{align*}
J_{T,G} & (\pi) \;=\; \int_0^Tdt\ \lan\nabla \rho_t, \nabla \widehat G_t\ran
- \frac{1}{2}\int_0^Tdt\ \lan\chi(\rho_t), (\nabla\widehat G_t)^2\ran \\
& - \int_0^Tdt\ \lan B(\rho_t), e^{\widehat G_t}+ \widehat G_t -1\ran
- \int_0^Tdt\ \lan D(\rho_t), e^{-\widehat G_t}- \widehat G_t -1\ran \;,
\end{align*}
where $\widehat G(t,\theta) = G(T-t,\theta)$.  The first line on the right hand
side is bounded above by
\begin{equation}
\label{f03}
\frac{1}{2} \int_{0}^{T}dt\ \int_{\T} d\theta\ 
\dfrac{|\nabla\rho(t,\theta)|^{2}}{\chi(\rho(t,\theta))} \;=\;
\frac{1}{2}\, \mc E_T(\rho)\;.
\end{equation}

Since for any $0<\rho<1$ and any $a\in \bb R$
\begin{align*}
-B(\rho)e^a +D(\rho) a +D(\rho) \;\le\; D(\rho) \log (D(\rho)/B(\rho)) \;, \\
-D(\rho)e^{-a} -B(\rho) a + B(\rho) \;\le\; B(\rho) \log (B(\rho)/D(\rho)) \;,
\end{align*}
and since $B(\rho) = (1-\rho) \widehat B(\rho)$, $D(\rho) = \rho
\widehat D(\rho)$, where $\widehat B(\rho)$, $\widehat D(\rho)$ are
the strictly positive functions introduced in \eqref{26}, the second
line on the right hand side is bounded above by
\begin{equation}
\label{f02}
\begin{split}
\int_0^Tdt\ &\lan D(\rho_t)\log(D(\rho_t)/B(\rho_t))
+ B(\rho_t)\log (B(\rho_t)/D(\rho_t))\ran\\
& \;\le\; C_0 T - \int_0^Tdt\ 
\lan D(\rho_t)\log (1-\rho_t)+ B(\rho_t)\log(\rho_t)\ran
\end{split}
\end{equation}
To estimate this last term, let $h(x) = x\log x + (1-x) \log (1-x)$
and note that
\begin{align*}
\partial_t h(\rho_t)
\;& =\; [\log(\rho_t)-\log(1-\rho_t)]\, \partial_t\rho_t\\
& =\; [\log(\rho_t)-\log(1-\rho_t)]\,
[\frac{1}{2}\Delta\rho_t + F(\rho_t)] \;.
\end{align*}
This equation is justified since, by Proposition \ref{reg}, any weak
solution of the equation \eqref{rdeq} is smooth for $t>0$.  Therefore,
\begin{align*}
& - \int_0^Tdt \ \lan D(\rho_t)\log (1-\rho_t)
+ B(\rho_t)\log(\rho_t)\ran \;=\; - \; \<h(\rho_T)\> \;+\; \<h(\rho_0)\>\\
&\qquad -\frac{1}{2}\, \mc E_T(\rho)- \int_0^Tdt \ 
\lan D(\rho_t)\log (\rho_t)+ B(\rho_t)\log(1-\rho_t) \ran   \;.
\end{align*}
Adding and subtracting $\<h(\bar\rho)\>$, since $\bar\rho$ takes value
in a compact interval of $(0,1)$, the first two terms on the right
hand side can be bounded above by
\begin{equation*}
C_0 \{\|\rho_T - \bar\rho\|_1 + \|\rho_0 - \bar\rho\|_{1}\}\;,
\end{equation*}
for some $C_0>0$. Since the last term in the penultimate displayed
formula is bounded by $C_0T$, we have shown that \eqref{f02} is less
than or equal to
\begin{equation*}
-\frac{1}{2}\, \mc E_T(\rho) \;+\; C_0 \big\{ 
T + \|\rho_T - \bar\rho\|_1 + \|\rho_0 - \bar\rho\|_{1}\big\}\;.
\end{equation*}
This estimate together with \eqref{f03} completes the proof of the
lemma.
\end{proof}

The last result of this section states that the cost to move inside a
set of static solutions is zero.

\begin{lemma}
\label{lem5}
Fix $1\le i\le l$. For all $\bar\vte_1, \bar\vte_2\in\mc M_i$,
\begin{equation*}
\inf \big\{ I_T (\pi|\bar\vte_1) : T>0 \,,\, \pi \in D([0,T], \mc M_+)
\,,\, \pi_T = \bar\vte_2 \}
\;=\; 0\;.
\end{equation*}
\end{lemma}

\begin{proof}
If $\mc M_i$ is a singleton, then the conclusion is clear.
Assume that $\mc M_i$ is not a singleton and fix $\bar\vte_1$,
$\bar\vte_2\in \mc M_i$ so that $\bar\vte_k(d\theta) = \bar\rho_k(\theta) d\theta$,
$k=1$, $2$, and $\bar\rho_2 (\theta) = \bar\rho_1(\theta+\theta_0)$ for some
$0<\theta_0<1$.  

Fix $a>0$ small, and let $\rho(t,\theta) = \bar\rho_1(\theta+ a t)$ so that
$\rho(0, \cdot ) = \bar\rho_1 (\cdot)$, $\rho(\theta_0/a, \cdot ) =
\bar\rho_2 (\cdot)$. Let $T=\theta_0/a$, $\pi_t(d\theta) = \rho(t,\theta) d\theta$. Since
$\partial_t \rho = a \nabla \rho$ and $(1/2) \Delta \rho + F(\rho)
= 0$, an integration by parts gives that for any smooth function
$G:[0,T] \times \bb T \to \bb R$,
\begin{align*}
J_{T,G}(\pi) \; =\; &- \, \int_0^T dt\,  \< a \rho_t , \nabla G_t \> - \frac 12 
\int_0^T dt\, \< \chi(\rho_t) , (\nabla G_t)^2 \> \\
& -\; \int_0^T dt\, \< B(\rho_t) , e^{G_t} -1 - G_t \> 
-\, \int_0^T dt\,  \< D(\rho_t) , e^{-G_t} -1 + G_t \> \;.
\end{align*}
The second line is negative, while the first one, by Young's
inequality and by Lemma \ref{lem1}, is less than or equal to
\begin{equation*}
\dfrac{a^2T}{2} \int_{\bb T} \frac {\bar\rho_1 (\theta)}{1- \bar\rho_1 (\theta)}\, d\theta
\;\le\; C_0 \, a\, \theta_0\;.
\end{equation*}
To complete the proof it remains to let $a\to 0$.
\end{proof}

\section{The static rate functional $W$}
\label{sec4b}

In this section, we present some properties of the quasi-potential
$W$. The main result asserts that the functional $W$, introduced in
\eqref{06}, is lower semicontinuous for the weak topology.

The first main result states that $W$ is continuous at each
measure $\bar\vte_i \in\mc M_i$ in the $L^2$-topology.  The
second one states that $W$ is lower semicontinuous in the weak
topology.

We start with an estimate of $V_i(\vte)$ for measures $\vte (d\theta)
= \gamma (\theta) d\theta$ whose density is close to $\mc M_i$ in the
$L^2$-topology.  This estimate together with Lemma \ref{bound4.1} will
allow us to prove that $V_i$ is continuous for the $L^2$-topology.

Let $\mathbb{D}$ be the space of measurable functions on $\T$ bounded
below by $0$ and bounded above by $1$, endowed with the $L^2$-topology:
\begin{equation*}
\mathbb{D} \;=\; \{\rho:\T\to[0,1] \;:\; 0\le\rho(\theta)
\le 1\;a.e.\  \theta\in\T\}\;.
\end{equation*}
For each $1\le i \le l$, let $\bb V_i:\bb D\to \bb [0,+\infty]$ be the
functional given by $\bb V_i(\rho) = V_i(\rho(\theta)d\theta)$.  Note
that the topology of $\M_{+,1}$ is the weak toplogy, while the one of
$\mathbb D$ is the $L^2$-toplogy.  Recall that we denote by $\mc H_1$
the Sobolev space of functions $G$ with generalized derivatives
$\nabla G$ in $L^{2}(\T)$. For each $h>0$ and each $\delta>0$, let
$\bb D^h_{\delta}$ be the subset of $\bb D$ consisting of those
profiles $\rho$ satisfying the following conditions:
\begin{enumerate}
\item[(A)] $\rho\in \mc H_1$ and $\displaystyle
  \int_{\T}(\nabla\rho(\theta))^2d\theta\le h$. 
\item[(B)] $\delta\le\rho(\theta)\le 1-\delta$ a.e. in $\T$.
\end{enumerate}

\begin{lemma}
\label{qpbound}
For each $1\le i\le l$, $h>0$, $\delta>0$ and an increasing
$C^1$-diffeomorphism $\alpha:[0,1]\to[0,1]$, there exist constants
$C_{1}=C_{1}(\delta, h)>0$ and $C_{2}=C_{2}(\delta, \a)>0$ such that for
any $\rho$ in $\mathbb D_{\delta}^{h}$ and $\bar\rho_i(\theta)d\theta$ in $\mc M_i$
\begin{equation*}
\bb V_i(\rho) \;\le\;
C_{1} \int_{0}^{1} \a(t)^{2} \, dt \;+\; 
C_{2} \, \|\rho-\bar\rho_{i}\|_{1} \;.
\end{equation*}
\end{lemma}

\begin{proof}
Fix $1\le i\le l$, $h>0$, $\delta>0$
and let $\a:[0,1]\to[0,1]$ be an increasing $C^1$-diffeomorphism. Let
$\rho\in\mathbb{D}_{\delta}^{h}$ and $\bar\rho_i(\theta)d\theta \in \mc M_i$.
Consider the path
$\pi_{t}^{\a}(d\theta)=\rho^{\a}(t,\theta)d\theta$ in $C([0,1], \M_{+})$ with density
given by $\rho_{t}^{\a} = (1-\a(t))\bar\rho_i+\a(t)\rho$.  It is clear
that $\pi^{\a}$ belongs to $D([0,1],\M_{+,1})$, and it follows from
condition (A) that $\mathcal{Q}_{1}(\pi^{\a})$ is finite.  From the
definition of $\rho^{\a}$ it follows that
$\nabla\rho_{t}^{\a}=\a(t)(\nabla\rho-\nabla\bar\rho_{i})+\nabla\bar\rho_{i}$
and that $\partial_{t}\rho_{t}^{\a}=\a'(t)(\rho-\bar\rho_{i})$.  Since
$(1/2) \Delta \bar\rho_{i} = -F(\bar\rho_{i})$, $J_{1,G}(\pi^{\a})$ can
be rewritten as
\begin{align}
\label{b4.1}
&J_{1,G}(\pi^\alpha)\;=\; \frac{1}{2} \int_{0}^{1}dt\  
\Big\{ \a(t) \lan(\nabla\rho-\nabla\bar\rho_{i}), \nabla G_{t}\ran
- \lan\chi(\rho_{t}^{\a}), (\nabla G_{t})^{2}\ran \Big\} \\
& \quad +\; \int_{0}^{1}dt \ \Big \lan
\big\{\a'(t)(\rho-\bar\rho_{i})+F(\bar\rho_{i}) \big\}
G_{t}-B(\rho_{t}^{\a})(e^{G_{t}}-1) - D(\rho_{t}^{\a}) (e^{-G_{t}}-1)
\Big\ran \;.
\notag
\end{align}

By Young's inequality, the first term on the right hand side of
\eqref{b4.1} is bounded by 
\begin{equation*}
\frac{1}{8} \int_{0}^{1}  \a(t)^2 \Big\lan
\frac{(\nabla\rho-\nabla\bar\rho_{i})^2}
{\chi(\rho_{t}^{\a})} \Big\ran \, dt 
\;\le\; C_{1}\int_{0}^{1}\a(t)^{2}\ dt
\end{equation*}
for some finite constant $C_{1}=C_{1}(\delta, h)$. To derive the last
inequality we used the fact that $\bar\rho_i$ is bounded away from $0$
and $1$ and conditions (A) and (B) on $\rho$.

To conclude the proof it is enough to show that the second term on the
right hand side of \eqref{b4.1} is bounded by
\begin{equation*}
C_{2}\|\rho-\bar\rho_{i}\|_{1} \,
\end{equation*}
for some constant $C_{2}=C_{2}(\delta, \a)$. 

Consider the function $\Phi:\R\times(0,1)\times\R\to\R$ defined by
\begin{equation*}
\Phi(H, \rho, G)\;=\; HG-B(\rho)(e^{G}-1)-D(\rho)(e^{-G}-1) \;.
\end{equation*}
If we set $H_{t}=\a'(t)(\rho-\bar\rho_{i})+F(\bar\rho_{i})$, 
it is clear that the second term on the right hand side of \eqref{b4.1}
can be expressed as
\begin{equation*}
\int_{0}^{1}\lan \Phi(H_{t}, \rho_{t}^{\a}, G_{t})\ran \ dt \;.
\end{equation*}
It follows from a straightforward computation that for any fixed $H
\in \R$ and $\rho \in (0,1)$, the function $\Phi(H, \rho, \cdot)$
reaches a maximum at
\begin{equation*}
G(H,\rho)\;=\; \log\left(\dfrac{H+\sqrt{H^{2}+4B(\rho)D(\rho)}}
{2B(\rho)}\right) \;.
\end{equation*}

From condition (B) and Lemma \ref{lem1}, there exists a constant
$c_{\delta}>0$ such that $c_{\delta}\le \rho^{\a} \le
1-c_{\delta}$. On the other hand, since $\Phi(H, \rho, 0)=0$ and since
$G(F(\rho), \rho)=0$, $\Phi(F(\rho), \rho, G(F(\rho), \rho))=0$ for any $\rho\in\R$.
Therefore, as $(H,\rho) \mapsto \Phi(H,\rho,G(H,\rho))$ is
a Lipschitz-continuous function on the interval $[-\Vert
\a'\Vert_\infty -\|F\|_{\infty}, \Vert \a'\Vert_\infty
+\|F\|_{\infty}]\times[c_{\delta}, 1-c_{\delta}]$,
\begin{align*}
\Phi(H_{t}, \rho_{t}^{\a}, G_{t})
&\;\le\; \Phi(H_{t}, \rho_{t}^{\a}, G(H_{t},\rho_{t}^{\a})) \\
&\;=\; \Phi(H_{t}, \rho_{t}^{\a}, G(H_{t},\rho_{t}^{\a}))
-\Phi(F(\rho_{t}^{\a}), \rho_{t}^{\a}, G(F(\rho_{t}^{\a}),\rho_{t}^{\a})) \\
&\;\le\; C_{2}|H_{t}-F(\rho_{t}^{\a})| \\
&\;\le\; C_{2} \big\{ \a'(t) + \|F'\|_{\infty}\a(t)\big\} |\rho-\bar\rho_{i}|
\end{align*}
for some finite constant $C_{2}=C_{2}(\delta, \a)$.  These bounds give
the desired conclusion.
\end{proof}

We are now in a position to prove that the functional $V_i$ is
continuous in the $L^2$-topology.

\begin{theorem}
\label{l2cont}
For each $1\le i\le l$, the function $\bb V_i$ is continuous
at $\bar\rho_i$ in $\bb D$.
\end{theorem}

\begin{proof}
Fix $1\le i \le l$, and let $\{\rho^{n} : n\ge 1\}$ be a sequence in $\bb D$
converging to $\bar\rho_{i}$. Denote by $\la^{n}$, $n\ge1$, the
weak solution to the equation \eqref{rdeq} with initial condition $\rho^{n}$.

By Lemma \ref{lem2}, there exists a constant $C_{0}>0$, independent of
$n$, such that
\begin{equation*}
\int_0^1 dt \int_{\T} |(\nabla \la_{t}^{n})(\theta)|^{2} d\theta \;\le\; C_{0} 
\end{equation*}
for all $n\ge 1$.  Fix $0<\zeta<1$. For each $n\ge1$, there exists
$\zeta\le T_{n}\le 2\zeta$ such that
\begin{equation*}
\int_{\T} |(\nabla \la_{T_{n}}^{n})(\theta)|^{2} d\theta \;\le\; C_{0}/\zeta \;.
\end{equation*}
Moreover, by Lemma \ref{lem1}, there exists a constant
$0<c_{\zeta}<1/2$, independent of $n$, such that $c_{\zeta}\le
\la_{T_{n}}^{n}(\theta) \le 1-c_{\zeta}$ for all $n\ge1$ and $\theta$ in $\T$.
Therefore, the density profiles $\la_{T_{n}}^{n}$, $n\ge1$, belong to
the set $\bb D_{c_{\zeta}}^{C_{0}/\zeta}$ introduced just above Lemma
\ref{qpbound}.

By definition \eqref{f01} of the functional $V_{i}$,
\begin{align*}
\bb V_{i}(\rho^{n}) \;\le\; 
\bb V_{i}(\la^{n}_{T_{n}}) \;+\; I_{T_{n}}(\pi^{n}) \;,
\end{align*}
where $\pi^{n}(t,d\theta)=\la^{n}(T_{n}-t ,\theta)d\theta$, $0\le t\le T_n$.
Therefore, it is enough to prove that
\begin{equation*}
\limsup_{n\to\infty} \bb V_{i}(\la^{n}_{T_{n}})
\;=\; \limsup_{\zeta\downarrow0}\limsup_{n\to\infty} I_{T_{n}}(\pi^{n}) \;=\; 0 \;.
\end{equation*}

Since $\la_{T_{n}}^{n}$ belongs to the set $\bb
D_{c_{\zeta}}^{C_{0}/\zeta}$, by Lemma \ref{qpbound}, there exist
constants $C_{1}=C_{1}(\zeta)>0$ and $C_{2}=C_{2}(\zeta,\a)>0$ such that
\begin{equation*}
\bb V_{i}(\la^{n}_{T_{n}}) \;\le\;
C_{1}\int_{0}^{1}\a(t)^{2}\, dt \;+\;
C_{2}\, \|\la_{T_{n}}^{n}-\bar\rho_{i}\|_{2} \;,
\end{equation*}
for any increasing $C^1$-diffeomorphism $\a:[0,1]\to[0,1]$.
By Lemma \ref{est3.1} and since $T_n\le 2\zeta\le 2$,
\begin{equation}
\label{b4.2}
\|\la_{T_{n}}^{n}-\bar\rho_{i}\|_{2} \;\le\; 
C_0\, \|\rho^{n}-\bar\rho_{i}\|_{2}
\end{equation}
for some finite constant $C_0>0$, independent of $n$, and whose value
may change from line to line.  Since $\rho^{n}$ converges to
$\bar\rho_{i}$ in $L^2$-topology,
\begin{equation*}
\limsup_{n\to\infty}\bb V_{i}(\la^{n}_{T_{n}})
\;\le\; C_{1}\inf_{\a} \Big\{\int_{0}^{1}\a(t)^{2}\, dt \Big\} \;=\; 0 \;.
\end{equation*}

It remains to prove that
\begin{equation}
\label{l4.1}
\limsup_{\zeta\downarrow0} \limsup_{n\to\infty} I_{T_{n}}(\pi^{n}) \;=\; 0 \;.
\end{equation}
By Lemma \ref{bound4.1}, by \eqref{b4.2}, and since $T_n\le 2\zeta$, 
\begin{equation*}
I_{T_{n}}(\pi^{n})
\;\le\; C_0\, \Big\{T_{n} \,+\, \|\la_{T_{n}}^{n}-\bar\rho_{i}\|_{2}
\,+\, \|\rho^{n}-\bar\rho_{i}\|_{2} \Big\}
\;\le\; C_0\, \big\{\zeta + \|\rho^{n}-\bar\rho_{i}\|_{2} \big\} \;.
\end{equation*}
To complete the proof of \eqref{l4.1} and the one of the lemma, it
remains to recall that $\rho^n \to \bar\rho_{i}$ in $L^2(\bb T)$.
\end{proof}

The proof of the next result is similar to the one of Proposition 4.9
in \cite{f}. The lemma asserts that the functional $V_i$ is uniformly
bounded in $\mc M_{+,1}$.

\begin{lemma}
\label{boundedness}
The function $W$ is finite if and only if $\vte$ belongs to
$\M_{+,1}$. Moreover, for all $1\le i\le l$,
\begin{equation*}
\sup_{\vte\in\M_{+,1}} V_i(\vte) \;<\; \infty \;,
\end{equation*}
so that $\sup_{\vte\in\M_{+,1}} W (\vte) < \infty$.
\end{lemma}

\begin{proof}
Fix $\vte\in \M_{+}$, and suppose that $W(\vte)<\infty$. By the
definition \eqref{06} of the functional $W$, there exists $1\le i\le
l$ such that $V_i(\vte)<\infty$. Hence, by \eqref{f01}, there exist
$\bar\rho\in \mc M_i$, $T<\infty$ and a trajectory $\pi_t(d\theta)$, $0\le
t\le T$, such that $\pi_T = \vte$, $I_T(\pi | \bar\rho)<\infty$. By
\eqref{i04} and \eqref{07}, $\pi \in D([0,T], \mc M_{+,1})$, proving
that $\vte = \pi_T$ belongs to $\mc M_{+,1}$, as claimed.

To prove the reciprocal assertion and the uniform bound, fix $1\le i
\le l$, $\vte\in\M_{+,1}$, $\vte (d\theta) = \rho(\theta) d\theta$, and denote by
$\la$ the weak solution of the equation \eqref{rdeq} with initial
condition $\rho$.  By Lemmata \ref{lem1} and \ref{lem2}, there exist
constants $0<a<1/2$ and $C_{0}>0$ such that $a \le \la_{t}(\theta) \le 1-a$
for all $t\ge 1$, $\theta\in \bb T$, and
\begin{equation*}
\int_0^2 dt \int_{\T} |(\nabla \la_{t})(\theta)|^{2} d\theta \;\le\; C_{0} \;.
\end{equation*}
In particular, there exists $1\le T\le 2$ such that $\la_{T}$ belongs
to the set $\bb D_{a}^{C_{0}}$.

By the definition \eqref{f01} of the functional $V_{i}$,
\begin{align*}
V_{i}(\rho) \;\le\; V_{i}(\la_{T}) \;+\; I_{T}(\pi) \;,
\end{align*}
where $\pi(t,d\theta)=\la(T-t ,\theta)d\theta$, $0\le t\le T$.  

As $\la_{T}$ belongs to the set $\bb D_{a}^{C_{0}}$, by Lemma
\ref{qpbound}, $V_{i}(\la_{T})\le C_1$ for some finite constant $C_1$
which depends only on $F$. On the other hand, by Lemma \ref{bound4.1}
and since $\Vert \gamma\Vert_1 \le 1$ for all density profile
$\gamma$, $I_{T}(\pi) \le C_2$, which completes the proof of the lemma
in view of the definition of the functional $W$.
\end{proof}

We now turn to the proof that the functionals $V_i$ are lower
semicontinuous for the weak topology. The idea of the proof is very
simple. Let $\vte^n$ be a sequence converging to $\vte$. Since, by
Lemma \ref{boundedness}, $V_i$ is finite, there exists a trajectory
$\pi^n_t$, $0\le t\le T_n$, such that $\pi^n_0 \in \mc M_i$,
$\pi^n_{T_n} = \vte^n$, $V_i(\vte^n) \le I_{T_n} (\pi^n) + 1/n \le
C$. We will now use the lower semicontinuity of $I_T$ and the fact
that the level sets are compact to conclude. If the sequence $T_n$ is
uniformly bounded, say by $T$, we may add a piece of length $T-T_n$ to
the trajectory $\pi^n$ letting it to stay at $\pi^n_0 \in \mc M_i$ in
the time interval $[0,T-T_n]$. In this way, we obtain a new sequence,
denoted by $\overline{\pi}^n_t$, of trajectories such that
$\overline{\pi}^n_0 \in \mc M_i$, $\overline{\pi}^n_{T} = \vte^n$,
$V_i(\vte^n) \le I_{T_n} (\pi^n) + 1/n = I_{T} (\overline{\pi}^n) +
1/n \le C$. Since the level sets are compact, we may extract a
converging subsequence.  Denote by $\pi$ the limit and observe that
$\pi_0 \in\mc M_i$, $\pi_T=\vte$. By the lower semicontinuity and by
definition of $V_i$, $V_i(\vte) \le I_{T} (\pi) \le \liminf_n
V_i(\vte^n)$, and we are done.

Of course, it might happen that the sequence $T_n$ is not bounded, and
this is the main difficulty. In this case, we will use Lemma
\ref{costs} to claim that the trajectory $\pi^n$ may not spend too much
time outside an $L^2$-neighborhood of a stationary solution. Hence,
all the proof consists in replacing the long intervals of time at which
the trajectory stays close to a stationary profile $\bar\rho_i$ by one
which remains only a time interval of length $2$. This is done by showing in
Lemma \ref{conn} below that it is possible to go from a neighborhood
of $\bar\rho_i$ to $\bar\rho_i$ in time $1$ by paying a small cost and
by using Theorem \ref{l2cont} to obtain a trajectory which goes from
$\bar\rho_i$ to a neighborhood of $\bar\rho_i$ in time $1$ by paying a
small cost.

\begin{lemma}
\label{conn}
Fix $1\le i\le l$ and $\bar\vte_i(d\theta)=\bar\rho_i(\theta)d\theta\in \mc M_i$. 
For any $\e>0$ there exists $\delta_{13}=\delta_{13}(\e)>0$ such that for any
$\vte(d\theta)=\ga(\theta)d\theta$ in $\mathcal{B}_{\delta_{13}}(\bar\vte_{i})$ there
exists a path $\pi(t,d\theta)=\rho(t,\theta)d\theta$ in $D([0,1], \M_{+})$ such that
$\pi_{0}=\vte$, $\pi_{1}=\bar\vte_{i}$, $I_{1}(\pi)\le\e$, and
$\pi_t \in \mathcal{B}_{\e}(\bar\vte_{i})$ for all $0\le t\le 1$.
\end{lemma}

\begin{proof}
Fix $1\le i\le l$, $\e>0$, and a density profile $\gamma: \bb T\to
\bb R$. Let $\rho(t,\theta)$ be the unique weak solution of \eqref{rdeq}
with initial condition $\gamma$, and let $\pi_t(d\theta)
=\rho(t,\theta) d\theta$. 

Fix $0<\zeta_1<\e$, to be chosen later, and set $T= 1/4$, $T'=1/2$. Let
$\zeta_2>0$ be the constant $\delta_{11}(\zeta_1, 1/4, 1/2)$ given by
Lemma \ref{rep} for $\bar\rho= \bar\rho_i$.  Assume that $\gamma\in
\mc B_{\zeta_2} (\bar\vte_{i})$.  According to Lemma \ref{rep},
$\rho_s(\theta)\, d\theta$ belongs to $\mc B_{\zeta_1} (\bar\vte_{i})$ for all $0\le
s\le 1/2$ and $\Vert \rho_s - \bar\rho_i\Vert_2 \le \zeta_1$ for all
$1/4\le s\le 1/2$.

By Proposition \ref{reg}, $\rho_t$ belongs to $C^\infty (\bb T)$ for
all $t>0$. By Lemma \ref{lem1}, there exists $a>0$, depending only on
$F$, such that $a\le \rho(s,\theta) \le 1-a$ for all $\theta\in\bb T$
and $1/4\le s\le 1/2$. On the other hand, by Lemma \ref{lem2}, there
exists a finite constant $C_0$, depending only on $F$, such that
\begin{equation*}
\int_0^{1/2} dt \int_{\bb T} |(\nabla \rho) (t,\theta)|^2 \, d\theta \;\le\;
C_0\;. 
\end{equation*}
In particular, there exists $T_1\in [1/4,1/2]$ such that
\begin{equation}
\label{04}
\int_{\bb T} |(\nabla \rho) (T_1,\theta)|^2 \, d\theta \;\le\; 4 C_0\;,
\end{equation}
so that $\rho_{T_1}$ belongs to $\bb D^{4C_0}_a$

Recall that $T_1\le 1/2$.  Let $\alpha: [0,1/2] \to [0,1]$ be an
increasing $C^1$-diffeomorphism, and define the trajectory
$\rho^\alpha_t$, $T_1\le t\le T_1+1/2$, by $\rho^\alpha(T_1+s, \theta) =
\alpha(s) \bar\rho_i + [1-\alpha(s)] \rho_{T_1}$, $0\le s \le$ 1/2.  By a similar
computation to the one presented in the proof of Lemma \ref{qpbound},
and since $\|\rho_{T_1}-\bar\rho_{i}\|_{1} \le
\|\rho_{T_1}-\bar\rho_{i}\|_{2} \le \zeta_1$
\begin{equation*}
I_{[T_1,T_1+1/2]} (\rho^\alpha) \;\le\; 
C_1  \int_{0}^{1/2} [1-\a(t)]^{2} \, dt \; \Vert \nabla \rho_{T_1} -
\nabla \bar\rho_{i}\|^2_{2}
\;+\; C_{2}(\alpha) \, \zeta_1 \;,
\end{equation*}
where $C_1$ is a finite constant depending only on $F$, and
$C_2(\alpha)$ one which also depends on $\alpha$.  In view of
\eqref{04},
\begin{equation*}
I_{[T_1,T_1+1/2]} (\rho^\alpha) \;\le\; 
C_3  \int_{0}^{1/2} [1-\a(t)]^{2} \, dt 
\;+\; C_{2}(\alpha) \, \zeta_1 \;,
\end{equation*}
for some finite constant $C_3$ independent of $\gamma$.

Choose an increasing $C^1$-diffeomorphism $\alpha: [0,1/2] \to [0,1]$
which turns the first term on the right hand side bounded by
$\varepsilon/2$. Note that this diffeomorphism does not depend on
$\gamma$. For this fixed $\alpha$, choose $\zeta_1$ small enough for
the second term to be less than or equal to $\varepsilon/2$. To
complete the proof of the first assertion of the lemma, juxtapose the
trajectories $\rho_t$, $0\le t\le T_1$, $\rho^\alpha_t$, $T_1\le t\le
T_1+1/2$, and the constant one $\bar\rho_i$, $T_1 +1/2\le t\le 1$.

We turn to the assertion that $\pi_t \in \mathcal{B}_{\varepsilon}
(\bar\vte_{i})$ for all $0\le t\le 1$. By the second paragraph of the
proof, and by definition of $T_1$, $\pi_t \in \mathcal{B}_{\zeta_1}
(\bar\vte_{i}) \subset \mathcal{B}_{\e} (\bar\vte_{i})$ for all
$0\le t\le T_1$. By definition of the trajectory $\rho^\alpha$,
$d(\pi_t, \bar \vte_{i}) \le d(\pi_{T_1}, \bar \vte_{i}) < \e$
for $T_1\le t \le T_1+1/2$. This completes the proof of the lemma
since $\pi_t = \bar \vte_{i}$ for $T_1+1/2 \le t\le 1$.
\end{proof}

We have now all the elements to prove the main result of the section. 

\begin{theorem}
\label{lsc}
For each $1\le i\le l$, $V_{i}$ is lower semicontinuous.
In particular, the rate function $W$ is also lower semicontinuous.
\end{theorem}

\begin{proof}
To keep notation simple, we prove the theorem in the case where all
solutions of \eqref{seeq} are constant in space, or equivalently, assume that
\begin{equation*}
\mc M_{\rm sol} \;=\; \{ \bar\vte_i(d\theta)=\bar\rho_id\theta : i=1,\cdots, l \} \;.
\end{equation*}
It is not difficult
to extend the argument to the general case by invoking Lemma
\ref{lem5}. 

Fix $1\le i \le l$, $q\in\R_{+}$, and let
\begin{equation*}
\mc V_i^{(q)} \;=\; \{ \vte\in\M_{+} : V_{i}(\vte)\le q \} \;.
\end{equation*}
By the proof of Lemma \ref{boundedness}, $\mc V_i^{(q)}\subset\M_{+,1}$.  We
claim that $\mc V_i^{(q)}$ is a closed subset of $\M_{+}$.  To see this, let
$\{\vte^{n}(d\theta)=\rho^{n}(\theta)d\theta : n\ge1\}$ be a sequence in $\mc V_i^{(q)}$
converging to some $\vte(d\theta)=\rho(\theta)d\theta$ in $\M_{+}$.

From \eqref{f01}, for each $n\ge1$, there exist $T_n>0$ and a path
$\pi^{n}$ in $C([0,T_{n}], \M_{+,1})$ such that
$\pi_{0}^{n}=\bar\vte_{i}$, $\pi_{T_{n}}^{n}=\vte^{n}$ and
\begin{equation}
\label{b4.3}
I_{T_{n}}(\pi^{n} | \bar\vte_{i} ) \;\le\; 
V_{i}(\vte^{n}) + 1/n \;\le\; q+1 \;.
\end{equation}

Assume first that the sequence $\{ T_{n} : n\ge1\}$ is bounded above
by some $T<\infty$. In this case, let $\widehat\pi^{n}$ be the trajectory
which remains at $\bar\vte_i$ in the time interval $[0, T-T_n]$ and
then follows the trajectory $\pi_n$:
\begin{equation*}
\widehat\pi^{n}_t \;=\;
\begin{cases}
\bar\vte_i & 0\le t\le T-T_n \;, \\
\pi_n(t-T+T_n) & T-T_n \le t\le T\;.
\end{cases}
\end{equation*}
Note that $\widehat\pi^{n}_T = \vte^{n}$ for all $n\ge 1$, and that
$I_{T}(\widehat \pi^{n} | \bar\vte_{i}) = I_{T_{n}}(\pi^{n} |
\bar\vte_{i}) \le q+1/n$ for all $n\ge 1$.

By Theorem \ref{dldp}, the functional $I_{T}(\cdot | \bar\vte_{i})$
has compact level sets and is lower semicontinuous. There exists, in
particular, a subsequence $n_j$ and a trajectory $\pi \in D([0,T], \mc
M_+)$ such that $\widehat \pi^{n_j} \to \pi$, $I_{T}(\pi | \bar\vte_{i})\le
q$. Since $\pi_T = \lim_j \widehat \pi^{n_j}_{T} = \lim_j \pi^{n_j}_{T_{n_j}} = \lim_j
\vte^{n_j} = \vte$, by definition of $V_i$, $V_i(\vte)\le q$, which
completes the proof of the theorem in the case where the sequence
$T_n$ is bounded.

Suppose now that the sequence $\{ T_{n} : n\ge1\}$ is unbounded.  In
view of \eqref{b4.3}, by Lemma \ref{costs}, the path $\pi^n$ may not
remain too long outside an $L^2$-neighborhood of one of the stationary
profiles. We use this observation, Lemma \ref{conn} and Theorem
\ref{l2cont}, to construct from $\pi^n$ a new path on a bounded time
interval by replacing the long intervals of time in which $\pi^n$
remained close to a stationary profile by a path defined in a time
interval of length $2$ which connects the entrance time in a
neighborhood of a stationary profile $\bar\rho_i$ to $\bar\rho_i$ and
from this profile to the exit time of the neighborhood. The details
are given below.

Fix $\e>0$. By Theorem \ref{l2cont}, there exists $\zeta_1$, such
that $V_j(\pi) \le \varepsilon$ if $\pi(d\theta) = \rho(\theta) d\theta$ and $\Vert
\rho-\bar\rho_j\Vert_2\le \zeta_1$ for any $1\le j\le l$. Let
$\zeta_2$ be the constant $\delta_{13}(\varepsilon)$ given by Lemma
\ref{conn} and set $\zeta = \min\{\zeta_1, \zeta_2\}$.

Let $\mc L_{j}$, $1 \le j \le l$, be the closed $L^2$-neighborhood of
$\bar\rho_{j}$:
$\mc L_{j} = \{\vte(d\theta)=\rho(\theta)d\theta :
\rho\in\mathbb{B}_{\zeta} [\bar\rho_{j}]\}$,
where $\mathbb{B}_{\zeta} [\bar\rho_{j}]$ represents the closure of
$\mathbb{B}_{\zeta} (\bar\rho_{j})$, and let
$\mc L = \cup_{1\le j\le l} \mc L_j$.  Assume that $\zeta$ is
sufficiently small so that $\{\mc L_j\}_{j=1}^l$ are mutually
disjoint.  Note that $\mc L_{j}$ is a closed subset of $\M_{+,1}$.  We
define a sequence of entrances and exit times associated to the sets
$\mc L_j$. Recall that $\pi^n_0= \bar\vte_i$, and set $\tau^n_1
=0$. Let $\sigma^n_1$ be the last exit time from $\mc L_i$:
\begin{equation*}
\sigma^n_1 \;=\; \sup \big\{t\le T_n : \pi^n_t \in \mc L_i
\big\}\;.
\end{equation*}
Note that the path $\pi^n_t$ may visit several neighborhoods $\mc
L_j$, $j\not = i$, in the time interval $[0, \sigma^n_1]$, and that it
does not return to $\mc L_i$ after $\sigma^n_1$.  Suppose that
$\tau^n_k$, $\sigma^n_k$, $1\le k < p$, have already been
introduced. Define
\begin{equation*}
\tau^n_p \;:=\; 
\inf \big\{\sigma^n_{p-1}\le t\le T_{n}: \pi^n_{t}\in\mc L \big\} \;,
\quad
\sigma^n_p \;:=\; 
\sup \big\{t\le T_{n}: \pi^n_{t}\in \mc L_{\mf j(p)} \big\} \;,
\end{equation*}
where $\mf j(p)$ is the index of the neighborhood visited at time
$\tau^n_p$: $\mf j(p)= a$ if $\pi^n_{\tau^n_p} \in \mc L_a$. By
convention, if $\pi^n_{t}\not \in\mc L$ for all $\sigma^n_{p-1}\le
t\le T_{n}$, we set $\tau^n_{p'}$ and $\sigma^n_{p'}$ to be $\infty$
for all $p'\ge p$ and we do not define $\mf j(p)$. Note that $\mf j(1)
= i$ and that $\mf j(p) \not = \mf j(q)$ if $q\not = p$.

Denote by $S_n$ the set of neighborhoods visited by $\pi^n$, $S_n = \{
\mf j(p) : \tau_p(\pi^{n})<\infty\} = : \{\mf j(1) , \dots, \mf
j(b)\}$. By the choice of $\zeta$ and by Lemma \ref{conn}, there
exist paths $\pi^{1,+}$, $\pi^{m,-}$, $\pi^{m,+}$, $2\le m< b$,
$\pi^{b,-}$ such that
\begin{equation*}
\pi_{0}^{k,-} \;=\; \pi^n(\tau_{\mf j(k)}) \;,\quad 
\pi_{1}^{k,-} \;=\; \bar\vte_{\mf j(k)} \;, \quad
\pi_{0}^{k,+} \;=\; \bar\vte_{\mf j (k)}\;, \quad 
\pi_{1}^{k,+}=\pi^n(\sigma_{\mf j(k)})\;,
\end{equation*}
and such that $I_{1}(\pi_{1}^{1,+}) \le 2\varepsilon$,
$I_{1}(\pi_{1}^{m,-}) + I_{1}(\pi_{1}^{m,+})\le 3\e$, $2\le m< b$,
$I_{1}(\pi_{1}^{b,-}) \le \varepsilon$. Here and below, to keep notation
simple, we denote sometimes $\pi^n_t$ by $\pi^n(t)$. If
$\sigma^n_b \le T_n$ there exists also a path $\pi^{b,+}$, such that
\begin{equation*}
\pi_{0}^{b,+} \;=\; \bar\vte_{\mf j (b)}\;, \quad 
\pi_{1}^{b,+}=\pi^n(\sigma_{\mf j(b)})\;,
\end{equation*}
and such that $I_{1}(\pi_{1}^{b,+}) \le 2\varepsilon$.

We construct below a path $\widetilde\pi^{n}$ from the previous paths and
from $\pi^n$ under the assumption that $\sigma^n_b \le T_n$. The
construction can be easily adapted to the cases
$\sigma^n_b=+\infty$. For $1\le c<b$, let
\begin{equation*}
R_{c} \;=\; 2c \;+\; \sum_{a=1}^{c} (\tau_{\mf j(a+1)} - \sigma_{\mf j(a)})\;.
\end{equation*}
This sequence represents the times at which the path $\widetilde\pi^{n}$
visits the measures $\bar\vte_{\mf j(c+1)}$. Set
\begin{equation*}
\widetilde T_n \;=\; R_{b-1} \;+\; 1 \;+\; T_n \;-\; \sigma_{\mf j(b)} \;,
\end{equation*}
and define the path $\widetilde\pi^{n}$ in $C([0,\widetilde{T}_{n}],
\M_{+,1})$ as follows:
\begin{equation*}
\widetilde\pi^{n} (t) \;=\; 
\begin{cases}
\pi^{1,+} (t) & 0\le t\le 1\;,\\
\pi^n(\sigma_{\mf j(1)} + t -1) & 1\le t\le 1+ \tau_{\mf j(2)} -
  \sigma_{\mf j(1)} \;, \\
\pi^{2,-} (t - 1 - \tau_{\mf j(2)} + \sigma_{\mf j(1)}) & 1+ \tau_{\mf
  j(2)} - \sigma_{\mf j(1)}\le t\le R_1 \;,\\
\pi^{2,+} (t - R_1) & 
 R_1 \le t\le  R_1 +1 \;, \\
\cdots & \cdots \\
\pi^{b,+} (t-R_{b-1}) & R_{b-1}\le t\le R_{b-1} + 1\;, \\
\pi^n(\sigma_{\mf j(b)} + t - R_{b-1} - 1) & R_{b-1} + 1\le t \le
R_{b-1} + 1 + T_n - \sigma_{\mf j(b)}\;.
\end{cases}
\end{equation*}

From the definition of $\widetilde\pi^{n}$ it is clear that
$\widetilde\pi_{0}^{n}=\bar\vte_{i}$, $\widetilde\pi^{n}(\widetilde{T}_{n})
=\pi^{n}(T_{n}) = \vte$ and that
\begin{equation}
\label{b4.3b}
I_{\widetilde{T}_{n}}(\widetilde\pi^{n}) \;\le\; I_{T_{n}}(\pi^{n}) \;+\;
3b\e \;\le\; I_{T_{n}}(\pi^{n}) \;+\;  3l\e\;.
\end{equation}
In particular, by \eqref{b4.3} and \eqref{b4.3b}, $I_{\widetilde{T}_{n}}(\widetilde\pi^{n})$ is
uniformly bounded. The time spent by $\widetilde\pi^{n}$ in $\mc L^{c}$ is
at least $\widetilde{T}_{n}-2l$. Therefore, by Lemma \ref{costs} and by
the previous uniform bound on $I_{\widetilde{T}_{n}}(\widetilde\pi^{n})$, the
sequence $\widetilde{T}_{n}$ is uniformly bounded. At this point, we may
repeat the arguments presented in the first part of the proof, which
are solely based on a uniform bound for the sequence
$I_{\widetilde{T}_{n}}(\widetilde\pi^{n})$, provided by \eqref{b4.3} and
\eqref{b4.3b}, and on a uniform bound of the sequence $\widetilde T_n$, to
conclude.
\end{proof}

We conclude this section with a result needed in the next one.
In Lemma \ref{conn} we constructed paths from $\mc B_\delta(\mc
M_i)$ to $\mc M_i$ whose costs are small. In the next result, we prove
a partial converse statement by showing that there are measures $\vte$
at distance $\delta$ from $\mc M_i$ for which there exist paths from
$\mc M_i$ to $\vte$ whose cost is small.

\begin{lemma}
\label{lem8}
For every $1\le i\le l$ and $\varepsilon>0$, there exists
$\delta_{14}=\delta_{14}(\varepsilon)>0$ such that for all $0\le
\delta<\delta_{14}$, there exists a measure $\vte \in \mc M_{+,1}$
such that $d(\vte, \mc M_i) = \delta$ and $V_i(\vte) \le
\varepsilon$.
\end{lemma}

\begin{proof}
Fix $1\le i\le l$, $\bar\vte_i\in\mc M_i$, and
$\varepsilon>0$. Recall the definition of the constant $c_1$
introduced in \eqref{11}, and let $0<a\le c_1/2$. Let $\pi(t,d\theta)$, $0\le
t\le 1$, be the trajectory $\pi(t,d\theta) = \rho(t,\theta) d\theta$, $\rho(t,\theta) =
\bar\rho_i(\theta) + at$. On the one hand, by definition \eqref{01} of the
distance, $d(\mc M_i, \pi_1) = d(\bar\vte_i, \pi_1) = a$. On the other
hand, for every smooth function $G:[0,1]\times \bb T \to \bb R$, since
$(1/2)\De\rho_t = (1/2)\De\bar\rho_i = -F(\bar\rho_i)$,
\begin{align*}
J_{1,G}(\pi) \; \le \; \int_0^1 dt\, \Big\{
\big\lan \big\{ a + [F(\bar\rho_i) - F(\rho_t)] \big\} \, G_t \big\ran 
\;& -\; \lan B(\rho_t), e^{G_t}-1 - G_t\ran \Big\}\;.
\end{align*}
Note that we omitted in the previous expression the terms
$\chi(\rho_t) [\nabla G_t]^2$ and $D(\rho_t) [e^{-G_t}-1 + G_t]$ which
are positive.  By definition of $\rho_t$, $F(\bar\rho_i) - F(\rho_t)$
is absolutely bounded by $C_F a$, where $C_F$ stands for the Lipschitz
constant of the function $F$.  By definition, $a\le c_1/2$ and by
\eqref{11}, $c_1 \le \bar\rho_i \le 1-c_1$. Thus, $c_1/2 \le \rho_t
\le 1-c_1/2 $. Let $b = \inf \{ B(x) : c_1/2 \le x \le 1-c_1/2\}>0$,
so that
\begin{equation*}
J_{1,G}(\pi) \; \le \; \int_0^1 dt\, \Big\lan (1+h_t)\, a \, G_t 
- b\, [ e^{G_t}-1 - G_t] \Big\ran \;,
\end{equation*}
where $h_t$ is absolutely bounded by $C_F$. Assume that $0\le a<a_1$
where $a_1$ is chosen so that $(C_F-1)a_1 < b$.  The right hand side
of the previous expression is bounded by $\psi_b((1\pm C_F) a)$,
uniformly in $G$, where $\psi_b(x) := (x+b) \log [1 + (x/b)] -
x$.
Therefore, $V_i(\vte) \le I_1(\pi) \le \psi_b((1\pm C_F) a)$. Since
$\psi_b (0)=0$, there exists $a_0$ such that
$\psi_b((1\pm C_F) a) \le \varepsilon$, for any $0 \le a < a_0$. The
assertion of the lemma holds provided we choose
$\delta_{14} = a_0 \wedge a_1$ and $\vte = \bar\vte_i + a d\theta$ for
$0\le a < \delta_{14}$.
\end{proof}

\begin{remark}
\label{rem2}
Actually, we proved the existence of trajectory $\pi_t$, $0\le t\le
1$, such that $I_1(\pi)\le \varepsilon$, $\pi_0 = \bar\vte_i \in \mc
M_{i}$, $d(\pi_1, \mc M_i) = \delta$.
\end{remark}

\section{The static large deviations principle}
\label{sec5}

We prove in this section Theorem \ref{sldp}. As we said before, the
proof is based on a representation of the stationary state of the
reaction-diffusion model in terms of the invariant probability measure
of a discrete-time Markov chain. In the first part of this section, we
introduce the discrete-time Markov chain and we prove in Proposition
\ref{main5.1} sharp upper and lower bounds for its invariant
probability measure.

For any $0<\beta_0<\beta_1$, let $B_i$ be the open neighborhoods, and
$\Gamma_i$ be the closed neighborhoods given by
\begin{gather*}
B \;=\; \bigcup_{i=1}^l B_i\;,
\quad \text{where}\quad B_i \;:=\; \mc B_{\beta_0}(\mc M_i) \;:=\;
\{\vte\in\M_{+}:\, \inf_{\bar\vte\in\M_{i}} 
d(\vte, \bar\vte) < \beta_0\} \;. \\
\Ga \;=\; \bigcup_{i=1}^l \Ga_i\;,
\quad \text{where}\quad
\Ga_i \;=\; \{\vte\in\M_{+}:\, 
\beta_1\le \inf_{\bar\vte\in\M_{i}} d(\vte, \bar\vte)\le 2\beta_1\} \;.
\end{gather*}
To stress the dependence of $B$ and $\Gamma$ on $\beta_0$, $\beta_1$
we sometimes denote $B$, $\Gamma$ by $B(\beta_0)$, $\Gamma(\beta_1)$,
respectively. 

For $N\ge 1$ and a subset $A$ of $\mc M_{+}$, let $A^N =
(\pi^N)^{-1}(A)$ and let $H_A^N:D(\bb R_+, X_N)\to[0,+\infty]$ be the
hitting time of $A^N$:
\begin{equation*}
H^N_A \;=\; \inf\big\{t\ge 0: \eta_t\in A^N\big\} \;.
\end{equation*}

The first result of the section states that the reaction-diffusion
model reaches the set $B$ in finite time with high probability.

\begin{lemma}
\label{lem4}
For every $\delta>0$, there exist $T_0$, $C_0$, $N_0 >0$, depending on
$\delta$, such that
\begin{equation*}
\sup_{\eta\in X_N} \bb P_{\eta}\left[H_{B(\delta)}^N \ge kT_0\right]
\;\le\; \exp\left\{-k C_0 N \right\} 
\end{equation*}
for all $N\ge N_0$ and all $k\ge 1$.
\end{lemma}

\begin{proof}
Fix $\delta>0$. By Corollary \ref{costw}, there exist $T_0>0$ and $C_0>0$,
which depend on $\delta$, such that
\begin{equation*}
\inf_{\pi\in\mc D_{T_0,\delta}} I_{T_0}(\pi) \;>\; C_0 \;,
\end{equation*}
where $\mc D_{T_0,\delta} = D([0,T_0],\mc M_{+}\backslash B(\delta))$.  For
each integer $N\ge 1$, denote by $\eta^N$ a configuration in $X_N$ such
that
\begin{equation*}
\bb P_{\eta^N}\left[H_B^N \ge T_0\right]
\;=\;\sup_{\eta\in X_N} 
\bb P_{\eta}\left[H_B^N \ge T_0\right] \;.
\end{equation*}

By the compactness of $\M_{+}$, every subsequence of $\pi^N(\eta^N)$
contains a sub-subsequence converging to some $\vte$ in $\M_{+}$.
Moreover, since each configuration in $X_N$ has at most one particle
per site, any limit point $\vte$ belongs to $\M_{+,1}$. From this
observation and since $\mc D_{T_0,\delta}$ is a closed subset of
$D([0,T_0],\M_{+})$, by the dynamical large deviations upper bound,
there exists a measure $\vte(d\theta) = \gamma(\theta)d\theta$ in $\M_{+,1}$ such
that
\begin{align*}
\limsup_{N\to\infty}\frac{1}{N}\log \bb P_{\eta^N}\left[H_B^N\ge  T_0\right]
\; & \le\; \limsup_{N\to\infty}\frac{1}{N}\log {Q}_{T_0, \eta^N}(\mc
D_{T_0,\delta})  \\
\; &\le\; -\inf_{\pi\in\mc D_{T_0,\delta}} I_{T_0}(\pi|\gamma) 
\;<\; -C_0 \;.
\end{align*}
In particular, there exists $N_0\ge 1$ such that for every integer
$N\ge N_0$,
\begin{align*}
\sup_{\eta\in X_N} \bb P_{\eta}\left[H^N_B\ge T_0\right] 
\;\le\; \exp\{-C_0N\} \;.
\end{align*}

To complete the proof, we proceed by induction, applying the strong Markov
property.  Suppose that the statement of the lemma is true for all
integers $j<k$.  Let $N\ge N_0$ and let $\hat{\eta}$ be a
configuration in $X_N$.  By the strong Markov property,
\begin{align*}
\bb P_{\hat{\eta}}\left[H^N_B\ge kT_0\right] \;=\;
& \bb E_{\hat{\eta}}\left[{\bf 1} \left\{H^N_B\ge T_0\right\}\, 
\bb P_{\eta_{ _{T_0}}}\left[H^N_B\ge(k-1)T_0\right]\right]\\
\;\le\; & \bb P_{\hat{\eta}}\left[H_B^N\ge T_0\right]\;
\sup_{\eta\in X_N} \bb P_{\eta}\left[H^N_B\ge(k-1)T_0\right] \\
\;\le\; & \exp\left\{-k\, C_0N\right\} \;,
\end{align*}
which completes the proof.
\end{proof}

\begin{corollary}
\label{lem7}
For every $\delta>0$, 
\begin{equation*}
\limsup_{N\to\infty}\frac{1}{N}
\log\sup_{\eta\in X_N} \bb E_{\eta}\left(H^N_{B(\delta)}\right)
\;\le\; 0 \;.
\end{equation*}
\end{corollary}

\begin{proof}
Fix $\delta>0$ and denote by $T_0$, $C_0$, $N_0\ge 1$ the constants
provided by Lemma \ref{lem4}. For every integer $N\ge N_0$ and for
every configuration $\eta$ in $X_N$, 
\begin{equation*}
\bb E_{\eta}\left(H_B^N\right)
\;\le\; T_0\sum_{k=0}^{\infty}\bb P_{\eta}\left(H_B^N\ge kT_0\right)
\;\le\; T_0\sum_{k=0}^{\infty}\exp\left\{-kC_0N\right\}
\;\le\; \frac{T_0}{1-e^{-C_0N}} \;,
\end{equation*}
which proves the corollary.
\end{proof}

We have now all elements to introduce the discrete-time Markov
chain. Let $\partial B^N$ (which depends on $\beta_0$)\footnote{Note
  that for $N$ large the set $\partial B^N$ does not depend on
  $\beta_1$.} be the set of configurations $\eta$ in $X_N$ for which
there exists a finite sequence of configurations
$\{\eta^i: 0\le i\le k\}$ in $X_N$ with $\eta^0$ in $\Gamma^N$,
$\eta^k = \eta$, and such that
\begin{enumerate}
\item[(a)] For every $1\le i\le k$, the configuration $\eta^{i}$ can
  be obtained from $\eta^{i-1}$ by a jump of the dynamics (either from
  the stirring mechanism or from the non-conservative spin flip
  dynamics).
\item[ (b)] The unique configuration of the sequence which can enter
  into $B^N$ after a jump is $\eta^k$.
\end{enumerate}
We similarly define the set $\partial B_{i}^{N}$,
$1\le i\le l$.  It is clear that for $N$ large
enough and $\beta_1$ small enough,
\begin{equation*}
\partial B^N \;=\; \bigcup_{i=1}^l\partial B^N_i \;.
\end{equation*} 

Let $\tau = \tau^N: D(\bb R_+, X_N)\to [0,\infty]$ be the stopping
time given by
\begin{equation}
\label{10}
\tau\;=\;\inf\left\{t> 0: \hbox{ there exist }
s < t \hbox{ such that } \eta_s\in \Ga^N \hbox{ and }
\eta_t\in\partial B^N\right\} \;.
\end{equation}
Set $\tau_1:=\tau$. We recursively define the sequence of stopping times $\{\tau_k : k\ge 1\}$ by
\begin{align*}
\tau_k\;=\;\inf\left\{\tau_k> 0: \hbox{ there exist }s < t \hbox{ such that }\tau_{k-1} < s, \eta_s\in \Ga^N \hbox{ and } \eta_t\in\partial B^N\right\} \;.
\end{align*}
This sequence generates a discrete-time
Markov chain $\xi_k$ on $\partial B^N$ by setting $\xi_k =
\eta_{\tau_k}$. The arguments presented before Lemma 5.1 in \cite{f}
show that this chain is irreducible. Denote by $\nu^{N}$ its unique
invariant probability measure.

Define $\tilde v_{ij}$ by
\begin{gather*}
\tilde v_{ij} \;=\; \inf \{I_T(\pi|\bar\rho): T> 0 \,,\,
\bar\rho(\theta)d\theta\in\M_{i} \,,\, \pi \in \mc A_T \,,\,
\pi_T \in \mc M_{j} \}\;, \\
\text{where}\quad 
\mc A_T \;=\; \big\{ \pi\in C([0,T],\mc M_{+}) : 
\pi_{t}\notin \M_{\rm sol} \text{ for all $0<t<T$ } \big\}\;.
\end{gather*}

The jumps of the Markov chain $\{\xi_k : k\ge 1\}$ correspond to paths
from $\partial B^N_i$ to $\partial B^N_j$ which do not visit other
boundaries. They are thus related to the dynamical large deviations
principle.  In Lemmata \ref{blb} and \ref{bub} we estimate the
probability of a jump from $\partial B^N_i$ to $\partial B^N_j$ for
$j\not = i$. These estimates provide sharp bounds for the invariant
probability measure of the chain $\xi_k$, alluded to at the beginning of
this section.

\begin{lemma}
\label{blb}
For every $1\le i \not = j\le l$, $\e>0$ and $0<\beta_1< (1/4)
\min_{a\not = b} d(\mc M_a, \mc M_b)$, there exists $0<\delta_{15}<
\beta_1$ such that for all $\beta_0<\delta_{15}$ 
\begin{equation*}
\liminf_{N\to\infty}\frac{1}{N}\log\inf_{\eta\in\partial B_i^N}\bb P_{\eta}
(\eta_{\tau}\in\partial B_j^N) \;\ge\; -\, \tilde v_{ij} \,-\, \e \;.
\end{equation*}
\end{lemma}

\begin{proof}
Fix $\varepsilon>0$, $1\le i \not = j\le l$, and $0<\beta_1< (1/4)
\min_{a\not = b} d(\mc M_a, \mc M_b)$. By definition of $\tilde v_{ij}
$, there exist $T>0$, $\bar\vte_i (d\theta) = \bar\rho_i (\theta) d\theta \in \mc
M_i$, $\bar\vte_j \in \mc M_j$, and $\pi \in \mc A_T$ such that $\pi_0
= \bar\vte_i$, $\pi_T = \bar\vte_j$, $I_T(\pi | \bar\rho_i) \le \tilde
v_{ij} + \varepsilon$. Note that $\pi_t\not\in \mc M_i$ for all
$0<t\le T$ because $\pi \in \mc A_T$.

Let $t_0$ be the first time the path $\pi_t$ is at distance $\beta_1$
from $\mc M_i$: $t_0 = \min\{ t\ge 0 : d(\pi_t, \mc M_i) \ge
\beta_1\}$. Since $\pi$ belongs to $\mc A_T$ and $\pi_T\in \mc M_j$,
$\zeta_2 = \inf_{t_0\le t\le T} d(\pi_t, \cup_{k\not =j} \mc M_k)
>0$. Let $\zeta_3$ be the constant $\delta_{13}(\min\{\beta_1,
\varepsilon\})$ given by Lemma \ref{conn}.  Set $\delta_{15}= (1/2)
\min\{\beta_1,\zeta_2, \zeta_3\}$ and fix $\beta_0<\delta_{15}$. 

For each integer $N>0$, let $\eta^N$ be a configuration in $\partial
B_i^N$ such that
\begin{equation}
\label{08}
\bb P_{\eta^N}(\eta_{\tau}\in\partial B_{j}^{N})
\;=\; \inf_{\eta\in \partial B_i^N}\bb P_{\eta}(\eta_{\tau}
\in\partial B_{j}^{N}) \;.
\end{equation}
Recall that every subsequence of $\pi^N(\eta^N)$ contains a
sub-subsequence converging in $\M_{+}$ to some measure $\vte$ which
belongs to $\M_{+,1}$. We may therefore assume that $\pi^N(\eta^N)$
converges to $\vte(d\theta)=\gamma(\theta)d\theta$ which belongs to the closure of
$B_{i}$: $\vte \in \overline{\mc B_{\beta_0}(\bar\vte'_{i})}$, for some
$\bar\vte'_{i} \in \mc M_i$.

Since $\beta_0<\zeta_3$, we may apply Lemma \ref{conn} to connect
$\vte$ to $\bar\vte_i'\in \mc M_i$.  Denote by $\pi'_t$, $0\le t\le
1$, the path given by Lemma \ref{conn} and such that $\pi'_0 = \vte$,
$\pi'_1 = \bar\vte_i'\in \mc M_i$, $I_1(\pi'| \gamma) \le \varepsilon$
and $\pi'_t \in \mc B_{\beta_1}(\bar\vte_i')$, $0\le t\le 1$.  We may
apply Lemma \ref{lem5} to connect $\bar\vte_i'$ to $\bar\vte_i$, the
initial point of the path introduced in the first paragraph of the
proof. By Lemma \ref{lem5}, there exists $T'$ and a path $\pi'' \in
C([0,T'], \mc M_+)$ such that $\pi''_0 = \bar\vte_i'$, $\pi''_{T'} =
\bar\vte_i$ and $I_{T'}(\pi''| \bar\vte_i') \le \varepsilon$.

Concatenate the paths $\pi'$, $\pi''$ and $\pi$ to obtain a path
$\widetilde \pi$ in $C([0, T+T'+1], \mc M_+)$ such that $\widetilde \pi_0 =
\vte$, $\widetilde \pi_{T+T'+1} \in \mc M_j$, $I_{T+T'+1}(\pi''| \vte) \le
\tilde v_{ij} + 3\varepsilon$. Moreover, $d(\pi_t, \mc M_i) <
\beta_1$ for $0\le t\le 1+T'$. In particular, $\widetilde \pi_t$ reaches
a distance $\beta_1$ from $\mc M_i$ for the first time at $t=
1+T'+t_0$, and $ \inf_{t_0\le t\le T} d(\widetilde\pi_{1+T'+t},
\cup_{k\not =j} \mc M_k) =\zeta_2> \beta_0$.

Denote by $\mc N$ a $\beta_0/2$-neighborhood in $D([0,T+T'+1], \mc
M_+)$ of the path $\widetilde\pi$. It follows from the last
observation of the previous paragraph and from the fact that $\widetilde
\pi_{T+T'+1}\in\mc M_j$ that $\mc N \subset \{\eta_{\tau}\in\partial
B_j^N\}$. Since $\mc N$ is an open set, by the lower bound of the
dynamical large deviations principle, and by \eqref{08}
\begin{align*}
& \liminf_{N\to\infty}\frac{1}{N} \inf_{\eta\in \partial B^N_i}
\log \bb P_{\eta}(\eta_{\tau}\in\partial B_{j}^{N})
\;=\; \liminf_{N\to\infty}\frac{1}{N} \log 
\bb P_{\eta^N}(\eta_{\tau}\in\partial B_{j}^{N}) \\
&\quad \;\ge\;
\liminf_{N\to\infty}\frac{1}{N}\log 
\bb P_{\eta^N}(\mc N)
\;\ge\; - \inf_{\bar \pi \in\mc N} I_{T+T'+1}(\bar \pi|\ga)
\;\ge\; - \, I_{T+T'+1}(\widetilde \pi|\ga)\;.
\end{align*}
This completes the proof of the lemma because
$I_{T+T'+1}(\widetilde\pi| \vte) \le \tilde v_{ij} + 3\varepsilon$.
\end{proof}

The proof of the upper bound for $\bb P_{\eta} (\eta_{\tau}\in\partial
B_j^N)$, $\eta \in\partial B_i^N$, requires a lower bound for the
dynamical large deviations rate functional.  For $T>0$ and $\zeta>0$,
let $\mc C_j = \mc C_j(T,\zeta)$ be the closed subset of $D([0,T],
\M_{+})$ consisting of all paths $\pi$ for which there exists some
time $t \in [0, T]$ such that $\pi(t)$ belongs to $\Gamma_j(\zeta)$ or
$\pi(t-)$ belongs to $\Gamma_j(\zeta)$.

\begin{lemma}
\label{lem6}
For every $1\le i\not = j\le l$, $\varepsilon >0$, there exist
$\delta_{16} = \delta_{16}(\varepsilon)>0$ and
$T=T(\varepsilon)>0$ such that for all $\delta'<\delta_{16}$,
$T'\ge T$, $\gamma(\theta)d\theta\in \Gamma_i(\delta')$,
\begin{equation*}
\inf_{\pi\in \mc C_j(T',\delta')} I_{T'}(\pi|\gamma)
\;\ge\; v_{ij} \;-\; \e \;.
\end{equation*}
\end{lemma}

\begin{proof}
Fix $1\le i\not = j\le l$ and assume that the assertion of the lemma
is false. In that case, there exists $\varepsilon>0$ such that for
for every $\zeta>0$ and $T>0$, there exist $\zeta'<\zeta$, $T'\ge T$,
$\gamma(\theta)d\theta\in\Gamma_i(\zeta')$ and $\pi\in \mc C_j(T',\zeta')$ with
\begin{equation*}
I_{T'}(\pi|\gamma) \; < \; v_{ij} \;-\; \e/2 \;.
\end{equation*}

In particular, taking the sequences $\zeta_n = 1/n$, $T_n=1$, $n\ge
1$, there exist $\zeta'_n <1/n$, $T'_n\ge 1$, $\gamma_n\in
\Gamma_i(\zeta'_n)$ and $\pi_n\in \mc C_j(T'_n, \zeta'_n) \cap
C([0,T'_n],\M_{+,1})$ with
\begin{equation}
\label{contra}
I_{T'_n}(\pi^n|\gamma_n)\;<\; v_{ij} \,-\, \e/2 \;.
\end{equation}
Since $\pi^n$ belongs to $\mc C_{j}(T'_n, \zeta'_n) \cap
C([0,T'_n],\mc M_{+,1})$, there exists $0<\widetilde T_n\le T'_n$ such
that $\pi^n(\widetilde T_n) \in \{\vte\in \M_{+,1}:\, \zeta'_n \le
\inf_{\bar\vte\in\M_{j}}d(\vte, \bar\vte)\le 2\zeta'_n\}$. 

Assume first that the sequence of times $\{\widetilde T_n: \, n\ge
1\}$ is bounded above by some $T>0$.  For each integer $n>0$, let
$\widehat\pi^n$ be the path in $C([0,T-\widetilde T_n],\M_{+,1})$ given by
$\widehat\pi^n_t(d\theta) = \widehat\rho^n(t,\theta) d\theta$, where $\widehat\rho^n$ is the
solution of the hydrodynamic equation \eqref{rdeq} with initial
condition $\rho^n(\widetilde T_n)$, where $\pi^n(\widetilde T_n, d\theta) =
\rho^n(\widetilde T_n, \theta) d\theta$. Since $d(\pi^n(\widetilde T_n), \mc
M_j) \le 2 \zeta'_n$, by Lemma \ref{rep}, $\widehat\pi^n (T-\widetilde T_n)$
converges to some element of $\mc M_j$.

Let $\widetilde\pi^n$ be the path in $C([0,T],\M_{+,1})$ given by
\begin{equation*}
\widetilde\pi^n_t \;=\;
\begin{cases}
\pi^n_t & \hbox{ if }\; 0\le t\le \widetilde T_n \;,\\
\widehat \pi^n(t-\widetilde T_n) & \hbox{ if }\; \widetilde T_n\le t\le T \;.
\end{cases}
\end{equation*}
By definition of $\widetilde\pi^n$, $I_{T}(\widehat \pi^n) = I_{T}(\widehat
\pi^n|\gamma_n) = I_{T'_n}(\pi^n|\gamma_n) < v_{ij} \,-\, \e/2$.
Since $I_T$ has compact level sets and since $\pi^n_0(d\theta) =
\gamma_{n}(\theta)d\theta$ belongs to $\Gamma_{i}(\zeta'_n)\cap \M_{+,1}$, there
exists a subsequence of $\widetilde\pi^n$ converging to some $\pi$ in
$C([0,T],\M_{+,1})$ such that $\pi_0 \in \mc M_i$, $\pi_T \in \mc
M_j$, and $I_T(\pi)\le v_{ij}-\e/2$, which contradicts the definition
of $v_{ij}$.

If the sequence $\{\widetilde T_n: \, n\ge 1\}$ is not bounded, we may
repeat the reasoning presented in the proof of Theorem \ref{lsc} to
replace the path $\pi^{n}$ by a path $\bar\pi^{n}$ which satisfies an
inequality analogous to \eqref{contra} (with extra factors of
$\varepsilon$) and whose entry time to the set $\{\vte\in \M_{+}:\,
\zeta'_n \le \inf_{\bar\vte\in\M_{j}} d(\vte, \bar\vte)\le 2 \zeta'_n \}$ is
uniformly bounded in $n$. This completes the proof of the lemma, since
the bounded case has been treated above.
\end{proof}

We now prove the upper bound.

\begin{lemma}
\label{bub}
For every $1\le i \not = j\le l$, $\e>0$, there exists $\delta_{17}=\de_{17}(\e)$
such that for all $0<\beta_0 < \beta_1< \delta_{17}$,
\begin{equation*}
\limsup_{N\to\infty}\frac{1}{N}\log\sup_{\eta\in\partial B_i^N}\bb P_{\eta}
(\eta_{\tau}\in\partial B_j^N) \;\le\; -\, v_{ij} + \e \;.
\end{equation*}
\end{lemma}

\begin{proof}
Fix $1\le i \not = j\le l$, $\varepsilon>0$. Let $\zeta_{1}>0$, $T>0$ be
chosen according to Lemma \ref{lem6}, and fix
$0<\beta_0<\beta_1<\zeta_1$. By the strong Markov property,
\begin{equation*}
\sup_{\eta\in\partial B_i^N}\bb P_{\eta}(\eta_{\tau}\in\partial B_j^N)
\;\le\; \sup_{\eta\in\Gamma_i^N}
\bb P_{\eta}(\eta (H_{\partial B^N}) \in\partial B_j^N) \;,
\end{equation*}
where $H_D$, $D\subset X_N$, represents the hitting time of the set
$D$ and $\eta(t) = \eta_t$.  For each integer $N>0$, fix a
configuration $\eta^N$ in $\Gamma_i^N$ such that
\begin{equation*}
\bb P_{\eta^N}(\eta (H_{\partial B^N}) \in\partial B_j^N)
\;=\; \sup_{\eta\in\Gamma_i^N}
\bb P_{\eta}(\eta (H_{\partial B^N}) \in\partial B_j^N) \;.
\end{equation*}

By Lemma \ref{boundedness}, $v_{ij}<\infty$. Thus, by Lemma
\ref{lem4}, there exists $T_{\beta_0}>0$ such that
\begin{equation*}
\limsup_{N\to\infty}
\frac{1}{N}\log\sup_{\eta\in X_N} \bb P_{\eta}\left[
H_{B}^N \ge T_{\beta_0} \right]  \;\le\; -\, v_{ij} \;.
\end{equation*}
We may assume that $T_{\beta_0}>T$, where $T$ is the time introduced at
the beginning of the proof.  On the other hand, since $\eta^N \in
\Gamma_i^N$,
\begin{equation*}
\limsup_{N\to\infty}\frac{1}{N}\log \bb P_{\eta^N}
\Big(\eta (H_{\partial B^N}) \in\partial B_j^N \,,\, 
H^N_B \le T_{\beta_0} \Big) \;\le\; \limsup_{N\to\infty}\frac{1}{N}\log
\bb P_{\eta^N}(H^N_{\Ga_j}\le T_{\beta_0})\;.
\end{equation*}
By intersecting the set $\{\eta (H_{\partial B^N}) \in\partial B_j^N
\}$ with the set $\{H^N_B \le T_{\beta_0}\}$ and its complement, since
\begin{equation}
\label{15}
\limsup_{N\to\infty}\frac{1}{N}\log\{a_N+b_N\}
\;\le\; \max\left\{\limsup_{N\to\infty}\frac{1}{N}\log a_N,
\limsup_{N\to\infty}\frac{1}{N}\log b_N\right\}\;,
\end{equation}
it follows from the two previous estimates that
\begin{equation}
\label{09}
\begin{aligned}
& \limsup_{N\to\infty} \frac{1}{N} \log \sup_{\eta\in\Gamma_i^N}
\bb P_{\eta} \big(\eta (H_{\partial B^N}) \in\partial B_j^N \big) \\
&\qquad \; \le\;  \max \Big\{\limsup_{N\to\infty}\frac{1}{N}\log
\bb P_{\eta^N}(H^N_{\Ga_j}\le T_{\beta_0}) \,,\, -\, v_{ij} \Big\} \;.
\end{aligned}
\end{equation}

Let $\mc C_j$ be the set introduced in Lemma \ref{lem6} associated to
the pair $(\beta_1, T_{\beta_0})$. Since $\mc C_j$ is a closed
set, and since $\{H^N_{\Ga_j}\le T_{\beta_0}\} \subset \mc C_j$, by the
dynamical large deviations upper bound and by the compactness of
$\M_{+}$, there exists $\gamma(\theta)d\theta \in \Gamma_i$ such that
\begin{equation*}
\limsup_{N\to\infty}\frac{1}{N}\log
\bb P_{\eta^N}(H^N_{\Ga_j}\le T_{\beta_0}) \;\le\;
\limsup_{N\to\infty}\frac{1}{N}\log Q_{T_{\beta_0}, \eta^N}(\mc C_j)
\;\le\; -\inf_{\pi\in \mc C_j} I_{T_{\beta_0}}(\pi|\gamma) \;.
\end{equation*}
By Lemma \ref{lem6}, the last term is bounded above by $-v_{ij} +
\varepsilon$. This completes the proof of the lemma in view of
\eqref{09}. 
\end{proof}

The proof of the next result is similar to the ones of Lemmata 3.1 and
3.2 of chapter 6 in \cite{fw1}. Recall the notation introduced above
equation \eqref{19}. Consider a set $\Omega$, which is not assumed to
be countable.  Denote by $\Omega_i$, $1\le i\le l$, a partition of
$\Omega$: $\Omega = \cup_{1\le i\le l} \Omega_i$, $\Omega_i \cap
\Omega_j = \varnothing$ for $i\not = j$. Let $(Z_n : n\ge 0)$ be a
discrete-time Markov chain on $\Omega$ and denote by $p(x,dy)$, $x\in
\Omega$, the transition probability of the chain $Z_n$. Assume that
any set $\Omega_j$ can be reached from any point $x\in\Omega$:
$\sum_{n\ge 0} P_x[Z_n \in \Omega_j]>0$.

\begin{lemma}
\label{invmea}
Suppose that there exist nonnegative numbers $p_{ij}, \tilde p_{ij}$,
$1\le i \not =j\le l$, and a number $a>1$ such that
\begin{equation*}
\frac 1a\, p_{ij} \;\le\; P(x,\Omega_j) \;\le\; a \, \tilde p_{ij}
\quad \text{for all } \; x\in \Omega_i \;, \, i\,\ne\, j \;.
\end{equation*}
Then,
\begin{equation*}
\frac 1{a^{2(l-1)}} \, \frac {Q_i}{\sum_{1\le j\le l}\tilde Q_j}
\; \le \; \nu(\Omega_i) \;\le\;
a^{2(l-1)} \frac {\tilde Q_i}{\sum_{1\le j\le l} Q_j}  
\end{equation*}
for any invariant probability measure $\nu$, where $Q_i$, $\tilde
Q_{i}$ are given by
\begin{align*}
Q_i \;=\; \sum_{g\in  \ms T(i)}\prod_{(m, n)\in g} p_{mn}\quad
\hbox{ and } \quad
\tilde Q_i \;=\; 
\sum_{g\in  \ms T (i)}\prod_{(m, n)\in g} \tilde p_{mn} \;.
\end{align*}
\end{lemma}

Let
\begin{equation*}
\tilde w_{i} \;=\; \min_{g\in \ms T(i)} \sum_{(m, n)\in g} \tilde v_{mn}\;.
\end{equation*}
By the argument presented in the proof of \cite[Lemma 4.1]{fw1}, we
have $w_{i}=\tilde w_{i}$ for all $1\le i\le l$.  We are now in a
position to state the main result of this subsection.

\begin{proposition}
\label{main5.1}
For every $\e>0$, there exists $\delta_{18} = \delta_{18}(\e)$ such
that for all $1\le i \le l$, $0< \beta_0 < \beta_1 <
\delta_{18}$, 
\begin{align*}
& \limsup_{N\to\infty}
\frac{1}{N}\log \nu^N(\partial B_i^N) \;\le\; -\, \overline{w}_i \,+\, \e \;, \\
&\quad \liminf_{N\to\infty}
\frac{1}{N}\log \nu^N(\partial B_i^N) \;\ge\; -\, \overline{w}_i \,-\, \e \;.
\end{align*}
\end{proposition}

\begin{proof}
Since $w_i = \tilde w_i$, the assertion of this proposition is a
straightforward consequence of Lemmata \ref{blb}, \ref{bub} and
\ref{invmea}.
\end{proof}

\subsection{Lower bound}

We prove in this subsection the large deviations lower bound, that is,
for any open subset $\mc O$ of $\M_{+}$,
\begin{equation*}
\liminf_{N\to\infty}\frac{1}{N}\log \mathcal P^{N}(\mc O)
\;\ge\; -\inf_{\vte\in\mc O} W(\vte) \;.
\end{equation*}

Recall the definition of the stopping time $\tau$ introduced in
\eqref{10}.  Following \cite{fw1,bg,f}, we represent the stationary
measure $\mu^{N}$ of a subset $A$ of $X_{N}$ as
\begin{equation}
\label{rep2}
\mu^{N}(A) \;=\; \frac{1}{C_N} \int_{\partial B^N}
\bb E_{\eta}\Big(\int_0^{\tau}
{\bf 1} \{\eta_s\in A\}\, ds \Big)\, d\nu^N(\eta) \;,
\end{equation}
where
\begin{equation*}
C_{N}\;=\;\int_{\partial B^{N}} \mathbb{E}_{\eta}(\tau)\, d\nu^{N}(\eta) \;.
\end{equation*}

The first lemma provides an estimate on the normalizing constant $C_{N}$.

\begin{lemma}
\label{nor}
For any $\e>0$, there exists $\delta_{19}=\de_{19}(\e)$ such that for all
$0<\beta_0<\beta_1<\delta_{19}$,
\begin{equation*}
\limsup_{N\to\infty}\frac{1}{N}\log C_{N} \;\le\; \e \;.
\end{equation*}
\end{lemma}

\begin{proof}
Fix $\varepsilon>0$ and let $\zeta_1$ be a positive number such that
$2\zeta_1$ is smaller than the constants $\delta_{13}(\varepsilon)$
introduced in Lemma \ref{conn} and smaller than the constant
$\delta_{14}(\varepsilon)$ introduced in Lemma \ref{lem8}. Fix
$0<\beta_0<\beta_1<\zeta_1$. Since
$H_{\Gamma_{i}}^{N}<\tau$ when the process starts from $\partial
B_{i}^{N}$, by the Strong Markov property,
\begin{align*}
C_{N}
&\;=\;\sum_{i=1}^{l}\int_{\partial B_{i}^{N}}
\mathbb{E}_{\eta}(\tau)\, d\nu^{N}(\eta)  
\;=\;\sum_{i=1}^{l}\int_{\partial B_{i}^{N}}
\mathbb{E}_{\eta} \big(\tau \, {\bf 1} \{
H_{\Gamma_{i}}^{N}<\tau\} \big) \, d\nu^{N}(\eta)  \\
&\qquad \;\le\; \sum_{i=1}^{l}\int_{\partial B_{i}^{N}}
\mathbb{E}_{\eta}(H_{\Gamma_{i}}^{N})\, d\nu^{N}(\eta)
\;+\; \sup_{\eta\in X_{N}} \mathbb{E}_{\eta} (H_{B}^{N}) \;.
\end{align*}
By Corollary \ref{lem7} and by \eqref{15}, it remains to show that for
every $1\le i\le l$,
\begin{equation}
\label{nor1}
\limsup_{N\to\infty} \frac{1}{N}\log 
\sup_{\eta\in \partial B_{i}^{N}} \mathbb{E}_{\eta} 
(H_{\Gamma_{i}}^{N}) \;\le\; 3\e \;.
\end{equation}

Fix $1\le i\le l$.  We claim that there exists $N_0$ such that for all
$N\ge N_0$,
\begin{equation}
\label{12}
\sup_{\eta\in \Ga^N_{i} \cup \partial B_{i}^{N}} 
\mathbb P_{\eta}(H_{\Ga_{i}}^{N} \ge 3) \;\le\; 1-\exp{\{-3N\e\}} \;.
\end{equation}
To prove this assertion, for each integer $N>0$, consider a
configuration $\eta^{N}$ such that $\pi^N(\eta^N)$ belongs to $\mc
B_{2\beta_1}[\M_{i}] := \overline{\mc B_{2\beta_1}(\M_{i})}$ (which
contains the set $\Ga^N_{i} \cup \partial B_{i}^{N}$) and such that
\begin{equation*}
\mathbb P_{\eta^{N}}(H_{\Ga_{i}}^{N}<3) \;=\;
\inf_{\eta\in\mc B_{2\beta_1}^{N}[\M_{i}]}
\mathbb P_{\eta}(H_{\Ga_{i}^{N}}<3) \;.
\end{equation*}
Recall that each subsequence of $\pi^N(\eta^N)$ contains a
sub-subsequence converging in $\M_{+}$ to some $\vte$ which belongs to
$\M_{+,1}$.  Therefore, we may assume that $\pi^N(\eta^N)$ converges
to $\vte (d\theta)=\gamma (\theta)d\theta$ and that $\vte$ belongs to the closure of
$\mc B_{2\beta_1}(\bar\vte_{i})$: $\vte \in \mc
B_{2\beta_1}[\bar\vte_{i}] := \overline{\mc
  B_{2\beta_1}(\bar\vte_{i})}$ for some $\bar\vte_{i} \in \mc M_i$.

Let $\overline{\pi}_0 = \vte\in \mc B_{2\beta_1}[\bar\vte_{i}] \subset
\mc B_{2\zeta_{1}}(\bar\vte_{i})$, and let $\overline{\pi}$ be a path
in $D([0,1], \M_{+,1})$ provided by Lemma \ref{conn}. Let $\widetilde
\pi$ be a path in $D([0,1], \M_{+,1})$ provided by Remark \ref{rem2}:
$\widetilde\pi_{0}=\bar\vte_{i}$, $\widetilde\pi_{1}\in\M_{+}\setminus
\mc B_{2\beta_1}[\M_{i}]$ and $I_{1}(\widetilde \pi)\le\e$.  Define
the path $\pi$ in $D([0,2], \M_{+,1})$ by concatenating the paths
$\overline{\pi}$ and $\widetilde \pi$:
\begin{equation*}
\pi_{t} \;=\; \overline{\pi}_{t} \quad\text{for}\quad
0\le t \le 1 \;, \qquad \pi_{t} \;=\; \widetilde\pi_{t-1}
\quad\text{for}\quad 1\le t \le 2   \;.
\end{equation*}
The path $\pi_t$, $0\le t\le 2$, starts from $\vte$, hits $\mc M_i$
and then $\mc B_{2\beta_1}[\M_{i}]^c$. Its cost $I_{2}(\pi)$ is
bounded by $2\varepsilon$.

Denote by $\Lambda_{\beta_1/2}(\pi)$ the $\beta_1/2$-open
neighborhood of the trajectory $\pi$ in $D([0,2], \M_{+,1})$. Since
$\Lambda_{\beta_1/2}(\pi) \subset \{H_{\Ga_{i}}^{N}<3\}$, by the
dynamical large deviations lower bound, by definition of the sequence
$\eta^N$ and since $I_{2}(\pi) \le 2\varepsilon$, for $N$ large enough,
\begin{equation*}
\mathbb P_{\eta}(H_{\Ga_{i}}^{N}<3) \;\ge\;
\exp-N \Big\{\inf_{\pi'\in\Lambda_{\beta_1/2}(\pi)}I_{2}(\pi'|\ga) \,
+\, \e \Big\} \;\ge\; \exp{\{-3N\e\}} 
\end{equation*}
for all $\eta\in\mc B_{2\beta_1}^{N}[\M_{i}]$, which proves
\eqref{12}. 

The estimate \eqref{12} together with the arguments presented in Lemma
\ref{lem4} and Corollary \ref{lem7} gives the bound \eqref{nor1},
which completes the proof.
\end{proof}

\noindent{\bf Proof of the lower bound of Theorem \ref{sldp}}.  We
first claim that for any open set $\mc O$ of $\M_{+}$ containing
some $\bar\vte_i \in \mc M_i$, $1\le i \le l$,
\begin{equation}
\label{llb}
\liminf_{N\to\infty}\frac{1}{N}\log\mc P^N(\mc O) \;\ge\; - \,\overline{w}_i \;.
\end{equation}
Indeed, fix $1\le i \le l$, $\e>0$ and choose first $\beta_1>0$ and
then $0<\beta_0<\beta_1$ satisfying two conditions: (a) $\beta_0
<\zeta_1$, where $\zeta_1$ is the positive constant
$\delta_{11}(\beta_1/2, 1/2,1)$ provided by Lemma \ref{rep} for
$\bar\rho=\bar\rho_i$, and (b) the pair $(\beta_0, \beta_1)$ fulfills
the lower bound of Proposition \ref{main5.1} and Lemma
\ref{nor}. Assume, moreover, that $\mc
B_{2\beta_1}(\overline\vte_i)\subset \mc O$. Note that condition (a)
entails that $\beta_0<\beta_1/2$.

By \eqref{rep2}, and since $\tau\ge H^N_{\Gamma_i}$ if the initial
configuration belongs to $\partial B_i^N$,
\begin{align*}
\mc P^N(\mc O) \; &=\; \frac{1}{C_N}\int_{\partial B^N}
\bb E_{\eta}\left(\int_0^{\tau}
{\bf 1} \{\eta_s\in \mathcal{O}^{N} \}\, ds\right)\, d\nu^N(\eta) \\
\; &\ge\; \frac{1}{C_{N}}\int_{\partial B_i^N}\bb
E_{\eta}(H^N_{\Gamma_i}) \, d\nu^N(\eta) 
\;\ge\; \frac{1}{C_{N}}\, \nu^N(\partial B_i^N)\,
\inf_{\eta\in \partial B_i^N}\bb P_{\eta}(H^N_{\Gamma_i}\ge 1) \;.
\end{align*}
By Lemma \ref{nor} and Proposition \ref{main5.1}, to conclude the
proof of claim \eqref{llb}, it remains to show that
\begin{equation*}
\liminf_{N\to\infty}\frac{1}{N}\log\inf_{\eta\in \partial B_i^N}
\bb P_{\eta}(H^N_{\Gamma_i}\ge 1) \;\ge\; 0 \;.
\end{equation*}

For each integer $N>0$, let $\eta^N$ be a configuration in $\partial
B_i^N$ such that
\begin{equation*}
\bb P_{\eta^N}(H^N_{\Gamma_i}\ge 1)
\;=\; \inf_{\eta\in \partial B_i^N}\bb P_{\eta}(H^N_{\Gamma_i}\ge 1) \;.
\end{equation*}
Denote by $\eta^{N_k}$ a subsequence of $\eta^N$ which transforms the
$\liminf$ in a limit and let $\pi_0 (d\theta) = \gamma(\theta) d\theta$ be a limit
point of $\pi^{N_k}(\eta^{N_k})$. Observe that $\pi_0 \in \overline{\mc
B_{\beta_0}(\M_{i})}$ and denote by $\rho:[0,1]\times\T\to[0,1]$ the unique weak
solution of the Cauchy problem \eqref{rdeq} starting from
$\gamma$. By Lemma \ref{rep} and by the definition of $\beta_0$,
$\beta_1$, $\pi_t(d\theta) := \rho(t,\theta) d\theta$ belongs to $\mc
B_{\beta_1/2}(\mc M_i)$ for all $0\le t\le 1$.

Let $\mc N$ be the subset of $D([0,1], \M_{+})$ given by all
trajectories $\pi'_t$, $0\le t\le 1$, such that
$\sup_{0\le t\le 1} d(\pi'_t,\pi_t) < \beta_1/2$. Note that the set
$\mc N$ is open because $\pi_t$ is continuous. In particular, since
$\pi_t$ belongs to $\mc B_{\beta_1/2}(\mc M_i)$ for any $0\le t\le 1$,
$\{H^N_{\Gamma_i}\ge 1\} \supset (\pi^N)^{-1}(\mathcal N) :=
\{\eta_\cdot\in D([0,1],X_N):\pi^N(\eta_\cdot)\in \mc N\}$.
Therefore,
\begin{align*}
\liminf_{N\to\infty}\frac{1}{N}\log \bb P_{\eta^N}(H^N_{\Gamma_i}\ge 1)
\; &\ge\; \liminf_{N\to\infty}\frac{1}{N}\log{Q}_{1,\eta^N}(\mc N) \\
\; & \ge\; -\inf_{\pi'\in\mc N} I_1(\pi'|\gamma) 
\;\ge\; - I_1(\pi|\gamma) \;=\; 0 \;,
\end{align*}
which completes the proof of the claim.

It follows from \eqref{llb} that there exists a sequence $\e_N\to 0$
such that
\begin{equation}
\label{14}
\liminf_{N\to\infty}\frac{1}{N}\log
\mc P^N(\mc B_{\e_N}(\M_{i})) \;\ge\; - \, \overline{w}_i  \;.
\end{equation}

Fix an open subset $\mc O$ of $\M_{+}$.  In order to prove the lower
bound, it is enough to show that for any measure $\vte$ in $\mc O$,
$1\le i\le l$, $T>0$, and any trajectory $\pi$ in $D([0,T],\mc M_{+})$
with $\pi_0\in\mc M_i$, $\pi_T = \vte$,
\begin{equation}
\label{13}
\liminf_{N\to\infty}\frac{1}{N}\log\mc P^N(\mc O)
\;\ge\; -\, \overline{w}_i \;-\; I_T(\pi) \;.
\end{equation}

To prove this claim, fix an open subset $\mc O$ of $\M_{+}$, a measure
$\vte$ in $\mc O$, $1\le i\le l$, $T>0$, and a trajectory $\pi$ in
$D([0,T],\mc M_{+})$ with $\pi_0\in\mc M_i$, $\pi_T = \vte$. Since
$\mu^N$ is the stationary measure,
\begin{align*}
\liminf_{N\to\infty}\frac{1}{N}\log\mc P^N(\mc O)
& \;=\; \liminf_{N\to\infty}\frac{1}{N}\log\bb E_{\mu^{N}}
\big[ \, \bb P_{\eta}\left(\pi^N_T\in\mc O\right)\big] \\
& \;\ge\; \liminf_{N\to\infty}\frac{1}{N}\log \Big(
\mc P^N(\mc B_{\e_N}(\M_{i}))\inf_{\eta\in \mc B_N} \bb
P_{\eta}[\pi^N_T\in \mc O] \Big)\;,
\end{align*}
where $\mc B_N = \{\eta : \pi^N(\eta) \in \mc B_{\e_N}(\M_{i}) \}$.
Let $\eta^N$ be a configuration in $\mc B_N$ such that
\begin{equation*}
\bb P_{\eta^N}\left[\pi^N_T\in \mc O\right]
\;=\; \inf_{\eta\in\mc B_N} \bb P_{\eta}\left[\pi^N_T\in \mc O\right] \;.
\end{equation*}
Since $\eta^N$ belongs to $\mc B_N$ and $\varepsilon_N\to 0$, we may
assume, taking a subsequence if necessary, that $\pi^N(\eta^N)$
converges to some $\bar\vte_i (d\theta) = \bar\rho_i (\theta) d\theta \in\mc M_i$.
By \eqref{14}, the expression appearing in the penultimate displayed
formula is bounded below by
\begin{align*}
- \, \overline{w}_i \,+\, \liminf_{N\to\infty}\frac{1}{N}
\log\bb P_{\eta^N}\left[\pi^N_T\in \mc O\right]
\;=\; -\, \overline{w}_i + \liminf_{N\to\infty}\frac{1}{N}
\log{Q}_{T, \eta^N}(\mc O_T)\;,
\end{align*}
where $\mc O_{T}=\{\pi'\in D([0,T], \mc M_+) : \pi_{T} \in \mc O\}$.
Since the set $\mc O_T$ is open, by the lower bound of the dynamical
large deviations principle, the previous expression is bounded below
by
\begin{equation*}
-\, \overline{w}_i \;-\; \inf_{\pi'\in \mc O_T} I_T(\pi')
\;\ge\; -\, \overline{w}_i \;-\;  I_T(\pi) \; .
\end{equation*}
In view of \eqref{27}, this completes the proof of \eqref{13} and the
one of the lower bound.

\subsection{Upper bound}
We prove in this subsection the large deviations upper bound. The
proof relies on the next two lemmata. The proof of the first one is
similar to the proof of Lemma \ref{lem6} and is left to the reader.

For a closed subset $\mc C$ of $\mc M_+$ and $T>0$, let $C_T$ be the
subset of $D([0, T],\M_{+})$ consisting of all paths $\pi$ for which
there exists $t$ in $[0,T]$ such that $\pi(t)$ or $\pi(t-)$ belongs to
$\mc C$. Note that $C_T$ is a closed subset of $D([0,T],\M_{+})$.

\begin{lemma}
\label{lem9}
Fix a closed subset $\mc C$ of $\mc M_+$ such that $\inf_{\vte \in \mc
  C} V_{i}(\vte)<\infty$. For every $\varepsilon >0$, there exist
$\delta_{20} = \delta_{20}(\mc C, \varepsilon)>0$ and $T_{20}>0$ such
that for all $1\le i\le l$, $0<\beta_1<\delta_{20}$, $T'\ge T_{20}$,
$\gamma(\theta)d\theta\in \Gamma_i$,
\begin{equation*}
\inf_{\pi\in C_{T'} } I_{T'}(\pi|\gamma)
\;\ge\; \inf_{\vte\in \mc C} V_{i}(\vte) -\e \;.
\end{equation*}
\end{lemma}

Recall the definition of the set $B$ introduced just before Lemma
\ref{lem4} and recall from Corollary \ref{lem7} that the set $B$ is
attained immediately. In particular, if the reaction-diffusion model
has to reach a set $\mc C$ before it hits $B$, it has to follow
straightforwardly the optimal trajectory to $\mc C$. The cost of such
trajectory has been estimated in the previous lemma, providing the
next result.

\begin{lemma}
\label{lem3}
Fix $1\le i\le l$ and a closed subset $\mc C$ of $\mc M_+$.  For every
$\e>0$, there exist $\delta_{21} = \delta_{21}(\mc C, \varepsilon)>0$
such that for all $0<\beta_0 < \beta_1<\delta_{21}$,
\begin{equation*}
\limsup_{N\to\infty}\frac{1}{N}
\log\sup_{\eta\in \Gamma^N_i}
\bb P_{\eta}\left[H_{\mc C}^N < H_{B}^N\right]
\;\le\; - \inf_{\vte \in\mc C} V_{i}(\vte) \,+\, \e \;.
\end{equation*}
\end{lemma}

\begin{proof}
Fix $\e>0$, $1\le i\le l$, and a closed subset $\mc C$ of $\mc
M_+$. We may assume that the left hand side of the inequality appearing
in the statement of the lemma is finite. This implies that
$\pi^{-1}_N(\mc C) \cap X_N \not = \varnothing$ for infinitely many
$N$'s. Let $\{\eta^{N_k} : k\ge 1\}$ be a sequence of configurations
such that $\pi^{N_k}(\eta^{N_k}) \in \mc C$. Since $\mc M_+$ is
compact, taking a subsequence, if necessary, we may assume that
$\pi^{N_k}(\eta^{N_k})$ converges to a measure, denoted by $\vte$, which
belongs to $\mc M_{+,1}$. Since $\mc C$ is closed, $\vte\in \mc C$ so
that $\mc C \cap \mc M_{+,1} \not = \varnothing$. In particular, By
Lemma \ref{boundedness}, $\inf_{\vte \in \mc C} V_{i}(\vte)<\infty$.

Let $\zeta_1$, $R_1$ be the constants $\delta_{20}$, $T_{20}$ given by
Lemma \ref{lem9}. Fix $0<\beta_0<\beta_1<\zeta_1$. Since $\inf_{\vte
  \in \mc C} V_{i}(\vte)<\infty$, by Lemma \ref{lem4}, there exists
$R_2>0$ such that
\begin{equation}
\label{16}
\limsup_{N\to\infty}\frac{1}{N}
\log\sup_{\eta\in X_N} 
\bb P_{\eta}\left[ H_{B}^N \ge R_2 \right] 
\;\le\; - \, \inf_{\vte \in \mc C} V_{i}(\vte) \;.
\end{equation}
Let $T=\max\{R_1, R_2\}$.

Recall the definition of the set $C_T$ introduced just above the
statement of Lemma \ref{lem9} and the fact that $C_T$ is a closed
subset of $D([0,T],\M_{+})$. Note also that $\{H^N_{\mc C}\le T\}
\subset C_T$.

Let $\eta^N$ be a configuration in $\Gamma_i^N$ such that
\begin{equation*}
\bb P_{\eta^N}\left[H^N_{\mc C}\le T\right]
\;=\; \sup_{\eta\in \Gamma_i^N}
\bb P_{\eta}\left[H^N_{\mc C}\le T\right] \;.
\end{equation*}
Taking a subsequence if necessary, we may assume that $\pi^N(\eta^N)$
converges to some measure $\vte(d\theta)=\gamma(\theta)d\theta$ in
$\Gamma_i\cap\M_{+,1}$.  Since $\{H^N_{\mc C}\le T\} \subset C_T$ and
since $C_T$ is a closed set, by the dynamical large deviations upper
bound and by Lemma \ref{lem9},
\begin{align*}
& \limsup_{N\to\infty}\frac{1}{N}\log \sup_{\eta\in \Gamma^N_i} 
\bb P_{\eta^N}\big(H^N_{\mc C} < H^N_B \le T \big)
\;\le\; \limsup_{N\to\infty}\frac{1}{N}\log
\bb P_{\eta^N}\big(H^N_{\mc C}\le T \big) \\
&\qquad \;\le\; \limsup_{N\to\infty}\frac{1}{N}
\log{Q}_{T,\eta^N}(C_T)
\;\le\; -\inf_{\pi\in C_T} I_{T}(\pi|\gamma) 
\;\le\; - \inf_{\vte \in\mc C} V_{i}(\vte) \;+\; \varepsilon \;.
\end{align*}
By \eqref{16}, by this estimate and by \eqref{15}, 
\begin{equation*}
\limsup_{N\to\infty}\frac{1}{N}
\log\sup_{\eta\in \Gamma^N_i}
\bb P_{\eta}\left[H_{\mc C}^N < H_{B}^N\right]
\;\le\;- \inf_{\vte \in\mc C} V_{i}(\vte) \;+\; \varepsilon \;,
\end{equation*}
which completes the proof of the lemma.
\end{proof}

\noindent{\bf Proof of the upper bound of Theorem \ref{sldp}.}
Let $\mc C$ be a closed subset of $\mc M_+$.  Assume first that
$\M_{i} \cap \mc C = \varnothing$ for any $1\le i \le l$. In this case,
let $\beta_1>0$ be such that $\cup_{1\le i\le l} \mc
B_{2\beta_1}(\M_{i})\cap \mc C = \varnothing$.

By the representation \eqref{rep2} of the stationary measure $\mu^N$, 
\begin{align*}
\mc P^{N}(\mc C) \; & =\; \mu^N(\mc C^N)
\;=\; \frac{1}{C_N}\int_{\partial B^N}
\bb E_{\eta}\left(\int_0^{\tau}
{\bf 1} \{\eta_s\in \mc C^N \}\, ds\right)d\nu^N(\eta) \\
& \;\le\; \frac{1}{C_N}\sum_{i=1}^l\nu^N(\partial B_i^N)
\sup_{\eta\in\partial B^N_i}\bb E_{\eta}\left(\int_0^{\tau}
{\bf 1} \{\eta_s\in \mc C^N \}\, ds\right) \;. 
\end{align*}

A configuration in $X_N$ can jump to at most $2N$ different
configurations and the jump rates are bounded by $N^2$.  Since any
trajectory in $D(\bb R_+,X_N)$ has to jump at least once before the
stopping time $\tau$, the constant $C_N$ appearing in the denominator
is bounded below by $c_0/N^{3}$ for some positive constant
$c_0$. Hence, by \eqref{15} and by Proposition \ref{main5.1}, in order
to prove the upper bound it is enough to show that for each $1\le i
\le l$,
\begin{equation}
\label{ubse}
\limsup_{N\to\infty}\frac{1}{N}\log\sup_{\eta\in\partial B^N_i}
\bb E_{\eta}\left(\int_0^{\tau}
{\bf 1} \{\eta_s\in \mc C^N\}\, ds\right)
\;\le\; -\, \inf_{\vte\in\mc C} V_{i}(\vte) \;+\; \e \;.
\end{equation}

The time integral appearing in the previous formula vanishes if $\tau
\le H^N_{\mc C}$. We may therefore introduce the indicator of the set
$H^N_{\mc C} \le \tau$. After doing this and applying the strong
Markov property, we obtain that the left hand side of the previous
inequality is less than or equal to
\begin{align*}
\limsup_{N\to\infty}\frac{1}{N}\log\sup_{\eta\in\partial B^N_i}
 \bb P_{\eta}\left[H_{\mc C}^N < \tau\right]
\sup_{\eta\in \mc C^N} \bb E_{\eta}\left(\tau\right) \;.
\end{align*}

Since the distance between the empirical measure before and after a
jump is bounded by $C/N$, and since $\mc B_{2\beta_1}(\M_{i})\cap \mc
C = \varnothing$, for $N$ large enough, any trajectory in $D(\bb R_+,
X_N)$ starting at some configuration in $\partial B_i^N$, $\mc C^N$,
satisfies $H^N_{\Gamma_i} \le H^N_{\mc C}$, $\tau\le H^N_B$,
respectively.  Hence, by the strong Markov property, the previous
expression is bounded above by
\begin{equation*}
\limsup_{N\to\infty}\frac{1}{N}
\log\sup_{\eta\in\Gamma^N_i}\bb P_{\eta}\left[H_{\mc C}^N < H_{B}^N\right]
\sup_{\eta\in \mc C^N}  \bb E_{\eta}\left(H_B^N\right) \;.
\end{equation*}
By Corollary \ref{lem7} and Lemma \ref{lem3}, the previous expression
is bounded by $-\, \inf_{\vte\in\mc C} V_{i}(\vte) \;+\; \e$, which
completes the proof of \eqref{ubse} and the one of the upper bound
in the case $\M_{i} \cap \mc C = \varnothing$ for $1\le i \le l$.

We turn to the general case. We first claim that for each $1\le i\le
l$, $\varepsilon>0$, there exists $\zeta_1>0$ such that for all
$0<\beta_0<\zeta_1$,
\begin{equation}
\label{17}
\limsup_{N\to\infty}\frac{1}{N}\log\mc P^{N}(\mc B_{\beta_0}(\mc M_{i}))
\;\le\; - \, \overline{w}_{i} \;+\; 2 \varepsilon \;.
\end{equation}

Indeed, fix $\varepsilon>0$ and set $\zeta_1 = \min\{\delta_{18},
\delta_{19}\}$, where $\delta_{18}>0$ is the constant provided by
Proposition \ref{main5.1} and $\delta_{19}>0$ is the one given by
Lemma \ref{nor}. Fix $\beta_0<\zeta_1$. By the representation
\eqref{rep2} of the stationary measure $\mu^N$,
\begin{align*}
\mc P^{N}(\mc B_{\beta_0}(\mc M_{i})) \; & =\; 
\mu^N(\mc B^N_{\beta_0}(\mc M_{i}))
\;=\; \frac{1}{C_N}\int_{\partial B^N}
\bb E_{\eta}\left(\int_0^{\tau}
{\bf 1} \{\eta_s\in  \mc B^N_{\beta_0}(\mc M_{i}) \}\, ds\right)d\nu^N(\eta) \\
& \;\le\; \frac{1}{C_N}\sum_{j=1}^l\nu^N(\partial B_j^N)
\sup_{\eta\in\partial B^N_j}\bb E_{\eta}\left(\int_0^{\tau}
{\bf 1} \{\eta_s\in \mc B^N_{\beta_0}(\mc M_{i}) \}\, ds\right) \;. 
\end{align*}

We have seen in the first part of the proof that $\limsup_N N^{-1}
\log C^{-1}_N \le 0$. On the other hand, for $\eta\in \partial B_j^N$,
$j\not = i$, $\tau \le H^N_{B_i}$, so that $\int_0^{\tau} {\bf 1}
\{\eta_s\in \mc B^N_{\beta_0}(\mc M_{i}) \}\, ds =0$. Finally, denote
by $\vartheta(t)$, $t>0$, the time translation of a trajectory by $t$.
For $\eta\in \partial B_i^N$, writing $\tau$ as $H^N_{\Gamma_i} +
H^N_B \circ \vartheta(H^N_{\Gamma_i})$, by the strong Markov property,
since $\beta_0<\zeta_1$, and by Proposition \ref{main5.1}, the left
hand side of \eqref{17} is bounded by
\begin{equation*}
-\, \overline{w}_i \;+\; \varepsilon \;+\; \limsup_{N\to\infty}\frac{1}{N}\log
\Big\{ \sup_{\eta\in\partial B^N_i}\bb E_{\eta} (H^N_{\Gamma_i}) \,+\,
\sup_{\eta\in\partial \Gamma^N_i} \bb E_{\eta} (H^N_{B}) \Big\}\;.
\end{equation*}
By \eqref{15}, \eqref{nor1} and Corollary \ref{lem7}, the limit superior of the
previous equation is bounded by $\varepsilon$, which completes the proof
of \eqref{17}.

Let $\mc C$ be a closed subset of $\mc M_+$ and fix $\varepsilon>0$. Let
$A$ be the set of indices $i$ such that $\mc C \cap \mc M_i \not =
\varnothing$. Let $\zeta_1$ be the positive constant introduced in
\eqref{17}, and choose $\beta_0 <\zeta_1$ such that $d(\mc C, \mc
M_j)>\beta_0$ for all $j\in A^c$. Since $\mc C \subset \cup_{i\in A}
\mc B_{\beta_0}(\mc M_i) \cup [ \mc C \setminus \{\cup_{i\in A} \mc
B_{\beta_0}(\mc M_i)\}]$, and since $\mc C \setminus \{\cup_{i\in A} \mc
B_{\beta_0}(\mc M_i)\}$ is a closed set which does not intersect the set
$\mc M_{\rm sol}$, by \eqref{15}, by \eqref{17} and by the first part of the
proof,
\begin{equation*}
\limsup_{N\to\infty}\frac{1}{N}\log\mc P^{N}(\mc C)
\;\le\; - \min \Big\{ \min_{i\in A} \overline{w}_i
\,,\, \inf_{\pi \in \mc C \setminus \{\cup_{i\in A} \mc
B_{\beta_0}(\mc M_i)\}} W(\pi) \Big\} \;+\; 2\varepsilon \;.
\end{equation*}
By \eqref{18}, $\overline{w}_i = W(\bar\vte_i)$ for $\bar\vte_i\in\mc
M_i$. On the other hand, since $\mc M_i \cap \mc C \not =
\varnothing$,
\begin{equation*}
\inf_{\pi\in \mc C} W(\pi) \;\le\; \min_{i\in A} W(\bar\vte_i)\; ,
\quad \inf_{\pi\in \mc C} W(\pi) \;\le\; 
\inf_{\pi \in \mc \mc C \setminus \{\cup_{i\in A} \mc
B_{\beta_0}(\mc M_i)\}} W(\pi)\;,
\end{equation*}
which completes the proof of the upper bound.

\section{Proof of Theorem \ref{stablesol}}
\label{sec7}

We first show that if there exists a heteroclinic orbit from
$\phi\in\mc M_i$ to $\psi\in\mc M_j$, then the cost of going from $\mc
M_i$ to $\mc M_j$ vanishes.

\begin{lemma}
\label{l-f1}
Suppose that there exists a heteroclinic orbit from $\phi\in\mc M_i$
to $\psi\in\mc M_j$. Then, $v_{ij}=0$.
\end{lemma}

\begin{proof}
Fix $i\not = j$ in $\{1, \dots, l\}$ and assume that there exists a
heteroclinic orbit from $\phi\in\mc M_i$ to $\psi\in\mc M_j$, denoted
by $\rho (t,\theta)$, $t\in \bb R$.  By Proposition \ref{reg}, $\phi$ is
smooth, and, by Lemma \ref{lem1}, there exists $0<c<1/2$ such that
$c\le \phi(\theta)\le 1-c$. Since $\rho(t)$, converges in $C^1(\bb T)$ to
$\phi$, $\psi$ as $t\to - \infty$, $t\to + \infty$ respectively, by
Lemmata \ref{qpbound}, \ref{conn} and since the dynamical large
deviations rate functional vanishes along the solution of the
hydrodynamic equation, $v_{ij} = 0$.
\end{proof}

We now prove that $\rho(\theta) = r$ is a stable solution of the
reaction-diffusion equation \eqref{seeq} if $r$ is a local minimum of
$V$. 

\begin{lemma}
\label{l-f2}
Fix $1\le i\le l$. Let $\bar\vte_i(d\theta)=\bar\rho_i(\theta)d\theta, \bar\rho_i(\theta) = r$,
where $r$ is a local minimum of $V$.
Then, for all $\varepsilon>0$ there exist $c>0$ such that
\begin{equation*}
\inf \big\{ V_i(\varrho): \varrho \not\in \mc
B_\varepsilon(\bar\varrho_i) \big\} \;\ge\; c\;.
\end{equation*}
\end{lemma}

\begin{proof}
Suppose that $\inf\{ V_i(\varrho) : \varrho \in \mc B_{\delta}(\bar\vte_i)^c\} =
0$ for some $\delta>0$. In this case there exists a sequence of
density profiles $\gamma_n$ and of trajectories $\pi^n (t,d\theta) =
\rho^n(t,\theta) d\theta$, $0\le t\le T_n$, such that $\rho^n(0,\theta) = \bar\rho_i$,
$\rho^n(T_n,\theta) = \gamma_n(\theta)$, $\gamma_n(\theta) d\theta \in \mc
B_{\delta}(\bar\vte_i)^c$ and $I_{T_n} (\pi^n) \le 1/n$.

By Lemma \ref{lem11}, there exists $0<\varepsilon<\delta$ such that
$\pi_t \in \mc B_\delta (\bar\vte_i)$ for all $t\ge 0$ if $\pi_0 \in
\mc B_{2\varepsilon} (\bar\vte_i)$. Let $\tau_n$ be the time the trajectory
$\pi^n$ leaves the set $\mc B_{\varepsilon}(\bar\vte_i)$ for ever, and let
$\sigma_n$ be the hitting time of the set $\mc B_{\delta}(\bar\vte_i)^c$
after $\tau_n$:
\begin{equation*}
\tau_n \;=\; \sup \{t\le T_n : \pi^n_t \in \mc
B_{\varepsilon}(\bar\vte_i)\}\;, \quad
\sigma_n \;=\; \inf \{ t\ge \tau_n : \pi^n_t \in \mc
B_{\delta}(\bar\vte_i)^c\}\;.
\end{equation*}
Since in the interval $[\tau_n, \sigma_n]$ the trajectory $\pi^n$
remains in the set $\mc B_{\delta}(\bar\vte_i) \setminus \mc
B_{\varepsilon}(\bar\vte_i)$, if $\delta$ is small enough for $\mc
B_{\delta}(\bar\vte_i) \cap \mc B_{\delta}(\mc M_j) = \varnothing$ for
all sets $\mc M_j$, $j\not = i$, by Corollary \ref{costw}, $\sigma_n -
\tau_n$ is uniformly bounded by a finite constant, denoted by $T$.

Extend the definition of $\pi^n$ from the interval $[0,T_n]$ to $\bb
R_+$ by following the hydrodynamic trajectory after $T_n$:
$\pi^n(T_n+t,d\theta) = \tilde \rho^n(t,\theta) d\theta$, where $\tilde \rho_n$ is
the solution of the hydrodynamic equation with initial condition
$\rho^n_{T_n}$.  Let $\bar\pi^n_t$, $0\le t\le T$, be the trajectory
defined by $\bar\pi^n_t = \pi^n(\tau_n +t)$. Since $\pi^n$ belongs to $C([0,T_n],\M_{+,1})$,
note that $\bar\pi^n_0\in \partial \mc B_{\varepsilon}(\bar\vte_i)$, that $\bar\pi^n$
hits the set $\mc B_{\delta}(\bar\vte_i)^c$ in the time interval
$[0,T]$ and that $I_T(\bar\pi^n) \le 1/n$.

By the compactness of the level sets of $I_T$, the lower
semi-continuity of this functional and the compactness of the space
$\mc M_+$, there exists a subsequence $\bar\pi^{n_k}$ which converges to
some trajectory $\pi$ such that $\pi_0\in \partial \mc
B_{\varepsilon}(\bar\vte_i)$, $\pi$ hits the set $\mc
B_{\delta}(\bar\vte_i)^c$ in the time interval $[0,T]$ and $I_T(\pi)
=0$. By Lemma \ref{zero}, the density of $\pi_t$, denoted by $\rho_t$,
is a solution of the hydrodynamic equation. This contradicts the
property of $\varepsilon$ and concludes the proof of the lemma.
\end{proof}

\begin{proof}[Proof of Theorem \ref{stablesol}]
Recall the definition of the set of indices $I_s$, $I_u$.  We claim
that $\overline{w}_a >0$ for all $a\in I_u$. To prove this statement,
it is enough to show that for each $a\in I_u$, there exists $b \in
I_s$ such that $w_a > w_b$. 

Fix $a\in I_u$. By assumption, there exists $b \in I_s$ such that
$v_{ab} =0$. We claim that $w_a > w_b$. Indeed, on the one hand, by
Lemma \ref{l-f2}, $v_{bc}>0$ for all $c\not = b$. On the other hand,
let $g$ be a graph in $\ms T(a)$ such that $w_a = \kappa(g)$.  Recall
that we denote by $(d,e)$, $e\not = d\in \ms V$, the oriented edge
where $d$ is the child and $e$ the parent.  Let $c$ be the parent of
$b$ in $g$. Of course, $c$ might be $a$. Denote by $g'$ the tree in
$\ms T(b)$ obtained from $g$ by adding the oriented edge $(a,b)$ and
removing the the edge $(b,c)$, and note that $\kappa(g) + v_{ab} =
\kappa(g') + v_{bc}$. Since $w_b$ is the minimal value of $\kappa
(\tilde g)$, $\tilde g\in \ms T(b)$, $w_b \le \kappa (g')$ so that
$w_b + v_{bc} \le \kappa(g) + v_{ab} = w_a + v_{ab} = w_a$. The last
identity follows from the fact that $v_{ab}=0$ and the next to last
from the fact that $\kappa(g) = w_a$. Since $v_{bc}>0$, we conclude
that $w_b< w_a$, as claimed.

We claim that for every $\delta>0$,
\begin{equation}
\label{28}
\inf\big\{ W(\pi) : 
\pi \not \in \bigcup_{i\in I_s} \mc B_{\delta}(\bar\vte_i) \big\} \;>\; 0 \;. 
\end{equation}
Fix $\delta>0$. Since $\overline{w}_a>0$ for all $a\in I_u$, in
view of the definition of $W$, we only need to check that
\begin{equation*}
\inf \big\{  V_j (\pi) : 
\pi \not \in \bigcup_{i\in I_s} \mc B_{\delta}(\bar\vte_i) \big\}
\;>\; 0 
\end{equation*}
for each $j\in I_s$.  This is the content of Lemma \ref{l-f2},
proving \eqref{28}

To complete the proof of the theorem, it remains to observe that the
complement of $\cup_{i\in I_s} \mc B_{\delta}(\bar\vte_i)$ is a closed
set and to apply the upper bound of the static large deviations
principle stated in Theorem \ref{sldp}. The theorem is proved.
\end{proof}

\section{The Chafee-Infante equation}
\label{sec6}

We present in this section an example of a reaction-diffusion model
which fulfills the hypotheses of Theorems \ref{sldp} and
\ref{stablesol}. Actually, in this model a complete description of the
stationary solutions and of the heteroclinic orbits is available.

Fix $0<a<b$ and recall the definition of the potential $V=V_{a,b}$
introduced in \eqref{20}. Denote by $\bb T_{1/2\pi}$ the one-dimensional
torus with length $(2\pi)^{-1}$.  Let $p = (1/2) \sqrt{a/b}$, $c = 2(2\pi)^2$
and define $\phi :\bb R_+ \times \bb T_{1/2\pi} \to \bb R$ as
\begin{equation*}
\phi(t,\theta) \;=\; \frac 1p \, \Big\{ \rho \big( c\, t \,,\, 2\pi  \theta \big)
\,-\, \frac 12 \Big\}\;. 
\end{equation*}
A simple computation shows that $\rho$ solves the equation \eqref{20}
if and only if $\phi$ solves
\begin{equation}
\label{21}
\partial_t \phi \;=\; \Delta\phi \;+\; \lambda
(1-\phi^2)\phi\;, 
\end{equation}
where $\lambda = 4ca= 8 (2\pi)^2 a$. Note that $\phi$ takes values in
the interval $[- \sqrt{b/a}, \sqrt{b/a}]$. This is the so-called Chafee-Infante
equation \cite{ci} with periodic boundary condition.

A complete characterization of the stationary solutions of the
Chafee-Infante equation with periodic boundary conditions is presented
in \cite[Proposition 1.1]{kmy}. In our context it can be stated as
follows. Let $\rho_\pm$ be the minima of $V$: $\rho_\pm = (1/2) \pm p
= (1/2) [1 \pm \sqrt{a/b}]$.

\begin{theorem}
\label{tkmy}
For all $0<a<b$, the equation \eqref{21} admits three constant
stationary solutions: $\psi_{\pm} = \rho_\pm$, $\psi_{1/2}=1/2$. For
all nonnegative integers $m$ such that
$1\le m^2 < \lambda = 32 \pi^2 a$, up to translations, there exists a
non-constant periodic stationary solution $\phi_m=\phi_{m,\lambda}$
with $m$ periods in $\bb T_{1/2\pi}$. Moreover,
$\lim_{\lambda \downarrow m^2} \phi_{m,\lambda} = 1/2$ in
$C^2(\bb T_{1/2\pi})$.  The reaction-diffusion equation \eqref{rdeq}
with potential $V_{a,b}$ has no other stationary solutions.
\end{theorem}

The heteroclinic orbits of the Chafee-Infante equation with periodic
boundary conditions have been characterized in \cite{frw}.  Next
result follows from Theorems 1.3 and 1.4 of \cite{frw},

\begin{theorem}
\label{tfrw}
There are heteroclinic orbits from $\psi_{1/2}$ to $\psi_{\pm}$, and
from $\psi_{1/2}$ to $\phi_m$ for all integers $1\le m^2 <\lambda$.
Fix $1\le n^2 <\lambda$. There are heteroclinic orbits from $\phi_n$
to $\psi_{\pm}$, and from $\phi_n$ to $\phi_m$ for all integers $1\le
m^2<n^2$. There are no other heteroclinic orbits.
\end{theorem}

Next proposition follows from the previous results and from Theorem
\ref{stablesol}.

\begin{proposition}
\label{p01}
Consider a reaction-diffusion model which satisfies the assumptions of
Theorem \ref{sldp} and which gives rise to the hydrodynamic equation
\eqref{20} with $0<a<b$. Let $\bar\vte_\pm (d\theta) = \rho_\pm d\theta$. Then,
for every $\delta>0$, there exist $c>0$ and $N_0\ge 1$ such that for
all $N\ge N_0$,
\begin{equation*}
\mc P^N \big(\mc B_\delta(\bar\vte_{-}) \cup \mc B_\delta(\bar\vte_{+}) \big)
\;\ge \; 1 \;-\; e^{-cN}\;.
\end{equation*}
\end{proposition}

If the jump rates are invariant under a global flipping of the
configuration: $c(\eta) = c(\bs 1-\eta)$, where $\bs 1$ is the
configuration with all sites occupied, there is a symmetry between
occupied and vacant sites so that $\mc P^N (\mc B_\delta(\bar\vte_+))
= \mc P^N (\mc B_\delta(\bar\vte_-))$ for all $\delta>0$. Hence, with
this additional assumption, we may refine the previous proposition:

\begin{corollary}
\label{cor2}
Under the assumptions of Proposition \ref{p01}, if $c(\eta) = c(\bs
1-\eta)$, for every $0<\delta< (1/2) d(\bar\vte_{-}, \bar\vte_{+})$,
there exist $c>0$ and $N_0\ge 1$ such that for all $N\ge N_0$,
\begin{equation*}
\Big| \mc P^N (\mc B_\delta(\bar\vte_-)) - \frac 12\, \Big| \;\le\; e^{-cN} \;, \quad
\Big| \mc P^N (\mc B_\delta(\bar\vte_+)) - \frac 12\, \Big| \;\le\; e^{-cN}
\;. 
\end{equation*}
\end{corollary}

\smallskip\noindent{\bf An example.} We conclude this section with an
example of a reaction-diffusion model satisfying the assumptions of
Theorem \ref{sldp} and whose hydrodynamic equation is given by
\eqref{20}. 

Consider the reaction-diffusion model whose jump rate $c(\eta)$ is
given by
\begin{align*}
c(\eta) \; =\; a_2 \bs 1\{ \eta_{-1} \not = \eta_1\} 
\;+\; a_1 \bs 1\{ \eta_{-1} = \eta_1 = \eta_0 \} 
\; +\; a_0 \bs 1\{ \eta_{-1} = \eta_1 \not = \eta_0 \} \;.
\end{align*}
Let $\xi$ be the configuration obtained from $\eta$ by flipping all
occupation variables: $\xi(x) = 1 - \eta(x)$, $x\in\bb T_N$.  Since
$[1-\eta(0)] c(\eta) = \xi(0) c(\xi)$, $B(\rho) =
D(1-\rho)$. Moreover,
\begin{equation*}
F(\rho) \;=\; \frac 14 \Big\{ (a_0-3a_1-2a_2) (2\rho-1) - (a_0 + a_1-2a_2)
(2\rho-1)^3 \Big\}\;.
\end{equation*}
Fix $0<\mf a<\mf b$, and set $a_1=\mf a>0$. Choose $a_2\ge \mf a +
2 \mf b>0$ and set $a_0 = 2 a_2 + 4\mf b - \mf a \ge  \mf a +
8\mf b >0$. Since the three parameters are positive, the jump rate is
strictly positive as required.

As $a_0-3a_1-2a_2 = 4(\mf b - \mf a)$ and $a_0 + a_1-2a_2 = 4\mf b$,
in the variables $\mf a$, $\mf b$ the function $F$ becomes
\begin{equation*}
F(\rho) \;=\; (\mf b-\mf a) (2\rho-1) \;-\; \mf b (2\rho-1)^3 \;,
\end{equation*}
in conformity with \eqref{20} for $a= (\mf b-\mf a)/2$, $b=\mf b/2$.

Since $B(\rho) = D(1-\rho)$, $D$ is concave if and only if $B$ is
concave. The functions $B$ is concave if $3a_1 +a_0\le 4a_2 \le 4a_0$.
We claim that these inequalities are in force. On the one hand, as $\mf
b > \mf a>0$ and $a_2>0$, we have that $a_2 < 2 a_2 + 4\mf b - \mf a =
a_0$. On the other hand, $3a_1 +a_0 -4a_2 = 2a_1 + 4\mf b - 2a_2 =
2\mf a + 4\mf b - 2a_2 \le 0$ from the definition of $a_2$.

This shows that all the assumptions of Theorem \ref{sldp} are
fulfilled.

\smallskip\noindent{\bf Acknowledgements}. The authors wish to express their
gratitude to the referee for a very careful reading which helped to
improve the presentation. C. Landim has been partially supported by
FAPERJ CNE E-26/201.207/2014, by CNPq Bolsa de Produtividade em
Pesquisa PQ 303538/2014-7, and by ANR-15-CE40-0020-01 LSD of the
French National Research Agency.
K. Tsunoda has been partially supported by Grant-in-Aid for Research Activity
Start-up JP16H07041.

\end{document}